\documentclass[lettersize,journal]{IEEEtran}
\usepackage{amsmath,amsfonts}
\usepackage{amsthm,bm}
\usepackage{amssymb}
\usepackage{algorithmic}
\usepackage{algorithm}
\usepackage{array}
\usepackage[caption=false,font=normalsize,labelfont=sf,textfont=sf]{subfig}
\usepackage{textcomp}
\usepackage{stfloats}
\usepackage{url}
\usepackage{verbatim}
\usepackage{graphics}
\usepackage{graphicx}
\usepackage{cite}
\usepackage[colorlinks,urlcolor=blue,citecolor=blue,linkcolor=blue]{hyperref}
\usepackage{color}
\usepackage[braket, qm]{qcircuit}
\usepackage{pgfplots}
\usepackage{xr}
\usepackage{xcolor}
\usepackage{caption}
\usepackage{epstopdf}
\usepackage{float}

\newtheorem{lemma}{Lemma}[section]
\newtheorem{definition}{Definition}[section]

\DeclareMathOperator*{\argmin}{argmin}
\DeclareFontFamily{U}{wncy}{}
\DeclareFontShape{U}{wncy}{m}{n}{<->wncyr10}{}
\DeclareSymbolFont{mcy}{U}{wncy}{m}{n}
\DeclareMathSymbol{\Sh}{\mathord}{mcy}{"58}

\usepackage{tikz}
\usetikzlibrary{positioning,shapes,snakes,arrows,backgrounds,external,fit,calc,matrix}
\usepackage{lmodern}
\tikzset{
    BlockCPU/.style={draw,thick, fill=blue!20, rectangle},
    BlockAltre/.style={draw,thick, fill=blue!35, rectangle},
    Periferic/.style={ellipse, draw, fill=blue!15},
    Registre/.style={rectangle, draw, fill=blue!5},
    Nom/.style={font=\normalsize\sffamily,text centered, minimum size=1cm,textwidth=1.5cm}
}

\begin{document}

\title{Multi-Function Multi-Way Analog 
Technology for Sustainable Machine Intelligence Computation}

\author{Vassilis Kalantzis,
Mark S.\ Squillante, % \IEEEmembership{Fellow, IEEE},
Shashanka Ubaru, 
Tayfun Gokmen,
Chai Wah Wu, %  \IEEEmembership{Fellow, IEEE},
Anshul Gupta,
Haim Avron,
Tomasz Nowicki,
Malte Rasch,
Murat Onen,
Vanessa Lopez-Marrero,
Effendi Leobandung,
Yasuteru Kohda,
Wilfried Haensch,
Lior Horesh
\thanks{V.\ Kalantzis, M.\ Squillante, S.\ Ubaru, T.\ Gokmen, C.W.\ Wu, A.\ Gupta, T.\ Nowicki, M.\ Rasch, L.\ Horesh are with IBM Research, Thomas J.\ Watson Research Center, Yorktown Heights, NY, USA}
\thanks{H.\ Avron is with Tel Aviv University, Tel Aviv, Israel}
\thanks{M.\ Onen is with Massachusetts Institute of Technology, Cambridge, MA, USA}
\thanks{V.\ Lopez-Marrero is with Brookhaven National Laboratory, Upton, NY, USA}
\thanks{E.\ Leobandung and Y.\ Kohda are with IBM Research, Tokyo Research Laboratory, Tokyo, Japan}
\thanks{Correspondence to: vkal@ibm.com, mss@us.ibm.com, lhoresh@us.ibm.com}
}

%\markboth{IEEE Transactions on Pattern Analysis and Machine Intelligence}
%{Kalantzis, Squillante,
%\MakeLowercase{\textit{(et al.)}:
%Multi-Function Multi-Way Analog Technology for Sustainable AI Computation}
%}

% % column for its text to clear the IEEEpubid mark.

\maketitle

\begin{abstract}
Numerical computation is essential to many areas of artificial intelligence (AI), whose computing demands continue to grow dramatically, yet their continued scaling is jeopardized by the slowdown in Moore's law.
\textit{Multi-function multi-way analog} (MFMWA) technology, a
% promising
computing architecture comprising arrays of memristors
% that support
supporting
in-memory computation of matrix operations, can offer tremendous
% advances with respect to computation and energy improvements,
improvements in computation and energy,
but at the expense of inherent unpredictability and noise.
We devise novel randomized algorithms tailored to MFMWA architectures that mitigate the detrimental impact of imperfect analog computations while realizing their potential benefits across
various
areas of
% AI and beyond.
AI,
such as applications in computer vision.
Through
% a combination of
%mathematical
analysis, measurements from analog devices,
% marking the first use of MFMWA operations in such an algorithm,
and simulations of larger systems, we demonstrate orders of magnitude reduction in both computation and energy with accuracy similar to digital computers.
% Together, our results can help to unlock the true potential of sustainable numerical computation and its impact on various areas of AI and their applications.
\end{abstract}

%%%%%%%%%%%%%%%%%%%%%%%%%%%%%%%%%%%%%%

\section{Introduction}
%
%Sustainable computing is critical to our society and our world, with motivating examples abound.
The environmental footprint of artificial intelligence (AI) computation today is becoming prohibitive, further exacerbating the great need for computational sustainability.
With numerical computation essential to many areas of AI, these computational trends continue to increase at a startling pace.
% grow and compound the desperate need for computational sustainability.
Such numerical computations can consume energy upwards of the hourly output of one or more nuclear plants, with current trends on the energy requirements of all numerical computations expected to exceed global energy production levels by the year $2040$~\cite{Sustainability2040,online}.
At the same time, 
% Meanwhile,
from the earliest computing facilities in the mid-1940s to the computing facilities of today, numerical computation has been central to our knowledge and understanding of complex phenomena and an integral part of successfully solving numerous grand challenges at the forefront of many
% disciplines---such as cosmology, genetics, climate science, computational chemistry, natural language processing and drug discovery,
disciplines---including AI, machine intelligence (MI), and pattern analysis in computer vision and imaging, 
to name just a few.

For decades, the need to solve larger computational problems at increasingly faster rates has been primarily met by shrinking the size and increasing the number of transistors in complementary metal-oxide-semiconductor (CMOS) integrated circuits.
This approach, unfortunately, can no longer support the exponential growth requirements of modern MI and AI applications for computational speed and energy consumption, far outpacing the rate of improvements in digital hardware~\cite{https://doi.org/10.1002/aisy.202000085}, especially in light of a slowing Moore's law~\cite{mooreslaw}.
In addition, the physical separation of processing units and main memory in most digital architectures (e.g., von Neumann) leads to wasted computing cycles and longer execution times due to memory latencies and data communication, as well as greater energy consumption~\cite{dongarra2017design}.
It is therefore imperative that we explore alternative architectures to enable numerical computations that handle the ever-increasing needs of
%modern 
MI and AI
applications, and do so in a manner that supports computational sustainability.

Meanwhile, rapid technological advances based on analog devices have made it feasible to significantly improve the computational speed and energy consumption of certain numerical
operations~\cite{yang2013memristive, agarwal2016energy}. 
Multi-Function Multi-Way Analog (MFMWA) is a promising novel computing architecture consisting of arrays of \emph{memristors}~\cite{chua:1971},
one advantage of which is the natural combination of computation and data storage that allows in-memory computation, thus eliminating the memory-hierarchy bottlenecks of digital computers.
Another major advantage is that MFMWA arrays can perform multiple operations, such as matrix-vector (MV) products and rank-$1$ or outer-product (OP) updates (whose computation depends on matrix size), in time that is independent of the matrix size (subject to physical constraints and architectural design), thus eliminating the critical computational and energy bottlenecks of digital computers for these operations.
% in combination with novel MFMWA-based randomized algorithms.
% %tailored to MFMWA architectures.
MFMWA hardware devices are becoming more available; see, e.g.,~\cite{ambrogio:analog:2023}.

However, the
% The 
promise of
%revolutionary 
tremendous
% computational sustainability through the 
numerical speedups and energy reductions of MFMWA technology,
and computational sustainability in turn, 
comes with a unique set of challenges.
Analog-based computing incurs
%random
noise in the results of the
computations due to hardware imperfections and temporal conductance variations.
%when updating the synaptic weights. \Lior{do we want to use terms from NN such as synaptic weights ? }
% As a result of these drawbacks,
As a result,
the full potential of MFMWA technology can only be realized through the design of appropriate algorithms tailored to MFMWA architectures that mitigate and tolerate the impact of such inherent noise.
% However, recent 
Recent
applications of analog technology,
however,
have been largely restricted to deep
neural network (NN)
training and inference where stochastic learning methods are thought to be resilient to analog variations.
Although a few other specific applications have been
%suggested~\cite{burr2017neuromorphic,sebastian2020memory},
suggested~\cite{sebastian2020memory},
%such as linear system solvers and random number generation,
% they have not exploited the full sustainable computation potential of MFMWA technology in our view.  
they have neither exploited the full computational potential of MFMWA technology nor 
considered applications across areas of MI and AI.

In this paper,
we
%make the critically important
%observation 
show
that it is possible to 
fulfill the great promise of
% %performance 
% computational
% and energy benefits 
numerical speedups, energy reductions and
computational sustainability
by
mitigating and tolerating the errors from analog hardware imperfections through
% hybrid digital-analog architectures and
% %%major advances in the core dimensions
% %a new class of
new
low-complexity algorithms
%and computational models
tailored to
% the capabilities of analog computing technology.
MFMWA in-memory architectures.
Specifically,
we design, analyze and implement a general class of \emph{randomized numerical linear algebra algorithms} for
%certain
core numerical methods involving
%methods that involve
%outer-product
OP
updates and 
%matrix-vector
MV
products that run on
%Multi-Function Multi-Way Analog
MFMWA
%memristor crossbar arrays
technology
and achieve solution
%%the same
%comparable
accuracy
%as
comparable to
digital computers.
%By addressing and mitigating in this way the negative impact of analog hardware imperfections in the results of the numerical computations,
Our solution approach
% based on this combination of architecture and algorithms 
resolves the negative impact of imperfect numerical computations on analog hardware and
%fulfills the promise of analog computing; this in turn
enables us to overcome the critical 
%performance 
computational and energy bottlenecks of modern digital
%computing.
architectures.
This renders tremendous computational speedups that support the acceleration and scalability of
% at scale and accelerate computational applications across
% %scientific, ML and other
% many
% disciplines.
numerical applications and tremendous energy reductions that support sustainability and low-power devices,
particularly for certain MI applications.
%, and to
%support the entire spectrum of computational environments from edge computing to cloud
%computing
%and high-performance computing.
% %Figure~\ref{fig:compscience} illustrates our perspective on the universe of analog-aware randomized algorithms, computing architectures, and computational applications.
% %across scientific, ML and other disciplines.
% %
% Figure~\ref{fig:cartoon1} illustrates such connections via the inherent dependence of several important tasks in computational science and ML on core numerical methods involving
% %outer-product
% OP
% updates and
% %matrix-vector
% MV
% products. The tasks in the leftmost and rightmost regions rely primarily on
% %outer product
% OP
% and
% %(matrix-vector)
% MV
% operations, respectively, and those in the middle region rely on both of these numerical kernels.
%
% Our focus herein is primarily on computational tasks that involve either OP operations or a combination of OP and MV operations, particularly those involving a majority of OP updates, which remain mostly unexplored.
Since matrix computation is the foundation of many important AI and MI methods, we focus on core matrix problems that involve ordinary linear least squares and principal component analysis with either OP operations or a combination of MV operations and (a majority of) OP operations.

% We present detailed simulation experiments to evaluate our new 
% analog-assisted randomized algorithms on artificial and real-world datasets from population genetics and video processing. These simulations demonstrate
% major improvements in computation time and energy consumption, including 
% speedups by more than factors of $20$ to $40$ and energy 
% reductions by more than a factor of $10$ in comparison to state-of-the-art digital
% computers. 
% We also present measurements from experiments conducted on a real analog 
% crossbar array device, marking the first time an 
% %outer-product
% OP
% update and randomized algorithm have been implemented on actual analog hardware.
% These physical experiments further demonstrate that our approach, in comparison with the corresponding digital approach,  achieves results of comparable quality together with computation and energy improvements consistent with our simulation experiments.

% {\bf In order to consider much larger problem instances and allow us to show results for analog devices on our planned roadmap. Supplement physical experiments with simulations}

%{\bf + Reorganize the red paragraph above as follows:}

To demonstrate and quantify these significant improvements, we derive
% a mathematical 
an
analysis of our randomized algorithms tailored to MFMWA architectures together with various experimental results. 
This includes measurements from experiments conducted on a real 
%Multi-Function Multi-Way Analog 
MFMWA
%array 
device with a randomized streaming algorithm, marking the first time an algorithm with
%outer-product
OP
update as the dominant computation has been implemented on actual physical analog hardware.
In order to consider much larger problem instances and
%to allow us 
to show results for MFMWA devices on
% our planned device roadmap, 
planned device roadmaps, 
we  supplement these physical experiments with simulation experiments of our analog-assisted randomized algorithms on synthetic and real-world datasets,
% which includes 
such as
the computer vision task of background subtraction in surveillance videos.
%from population genetics  and video processing applications.
This collection of analytical and experimental results demonstrate major advances in 
% computation time and sustainability, 
numerical speedups, energy reductions and
computational sustainability,
including
% speedups by more than factors of $20$ to $40$ and energy reductions by more than a factor of $10$ in comparison 
orders of magnitude improvements
% in computational speedups and energy reductions compared with
over state-of-the-art algorithms on digital computers.
Together, our algorithms and results can help to unlock the true potential of sustainable numerical computation and its impact on various 
% areas of AI and their applications.
application areas of MI and AI.

The remainder of the paper is organized as follows.
After providing some technical preliminaries in Section~\ref{sec:prelim}, we present our randomized numerical linear algebra algorithms in Section~\ref{sec:alg}.
We evaluate the performance of our algorithms across various applications in Section~\ref{sec:results},
followed by concluding remarks in Section~\ref{sec:concl}.
The supplement contains additional details and results.

\section{Technical Preliminaries}
\label{sec:prelim}
Before presenting our randomized algorithms tailored to MFMWA architectures, we provide some background information and technical details on analog computing technology and randomized numerical linear algebra.
Additional details on both are provided in the supplement.

% %\section*{Background}
% \section*{Background: Analog Computing Technology}
% % %
% % We present some background on numerical computation across various disciplines and related analog computing technology.
% % % , with additional details provided in
% % % %\ref{app:ssec:numerical}, \ref{app:ssec:analog} and \ref{app:ssec:sim}.
% % % the supplement.
% % %of the supplement.

\begin{figure*}[htbp]
\centering
% \includegraphics[width=0.45\textwidth,trim={20cm 0cm 0cm 20cm}]{mat-vec.png}\\
% \includegraphics[width=0.45\textwidth,trim={3cm 0cm 3cm 0cm}]{outer-product.png}
% \qquad
% \includegraphics[width=0.45\textwidth,trim={6cm 0cm 6cm 0cm}]{pulse-conicidence-color-rot.png}
\includegraphics[width=0.7\textwidth,trim=0cm 0cm 0cm 0cm,clip]{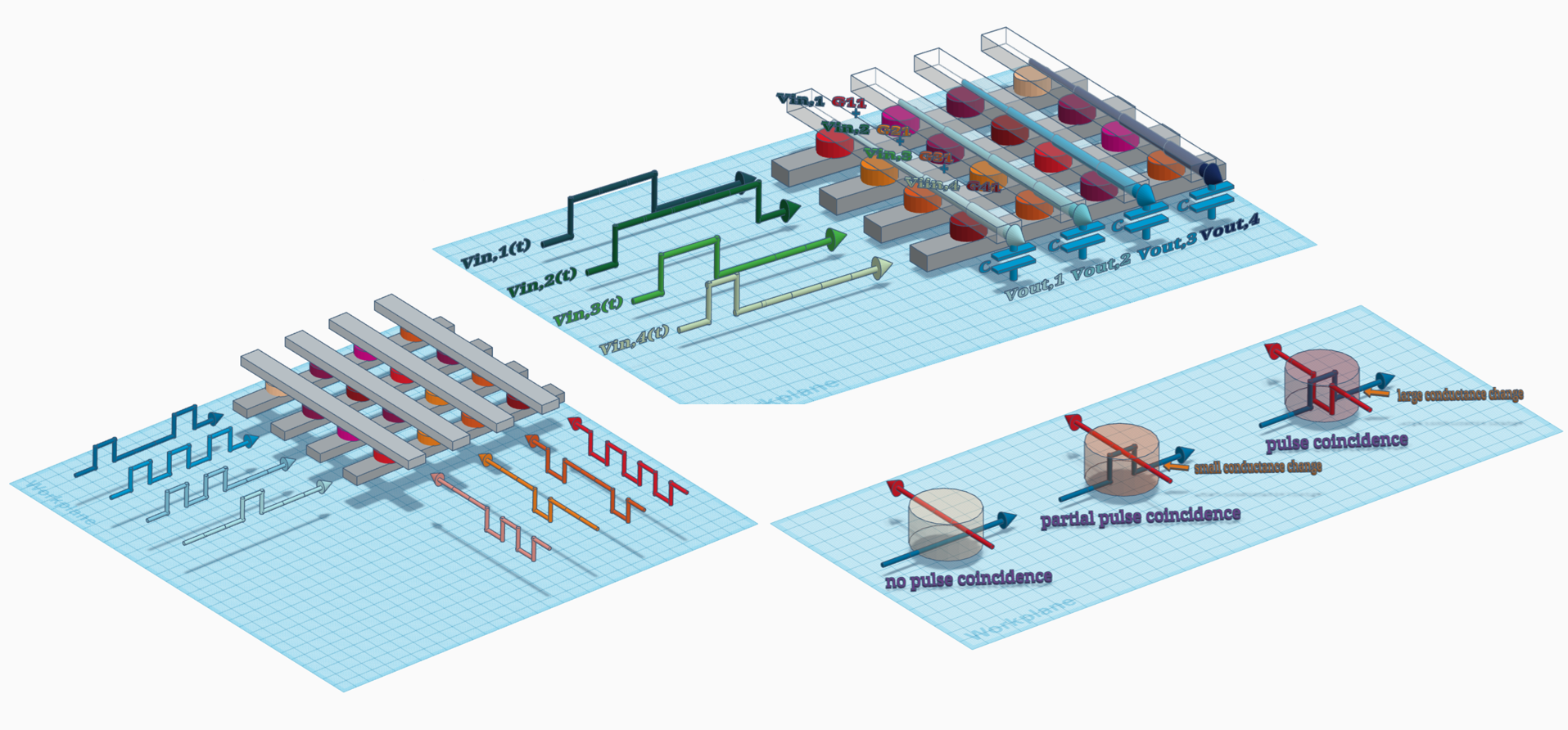}\\
(a)\\
\includegraphics[width=0.55\textwidth]{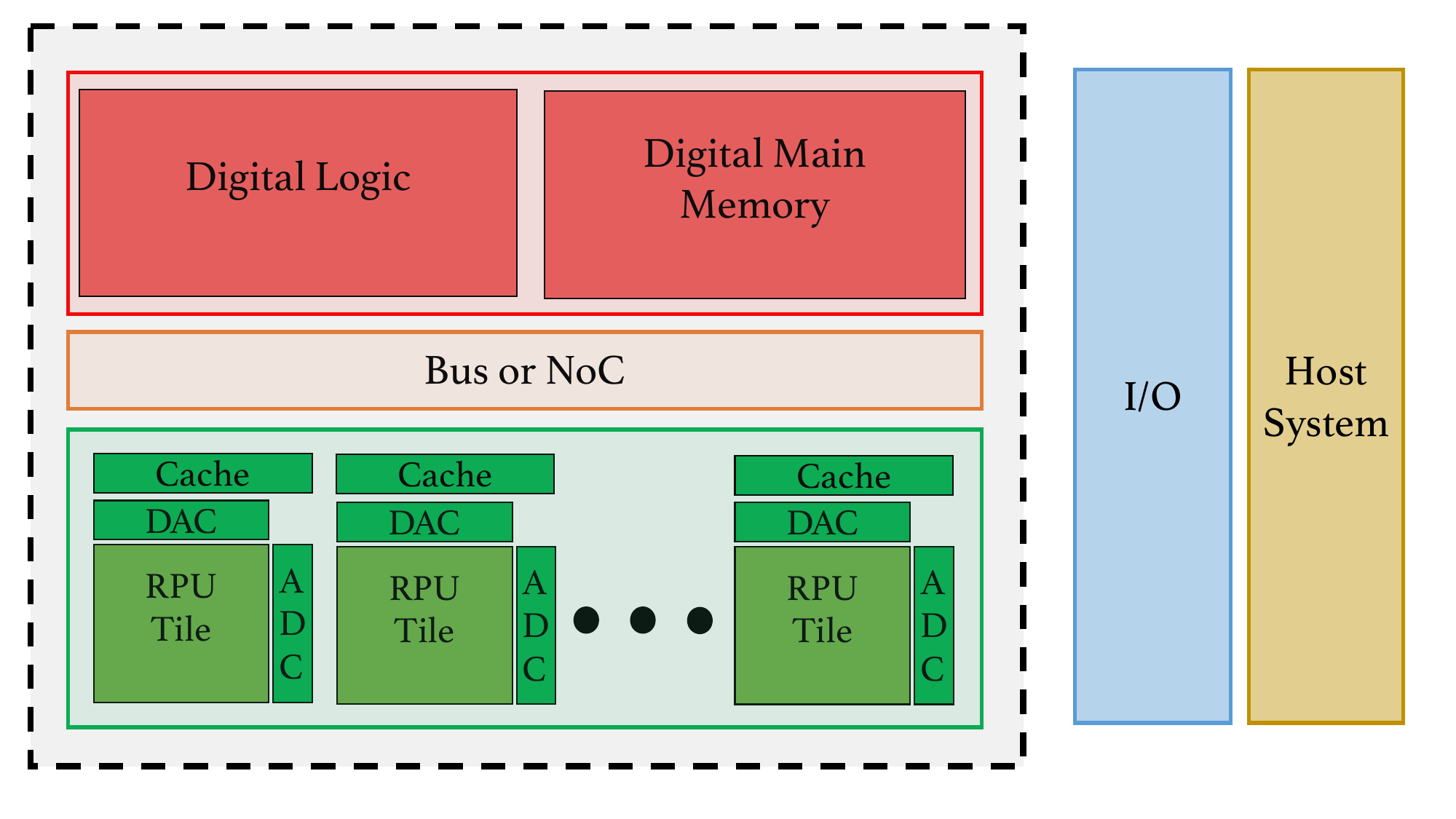}\\
(b)
\caption{\textit{(a) Multi-Function Multi-Way Analog array technology. Center: MV product. The output voltages $V_{out,k}$ consist of the integral of the currents $I_k$ over time divided by the capacitance $C$ per $V_{out,k} = \frac{1}{C} \int_0^{\top} I_k(\tau) d \tau$. At each cross point, an  element with conductance $g_{kl}$ interacts with the input signal $V_{in,j}(t)$ per Ohm's law and determines the current signal $I_k(t)$.
Left: OP update alters the conductance of each cross point based on coincidence of the input pulse signals from both sides. Right: zoom-in view over the pulse coincidence update mechanism. The conductance changes whenever the coincidence of the input pulses occurs.  \\
(b) Hybrid architecture consisting of a digital chip including a (multi-core) CPU and 
RAM system memory, connected to one or more 
%Multi-Function Multi-Way Analog 
MFMWA
arrays through a dedicated system bus. The CPU is focused on executing low-complexity 
operations as well as transferring data to the
%Multi-Function Multi-Way Analog 
MFMWA
%crossbar 
array(s) through the dedicated 
system bus.
}
}
\label{revised-figna1}
\end{figure*}

\subsection{Analog Computing Technology}
\label{sec:analog}
The general concept of arrays of analog devices as the basis of computing machines dates back to before the 1960s~\cite{analogjackson60}. 
%The rise of the digital computer has diminished the competitiveness of analog computing technology and only in recent years is there a resurgence in using analog technology for low-power high speed computations.
Once a matrix is initialized on an
%Multi-Function Multi-Way Analog 
MFMWA
array, MV and OP operations can be executed orders of magnitude faster on analog hardware than on hardware based on digital memory~\cite{agarwal2016energy}. 
Moreover, the energy efficiency of 
analog
devices can scale significantly better than that of digital hardware for a variety of tasks ranging from training a
%neural network 
NN
($\sim 250\times$) \cite{ambrogio2018equivalent} to approximating the solution of a linear system ($\sim 10\times$) \cite{feinberg2018enabling}. 
The starting point of MFMWA arrays for numerical computation considers each memory element as an electrically programmable resistor. 
% Many different materials have been explored as instantiations of such memory devices---e.g., charge-based Static Random Access Memory (SRAM), or memristive Phase Change Material (PCM) and Resistive Random Access Memory (RRAM) \cite{sebastian2020memory}---each having different tradeoffs in terms of noise, volatility, switching characteristics, and energy usage. 
% Independent of the particular choice of memory device, however,
Independent of the materials used to implement such memory devices, Kirchhoff's laws and Ohm's law are leveraged to obtain massively parallel analog multipliers by considering the voltages and currents of the circuit as informational data.
This analog computation can be well integrated into the otherwise traditional digital compute flow, including parallel analog-to-digital converters (ADCs) and digital-to-analog converters (DACs) 
to convert analog values to digital and vice versa.
%of output charges back to digital. 

Such
%analog 
MFMWA
technology enables the execution of certain matrix operations
% The benefits of Multi-Function Multi-Way Analog crossbar arrays come from the fact that certain matrix operations can be performed
in a remarkably efficient manner.
% This is illustrated in Figure~\ref{revised-figna1}.  
Specifically, 
as illustrated in Figure~\ref{revised-figna1}(a),
the basic
% mathematical 
numerical
operations that can be performed on these
MFMWA
arrays are: $(i)$
%matrix-vector (MV)
MV
%multiplication
product by directly utilizing multiplication and addition via Ohm's law and Kirchhoff's laws, respectively;
and $(ii)$ 
%outer-product (OP)
OP
update by stochastic computing and voltage pulse coincidences.
%(see below). 
For
%matrix-vector (MV)
an MV
%multiplication
product
with a dense $m\times m$ matrix, roughly $2m^2$ arithmetic operations are required on digital hardware;
%a rank-$1$
an OP
update on such a matrix requires the same number of digital operations. In strong contrast, if these
% mathematical 
numerical
operations are performed on a hybrid architecture as in Figure~\ref{revised-figna1}(b), the MV
%multiplication
product
and
%rank-$1$
OP
update can be executed in constant time independent of the matrix size, subject to physical constraints on the design and size of the MFMWA arrays.

%\section*{Hybrid MFMWA Computing Architectures}
%\label{sec:arch}
%\subsection*{Stochastic Computing and Rounding}
%
%Figure \ref{fig:schematic}
We 
%propose to use
consider the use of multiple tiles of memristive arrays of the type illustrated in
Figure~\ref{revised-figna1}(a) 
% Figure~1(a)
in combination with a large digital main-memory, digital logic, ADCs and digital-to-analog converters (DACs), and digital cache storage to construct a hybrid MFMWA accelerator architecture as illustrated in
% Figure~{\color{red} 6}
% (supplement).
Figure~\ref{revised-figna1}(b).
% Figure~1(b).
%illustrates a hybrid analog accelerator architecture consisting 
%of a large digital main-memory, digital logic, and a collection of analog devices each with multiple 
%tiles of memresistive crossbar arrays, an ADC and DAC, and additional digital storage. 
The digital logic supports low-complexity matrix operations (e.g, dot products and AXPY operations) that do not benefit from MFMWA technology. The hybrid MFMWA accelerator interfaces with a digital host system comprising a conventional von Neumann architecture.% with a central processing unit (CPU) and a system memory. 

% As previously noted, the key ideas behind using MFMWA arrays for numerical computation start with considering each memory element as an electrically programmable resistor. Many
With respect to each memory element as an electrically programmable resistor, many
different materials have been explored as instantiations of such memory devices---e.g., charge-based Static Random Access Memory (SRAM), memristive Phase Change Material (PCM) and Resistive Random Access Memory (RRAM) \cite{sebastian2020memory}---each having different tradeoffs in terms of noise, volatility, switching characteristics, and energy usage.
In each case, Kirchhoff's laws and Ohm's law
are leveraged to obtain massively parallel analog multipliers by considering the voltages and currents of the circuit as informational data.
This analog computation can then be well integrated into the traditional digital compute flow outlined above.
For instance, capacitors can be used to integrate the output current during the fixed time-span of applying width-modulated voltage pulses of the digital inputs, with parallel ADCs converting the output charge back to digital.

To encode information for input to the MFMWA arrays when executing the OP operation, we use stochastic computing where information is represented and processed as binary pulses that allow simple logic gates to perform computations such as addition and multiplication. 
%This coding format is used in analog computing to provide the input data to the analog crossbar array during the OP operation. %However, the decoding of the stochastic computing output is done via the physics of each crossbar
%element \cite{gokmen2016acceleration}.
%element.
For instance, as illustrated in 
Figure~\ref{revised-figna1}(a), 
% Figure~1(a), 
when there is pulse coincidence due to an overlap in the voltage pulses coming from a row and a column, the voltage difference at the intersection exceeds a certain threshold and, as a result, the stored conductance value changes by a small amount. Therefore, a pulse coincidence followed by an incremental conductance change effectively performs the multiplication and summation operation needed for OP updates using stochastic computing and the device physics on MFMWA arrays.
%Stochastic rounding \cite{Gupta2015a},
In addition, we use stochastic rounding---a form of multibit stochastic computing where a stochastic binary sequence is added to the traditional rounded value---to encode the input to the MFMWA arrays when performing the MV operation and in the modeling of the  %analog-to-digital converter
ADC
noise.

% We refer the reader to the supplement for additional details on MFMWA technology and hybrid MFMWA architecture technology.

\subsection{Randomized Numerical Linear Algebra}
%\section*{Randomized Algorithms on MFMWA Architectures}
%
Matrices are among the most common objects used to represent interactions between data,
% in many disciplines, 
making matrix computation the foundation of many important MI and AI methods.
For example, a matrix $A$ can capture information related to a set of objects, each with a set of corresponding features;
or the discretization of Poisson's differential equation in space;
or the encoding of weights between successive layers of NNs.
Hence, the design of efficient
% numerical
algorithms for core matrix problems can accelerate numerous computational applications in MI and AI.
We focus here on two such representative matrix problems.

\vspace*{0.25cm}
% \begin{enumerate}
%     \item 
\noindent
$1$).
    {\bf Ordinary linear least-squares} (OLLS): Given an $m\times n$ matrix $A$ and an 
    $m\times 1$ vector $b$, compute an $n\times 1$ vector $\widehat{x}$ such that 
    $\widehat{x}=\argmin\limits_{x\in\mathbb{R}^n}
    \dfrac{1}{m}(b-Ax)^{\top}(b-Ax)$.\\
    % \item 
$2$).
    {\bf Principal component analysis} (PCA): Given an $m\times n$ matrix $A$, compute 
    a basis of the linear subspace spanned by its $k\in \mathbb{Z}$ dominant left singular vectors.
% \end{enumerate}
\vspace*{0.25cm}

Many efficient numerical algorithms have been designed for the above matrix problems when $A$ is small~\cite{golub2013matrix}. Nevertheless, the ever-increasing sizes of
% modern 
MI and AI
problems create challenges for state-of-the-art algorithms; e.g., the data matrix $A$ may not fit in system memory, the data comprising $A$ is only available over time, and/or repeated access to data from secondary storage incurs tremendous time and energy for input/output communication. 
Specific examples of applications involving OLLS and PCA with tremendous and ever growing numerical computational requirements include
the analysis of
% features and 
patterns in computer imaging and vision~\cite{gittens2016matrix,LiHuGuTi04,candes2011robust},
the analysis of
%statistical 
patterns within privacy-preserving collaborative genome-wide association studies
(GWAS)~\cite{bush2012genome,Uffelmann2021,Mbatchou2021},
% (GWAS)~\cite{Uffelmann2021,Mbatchou2021},
and 
the analysis of population structure in human genetics~\cite{paschou2007pca,bose2019terapca}.
This motivates the need to design next-generation algorithms that: $(i)$ provide an approximate solution when the data becomes available over time; $(ii)$ trade high accuracy for reduced computational complexity; and $(iii)$ execute on in-memory hybrid computing architectures such as Figure~\ref{revised-figna1}(b).
% We exploit randomized numerical linear algebra to effectively address these challenges~\cite{drineas2016randnla}. 
% \Malte{One possible criticism we might get here is whether these two problems are actually the main drivers of the recent compute increase in the world mentioned at the beginning. So does it even matter to make these more computational efficient? But I guess we do not have any data of what the dominating type of computation is when looking of the world energy budget. Maybe bitcoin mining is a candidate...  } 
% \mss{what about simply mentioning a few areas of computational science that are well known to require tremendous computation and make great use of OLLS and/or PCA?}
To
% effectively
address these challenges of large-scale numerical computations,
we exploit the important and powerful multi-disciplinary field of randomized numerical linear algebra~\cite{drineas2016randnla}.

\section{Algorithms and Analysis}
\label{sec:alg}
In this section we first present our framework for randomized numerical linear algebra within the context of MFMWA architectures.
We then present our design and implementation of a general class of numerical linear algebra algorithms for randomized OLLS and PCA tailored to hybrid MFMWA-digital environments, together with a corresponding analysis related to these algorithms.
Additional details and results are provided in the supplement.

\begin{table*}[htb]
\caption{\it Summary of time complexity (in microseconds, or $\mu$s) and energy complexity (in micro-Joules, or $\mu$J) comparisons between the hybrid MFMWA accelerator architecture and purely digital accelerator architecture for the computational core primitives of our algorithms with respect to an MV product of an $(n \sqrt{m}) \times (n \sqrt{m})$ matrix and an $(n \sqrt{m}) \times 1$ vector, and an OP update of an
% order-$(n \sqrt{m})$ 
$\mathcal{O}(n \sqrt{m})$ 
column vector and an 
% order-$(n \sqrt{m})$ 
$\mathcal{O}(n \sqrt{m})$ 
row vector. The MFMWA accelerator %architecture 
%%consists of
includes
$64$ tiles of $2K \times 2K$ of memresistive arrays, additional digital storage,
%on the order of $n\times m$,
ADC and DAC logic, and additional digital logic to support the MV product and OP update operations;
the matrix write primitive is from the digital main-memory to the arrays, and the matrix read primitive is from the arrays to the digital main-memory. The digital accelerator
%architecture
is based on an A100 GPU with $40$-$80$MB digital cache-memory.
% ;
% the matrix write primitive is from the digital main-memory to the digital cache-memory, and the matrix read primitive is from the digital cache-memory to the digital main-memory.
}\label{tab:analog-digital}
%{\small
\begin{center} 
\begin{tabular}{||l|l|l||}
\hline\hline
\textbf{Core Primitive} & \textbf{Time Complexity} & \textbf{Energy Complexity} \\
% \hline\hline
% MFMWA Matrix Write &
% $T_{\mathcal{MW}} \in [2048\mu\mbox{s},20480\mu\mbox{s}]$ 
% &
% $E_{\mathcal{MW}} \in [262144\mu\mbox{J},13107200\mu\mbox{J}]$
% \\
% \hline
% Digital Matrix Write &
% $T_{\mathcal{MW}}^d \approx 250\mu\mbox{s}$ 
% &
% $E_{\mathcal{MW}}^d \approx 12000\mu\mbox{J}$
% \\
\hline\hline
MFMWA MV Product &
$T_{\mathcal{MV}} \in [0.135\mu\mbox{s},0.240\mu\mbox{s}]$
&
$E_{\mathcal{MV}} \in [12.928\mu\mbox{J},33.28\mu\mbox{J}]$
\\
\hline
Digital MV Product &
$T_{\mathcal{MV}}^d \in [250.005\mu\mbox{s},250.02\mu\mbox{s}]$
&
$E_{\mathcal{MV}}^d \in [12000.064\mu\mbox{J},12000.64\mu\mbox{J}]$
\\
\hline\hline
MFMWA OP Update &
$T_{\mathcal{OP}} \in [0.11\mu\mbox{s},0.14\mu\mbox{s}]$
&
$E_{\mathcal{OP}} \in [12.928\mu\mbox{J},33.28\mu\mbox{J}]$
\\
\hline
Digital OP Update &
$T_{\mathcal{OP}}^d \in [250.01\mu\mbox{s},250.04\mu\mbox{s}]$
&
$E_{\mathcal{OP}}^d \in [12000.128\mu\mbox{J},12001.28\mu\mbox{J}]$
% \\
% \hline\hline
% MFMWA Vector Read &
% $T_{\mathcal{RV}} \in [0.005\mu\mbox{s},0.020\mu\mbox{s}]$
% &
% $E_{\mathcal{RV}} \in [64000\mu\mbox{J},640000\mu\mbox{J}]$
% \\
% \hline
% Digital Vector Read &
% $T_{\mathcal{RV}}^d \approx 0.0153\mu\mbox{s}$
% &
% $E_{\mathcal{RV}}^d \approx 0.7324\mu\mbox{J}$
\\
% \hline\hline
% MFMWA Matrix Read &
% $T_{\mathcal{RM}} \in [2293.76\mu\mbox{s},4259.84\mu\mbox{s}]$
% &
% $E_{\mathcal{RM}} \in [212860.928\mu\mbox{J},555745.28\mu\mbox{J}]$
% \\
% \hline
% Digital Matrix Read &
% $T_{\mathcal{RM}}^d \approx 250\mu\mbox{s}$
% &
% $E_{\mathcal{RM}}^d \approx 12000\mu\mbox{J}$
% \\
\hline\hline
\end{tabular}
\end{center} 
\end{table*}

%%\section*{Key Algorithms}
%%\subsection*{Randomized Algorithms on Analog Crossbar Arrays}
\subsection{Randomized Algorithms on MFMWA Architectures}
%\label{sec:algs}
%
One such well-known example of
%randomization 
randomized numerical linear algebra
concerns \emph{matrix sketching} algorithms that have been successfully applied to the 
% cryptography of
analysis of statistical patterns in 
privacy-preserving collaborative
GWAS~\cite{bush2012genome,Uffelmann2021,Mbatchou2021},
% genome-wide association studies (GWAS),
improving the runtime of the most advanced cryptographic protocols by two orders of magnitude~\cite{kockan2020sketching}.
Similarly, randomization is a key component of algorithms that lead to dramatic reductions in the runtime of
%\emph{principal component analysis} (PCA)
PCA
applied to 
% population genetics~\cite{bose2019terapca}
the analysis of population structure in human genetics~\cite{paschou2007pca,bose2019terapca}
as well as to 
% cosmology, climate science, and mass spectrometry imaging~\cite{gittens2016matrix}.
the analysis of features and patterns in mass spectrometry imaging and climate science~\cite{gittens2016matrix} and in computer vision~\cite{LiHuGuTi04,candes2011robust}.

By exploiting the above hybrid MFMWA architecture to offload some computation on very large data matrices or on data that becomes available one chunk at a time,
% these chunks to analog computers, 
we can execute high performance computations involving matrix sketching and PCA while reaping the power/speed/area benefits of MFMWA arrays.
The occurrence of inherent stochasticity in these arrays necessitates a delicate interplay between the analog and digital components of
% our 
the
hybrid MFMWA architecture in order to achieve the desired accuracy and performance metrics of the computational task.
% We therefore design a new class of randomized algorithms 
% that are tailored to analog crossbar arrays in order to break the existing performance and 
% energy bottlenecks in computational applications.
%% While our focus 
%% does not extend beyond the interest of randomized algorithms, we note that 
%% the topic of leveraging analog crossbar arrays to accelerate deterministic 
%% numerical algorithms is an active one as well, e.g., see %%\cite{sun2019solving,le2018mixed,liu2018memristor,huang2016evaluation}. 
To understand this delicate interplay, we  derived 
% a mathematical 
an
analysis of 
the communication costs between analog
and digital components (refer to the supplement) as well as 
% a mathematical 
an
analysis of the relationship between noise and performance when an approximate inverse preconditioner on an MFMWA array is used to solve sparse linear systems, deriving bounds on the stochastic analog noise to ensure convergence of the solver~\cite{Kalantzis2021-nb}.
In particular, for the moderately conditioned practical
problems that we consider,
reductions in computation and energy significantly outweigh the additional number of iterations required.

A key result that spurred research on randomized matrix sketching is the famous Johnson-Lindenstrauss Lemma, which is expressed as follows. 
\begin{lemma}[Johnson-Lindenstrauss Lemma]
	Given $0 <\epsilon < 1/2$, a set ${\cal A}$ of $m$ points in $\mathbb{R}^n$ and a 
	parameter $\ell \geq 8\ln(m)/\epsilon^2$, there exists a map $f: \mathbb{R}^n\rightarrow 
	\mathbb{R}^\ell$ such that
	\begin{equation*}
	(1-\epsilon) \|u-v\|_2^2 \leq \|f(u)-f(v)\|_2^2 \leq (1+\epsilon)\|u-v\|_2^2,
	\end{equation*}
	for all $u,v\in {\cal A}$.
	\label{app:lem4.1}
\end{lemma}
The above lemma states that any $n$-point subset of the Euclidean space can be embedded in $\ell < m$ dimensions without distorting the distances between any pair of points by more than a factor of $1\pm \epsilon$, for any $\epsilon \in (0,1/2)$. Additionally, it provides an upper bound on the dimensionality $\ell$ required for a mapping to ensure that it will approximately preserve distances. Though Lemma~\ref{app:lem4.1} is stated in an existential manner, there are known randomized constructions for $f$ that work for a given set of points with high probability. For example, the map function $f$ can be implemented  through an MV product of the form $f(u)=Su$, where $S\in \mathbb{R}^{\ell \times m}$ is a random ``sketching'' matrix chosen from an appropriate distribution. 

In randomized matrix algorithms, we actually do not use embeddings for a {\em finite} set of points (as provided by the Johnson-Lindenstrauss Lemma~\ref{app:lem4.1}), but rather we use embeddings for a full subspace of points. Furthermore, such algorithms use embeddings that are linear transformations defined by some matrix $S$.
\begin{definition}[$\epsilon$-subspace embedding ($L_2$ norm)]
	For a given $m\times n$ matrix $A$, an $\ell \times m$ random matrix $S$ is an $\epsilon$-subspace for $A$ if, for all vectors $x\in \mathbb{R}^n$,
	\begin{equation*}
	(1-\epsilon) \|Ax\|^2 \leq \|SAx\|^2 \leq (1+\epsilon)\|Ax\|^2.
	\end{equation*}
	\label{app:lem4.2}
\end{definition}
The above definition tells us that, for the appropriate matrix $S$, we can compute a column-wise 
sketch of matrix $A$ such that the squared $L_2$ norm of elements of its column space 
are within $\epsilon$ distance of the corresponding elements of the column space of  its sketch $SA$.
An example of a randomized construction of $S$ is the {\em Gaussian Sketch}, which is essentially a randomized matrix with independent and identically distributed (i.i.d.) Gaussian entries.
\begin{lemma}
Let $S \in \mathbb{R}^{\ell \times m}$ be a random matrix where each entry is sampled i.i.d.\ from ${\mathcal N}(0, \ell^{-1})$. Given $\epsilon\in(0,1)$ and $\delta\in(0,1)$, if $\ell = \Omega(n \epsilon^{-2}\log(1/\delta))$, then $S$ is an $\epsilon$-subspace embedding for a given $A\in\mathbb{R}^{m \times n}$ with probability at least $1-\delta$.
\end{lemma}
Note that there are other constructions of $L_2$-subspace embedding matrices~\cite{woodruff2014sketching}, as well as constructions of subspace embeddings for other norms. 

We exploit here randomization in numerical linear algebra to devise randomized algorithms, tailored to hybrid MFMWA architectures, that achieve the goals of 
% sustainable computation 
numerical speedups, energy reductions and computational sustainability
while mitigating the detrimental impact of imperfect analog computations.
More specifically, we now show that our randomized algorithmic framework, including randomized LLS and randomized PCA, fits well within the context of MFMWA arrays. In particular, our randomized algorithms tailored to hybrid MFMWA architectures can provide orders of magnitude reduction in both computation and energy across many applications, such as those arising in 
% genetics and surveillance, 
the pattern analysis of GWAS and the pattern analysis of computer imaging and vision studies,
with accuracy similar to digital computers.
This includes our
% mathematical 
analysis, derived in the supplement, of the time complexity and energy complexity of core numerical primitives for randomized LLS and randomized PCA under the hybrid MFMWA accelerator architecture of Figure~\ref{revised-figna1}(b) and a corresponding purely digital accelerator architecture.
Table~\ref{tab:analog-digital} summarizes these results, providing a quantitative time and energy complexity comparison for such computational core primitives across both architectures.

\subsection{Randomized
%ordinary linear least-squares 
OLLS
solver}
Let an $m \times n$ matrix $A$ represent $m$ linear equations on $n$ variables, where we assume without loss of generality that $n \ll m$ and $b\in \mathbb{R}^m$ represents a vector of $m$ observed values. 
When $A$ has full rank, the solution of the 
%linear least-squares 
LLS
problem is given by the vector $\widehat{x}=(A^{\top}A)^{-1}A^{\top}b$,
%. This approach, also known as ``normal equations'',
which requires
$\mathcal{O}(mn^2)$ 
% on the order of $mn^2$ 
%floating-point operations (FLOPs). 
arithmetic operations.
For large values of $m$ and $n$, the latter runtime complexity 
can be prohibitive. 
Randomized matrix sketching instead approximates the solution $\widehat{x}$ by that of the sketched problem
$$
    \widehat{x}_S=\argmin\limits_{x\in\mathbb{R}^n}
    \dfrac{1}{\ell}(Sb-SAx)^{\top}(Sb-SAx),
$$
where $S\in \mathbb{R}^{\ell \times m}$. For appropriate matrices $S$, the approximate solution $\widehat{x}_S$ is known to be within
%an $L_2$-norm $\epsilon$ distance of
a distance of $\epsilon \in \mathbb{R}_{>0}$ from the exact solution $\widehat{x}$, where $\epsilon$ 
is an accuracy parameter on which $\ell$ generally depends \cite{woodruff2014sketching}. 
While randomized sketching reduces the size of the problem, choosing the sketching matrix $S$ to be dense and unstructured requires 
$\mathcal{O}(mn\ell)$ 
% on the order of $mn\ell$
% % floating-point operations per second (FLOPS)
% FLOPs
arithmetic operations
to compute the sketch $SA$. When $\ell >n$, this exceeds the
$\mathcal{O}(mn^2)$ 
% order of $mn^2$ 
runtime complexity required to directly compute $\widehat{x}=(A^{\top}A)^{-1}A^{\top}b$. Practical instances of randomized sketching can choose $S$ in a manner that computes a sketch of $A$ with a runtime complexity
of $\mathcal{O}(m\ell \mathtt{log}\ell) + \mathtt{poly}(\ell/\epsilon)$~\cite{sarlos2006improved,ubaru2017low}
% on the order of $m\ell \mathtt{log}\ell + \mathtt{poly}(\ell/\epsilon)$~\cite{sarlos2006improved,ubaru2017low}
or
$\mathcal{O}(\mathtt{nnz(A)})+\mathtt{poly}(\ell/\epsilon)$~\cite{avron2010blendenpik},
% $\mathtt{nnz(A)}+\mathtt{poly}(\ell/\epsilon)$~\cite{avron2010blendenpik},
depending on the approach taken;
here ${\tt nnz}(\cdot)$ denotes the number of nonzero entries and $\mathtt{poly}(\cdot)$ indicates a 
polynomial value of the input. 

Our general algorithmic framework tailored to MFMWA architectures enables extremely efficient randomized sketching of datasets and overcomes existing computational and energy bottlenecks. Since the data matrix $A$ becomes available only partially over time in many real-world applications, we focus on the streaming case. 
The $\ell \times n$ matrix $Z=SA$, referred to as a (linear) sketch of $A$, can then be formed using a series of $m$ OP updates. 
More specifically, we can write
$$Z=SA=\sum_{j=1}^{j=m}S_{:,j}A_{j,:} , $$
where $S_{:,j}\in \mathbb{R}^{\ell}$ and $A_{j,:}^\top\in \mathbb{R}^n$ denote the $j$th column of matrix $S$ and $A^\top$, respectively. 
The $m$ OP updates can be computed using an $\ell\times n$
%Multi-Function Multi-Way Analog 
MFMWA
array in which the conductance of each crossbar is initialized to zero. 
Since the runtime complexity of each OP update is fixed regardless of the values of $\ell$ and $n$ (within the physical constraints and architectural design), the total runtime complexity to obtain the sketch of matrix $A$ is 
$\mathcal{O}(m)$, 
% on the order of $m$, 
thus allowing the use of a larger sketch size to reduce the error in the approximation of $\widehat{x}$ by $\widehat{x}_S$.

From a digital computer perspective, the computational complexity of computing an 
$\ell$-dimensional Gaussian sketch of an $m \times n$ matrix $A$ comprises $\ell n (2m-1)$ 
arithmetic operations. On the other hand, the $m$ OP updates required to sketch $A$ can be performed much more efficiently on MFMWA devices, as noted above. Following 
Table~\ref{tab:analog-digital}, 
%Table~{\color{red} 4} of the supplement,
assuming no array initialization (as was the case in our experiments),
the latter OP updates imply an estimated runtime not exceeding $m\times 0.14\mu\mbox{s}$ 
for the MFMWA-based approach, versus $m\times 250\mu\mbox{s}$ for the purely digital-based approach. 
Therefore, a computational improvement of up to three orders of magnitude can be realized, with similar improvements in energy consumption. The analog storage requirement of the algorithm on the MFMWA array is only $\ell n$, i.e., the size of the sketched matrix (the output).
% A graphical illustration is provided in
% %Figure \ref{fig:sketching}
% Figure~\ref{fig:simple-sketching}
% for a toy problem in which the columns and rows  of matrices $S$ and $A$ that form the OPs are colored alike. Each of these OP updates can be executed on the analog device
% %in $\mathcal{O}(1)$
% after performing digital-to-analog conversions on the corresponding vectors. 
The complete procedure is summarized in Algorithm~\ref{app:alg:algo0}, providing a pseudo-code description of our hybrid MFMWA implementation of randomized Gaussian sketching.

We remark that we could have implemented the Gaussian sketch by loading the sketching matrix $S$ onto the MFMWA device, and sketched the columns of $A$ one by one via the MV-hybrid operation. However, that would have resulted in an algorithm that requires $\ell m$ storage on the MFMWA device. Since the typical use-case of the Gaussian sketch is when $m\gg n$, this would have been wasteful.  
In fact, our hybrid MFMWA sketching algorithm is well suited for the streaming model of computation where, in many practical applications, the data matrix $A$ becomes available in chunks, e.g., one sample at a time (i.e., row). Moreover, the dataset might be so large that it can  not be stored in system memory even if we could afford to wait and process all data samples simultaneously. Randomized sketching is an important 
practical alternative since at each step of the algorithm we store the sketch in memory, and continuously update it as new samples arrive, without storing previous samples. In the hybrid randomized streaming algorithm the sketch, which should be small, is kept on the MFMWA device. It is continuously updated as new (digital) data becomes available.
We observe that, 
given an MFMWA OP update primitive, sketching and MFMWA computation work in synergy within the context of streamed Gaussian sketching.

Algorithm~\ref{app:alg:algo0} presents our general problem-solving framework of randomized sketching on MFMWA array devices.
Once we computed a Gaussian sketch we can employ it in a multitude of algorithms. As an example, 
Algorithm~\ref{app:alg:algo0} can be used in the context of linear regression, i.e., 
%minimize 
$\displaystyle \min_x \|Ax - b\|_2$. More specifically, we apply
% the 
Algorithm~\ref{app:alg:algo0} on the augmented matrix 
% $\left[\begin{array}{cc} A & b \end{array}\right]$
$[\: A \;\; b \:]$
to compute $Z$. Once we have $Z$, we split it as 
% $Z=\left[\begin{array}{cc} Z_A & z_b \end{array}\right]$. 
$Z=[\: Z_A \;\; z_b \,]$. 
Finally we output
$$
\tilde{x} = \arg\min_x\|Z_A x - z_b\|_2.
$$
Other problems that fit within the context of Algorithm~\ref{app:alg:algo0} include those involving an embedding of the column space of matrix $A$.
Algorithm~\ref{app:alg:algo0}
provides significant runtime and energy complexity improvements over the most efficient digital computer algorithms.

% Algorithm~\ref{app:alg:algo0} 
% presents our general problem-solving framework of randomized sketching on MFMWA array devices. 
% Problems that fit within the context of
% Algorithm~\ref{app:alg:algo0}
% include those involving
% ordinary least-squares 
% or an embedding of the column space of matrix $A$.

\begin{algorithm}[H]
\centering
\caption{\it Hybrid Analog-Digital  Gaussian Sketch }\label{app:alg:algo0}
\footnotesize
\begin{algorithmic}[1]
\STATE {\bfseries Input:} Matrix $A\in\mathbb{R}^{m\times n}$, sketch size $\ell\in \mathbb{Z}$
 \STATE {\bfseries Output:} Sketched matrix $Z\in\mathbb{R}^{\ell\times n}$
 %\newline
 \STATE Allocate on analog device: $\bar{Z}\in\mathbb{R}^{\ell\times n}$ 
\STATE  $\textrm{MReset}(\bar{Z})$
\FOR {$j=1,\ldots,m$}
\STATE Generate random vector $s\in\mathbb{R}^{\ell}$ of i.i.d.\ standard Gaussian entries $\mathcal{N}(0, 1)$
%\STATE $\textrm{OP-hybrid}(\bar{Z}, s, A^\top_{j,:})$
\STATE Update $\bar{Z} = \bar{Z} + s A^\top_{j,:}$
\ENDFOR
\STATE $\textrm{MRead}(Z, \bar{Z})$
\RETURN $Z$.
\end{algorithmic}
\end{algorithm}
%\end{minipage}
%\hfill

Table \ref{tab:cost0-olls}
%Table~{\color{red} 2} of the main paper
lists an asymptotic analysis of the number of 
arithmetic operations required by the above two 
computational OLLS tasks for: $(i)$ the standard approach typically used to perform the 
corresponding task in the published literature; $(ii)$ randomized sketching with 
JLT transformations, implemented on digital hardware; and $(iii)$ randomized 
sketching implemented on digital-analog hardware. In the latter case, 
the OP updates required to sketch $A$ are executed on an MFMWA 
array. For both problems considered, the number of arithmetic 
operations required by the hybrid variant is independent of the precision $\epsilon \in \mathbb{R}$ featured in the digital implementation. 
In addition, following the analysis presented in
% Section~\ref{app:ssec:model},
Appendix~\ref{app:ssec:model},
% % Section~{\color{red} A.5},
% Section~{\color{red} A.E},
the hybrid variant reduces the wall-clock time 
and energy consumption by up to three orders of magnitude each.

\begin{table*}[htb]
\caption{\it Asymptotic analysis of the number of arithmetic operations required to complete various computational tasks associated with randomized OLLS using purely digital-based and hybrid-based methods.}\label{tab:cost0-olls} 
%{\small
\begin{center} 
%\resizebox{1\textwidth}{!}{
\begin{tabular}{|c|c|}
\hline
Methods & Number of arithmetic operations \\
\hline
{\bf Computation of column subspace embedding} & \\
Direct digital SVD solution & $\mathcal{O}(mn^2)$ \\ 
Fully digital streaming solution & $\mathcal{O}(m\ell \mathtt{log}\ \ell+n\mathtt{poly}(\ell/\epsilon))$ \\ 
Hybrid MFMWA streaming solution & $\mathcal{O}(m +n \ell^2)$\\ 
%\hline
{\bf Computation of ordinary linear least-squares} & \\
Direct solution of equations & $\mathcal{O}(mn^2)$ \\ 
Fully digital streaming solution & $\mathcal{O}(m\ell \mathtt{log}\ \ell+\mathtt{poly}(\ell/\epsilon) + \ell n^2+n^3)$ \\ 
Hybrid MFMWA streaming solution & $\mathcal{O}(m+\ell n^2+n^3)$\\
%\hline
\hline
\end{tabular}
%}
\end{center} 
\end{table*}

\subsection{Randomized
%Principal Component Analysis}
PCA}\label{ssc:randPCA}
PCA is a fundamental unsupervised linear dimensionality reduction technique that 
transforms a set of correlated variables into another smaller set of linearly 
uncorrelated variables, called ``principal components'' (PCs) \cite{pearson1901liii,hotelling1933analysis}. 
From a mathematical viewpoint, PCA is directly related to rank-$k$ truncated Singular Value Decomposition (SVD). 
More formally, PCA with $k\in\mathbb{N}$ is equivalent to computing the $k$ leading left singular vectors of a matrix $A\in \mathbb{R}^{m\times n}$, which represents the \emph{column-centered} version of a data matrix $\widetilde{A}\in \mathbb{R}^{m\times n}$ encoding a dataset $\mathcal{A}$ of $m$ samples each having $n$ features;
i.e., $A$ stems from $\widetilde{A}$ after subtracting the column mean from all $n$ columns, followed by division with $\sqrt{m}$. 
The standard approach to compute these $k$ leading left singular vectors is to compute all $\min(m,n)$ singular triplets of $A$ with a runtime complexity of
% order $mn^2$.
$\mathcal{O}(mn^2)$.
As an alternative to reduce this complexity, a recent class of algorithms, which has been gaining prominence for the rapid approximation of a truncated SVD of large matrices, is that of randomized 
(probabilistic) subspace iteration \cite{gu2015subspace}. 
Randomized subspace iteration has been shown to outperform some of the best 
deterministic methods while also reducing the number of passes over large 
datasets, thus enabling the solution of previously infeasible problems. 
For example, randomized subspace iteration led to considerable speedups in PCA 
analysis of large-scale population genetics datasets~\cite{bose2019terapca}. 

We begin the description of randomized subspace iteration with its non-iterative variant, randomized SVD~\cite{halko2011finding}, a probabilistic approach that reduces the number of passes over large datasets, thus making it  possible to solve previously infeasible problems.
Given a matrix $A\in\mathbb{R}^{m\times n}$ and a target 
rank $k\ll \min(m,n)$, randomized SVD first constructs a random sketch matrix 
$R\in\mathbb{R}^{n\times \ell}$, $\ell\geq k$, whose entries are chosen from a distribution 
of choice, e.g., Gaussian or Count-Sketch. We 
then map $A$ to a low-dimensional subspace by computing the Matrix-Matrix product
\begin{equation*}
  Y := AR.
\end{equation*}
Once $Y\in \mathbb{R}^{m \times \ell}$ is formed, we compute an orthonormal 
basis of its column subspace by QR factorization, i.e., we compute matrices 
$Q\in \mathbb{R}^{m \times \ell}$ and $F\in \mathbb{R}^{\ell \times \ell}$ 
such that $[Q,F]=\mathtt{QR}(Y)$ and $\mathtt{range}(Q)=\mathtt{range}(Y)$. 
After forming matrix $Q$, we finally compute a low-rank projection of $A$, namely
\begin{equation*}
  B=Q^{\top}A,
\end{equation*}
and write 
\begin{equation*}
  A \approx QB=QQ^{\top}A=Q\widetilde{U}_k \widetilde{\Sigma}_k \widetilde{V}_k^{\top}=\widehat{U}_k\widetilde{\Sigma}_k \widetilde{V}_k^{\top},
\end{equation*}
where $\widetilde{U}_k \widetilde{\Sigma}_k \widetilde{V}_k^{\top}$ denotes the rank-$k$ 
truncated SVD of the matrix $B=Q^{\top}A$,
and thus the PCs of $A$ are finally approximated as $Q\tilde{U_k}$.

Randomized subspace iteration follows the 
same principles as randomized SVD, except that the matrix $A$ is now replaced by 
the matrix $(AA^{\top})^qA$, for some $q\in \mathbb{N}$. Indeed, the two matrices have the 
same left/right singular vectors, except that the singular values of the latter 
matrix decay much more rapidly. When implemented on a digital architecture, 
the computational cost of randomized subspace iteration consists of roughly 
$2m(nq \ell+n^2)$ arithmetic operations.

We present a pseudo-code description of our hybrid MFMWA implementation of randomized PCA in Algorithm~\ref{app:alg:algo1}. The crux of the algorithm is to compute $Y=(AA^{\top})^q A R$ column-by-column on the MFMWA device which, following
Table~\ref{tab:analog-digital}, 
%Table~{\color{red} 4} of the supplement,
can be achieved more than three orders of magnitude faster than the same operation performed in double precision on digital hardware, with similar improvements in energy consumption. In addition, this computational and energy performance gap widens as the value of $q$
increases. 
Furthermore, in dynamic settings, where matrix $A$ is updated as $A \gets A + C D{^\top}$, 
$C\in\mathbb{R}^{m\times r}$, $D\in \mathbb{R}^{n \times r}$, 
we can exploit MFMWA-based OP updates to further increase the asymptotic performance gap between the analog and digital algorithm. 
Therefore, the use of MFMWA devices allows us to perform PCA analysis of evolving datasets 
in a rapid fashion.

\if0 
A detailed sketch of our variant of randomized subspace iteration, which also considers 
the case of where the rows and columns of matrix $A$ are updated dynamically, is summarized 
in Algorithm~\ref{app:alg:algo1}. The procedure is divided into two separate phases, where 
the first phase is executed on an analog crossbar array device. More specifically, the 
matrix $A$ is loaded on the device, and the matrix product $AR=\sum\limits_{j=1}^{j=n}A_{:,j}R_{j,:}$ is computed as 
a series of $n$ OP updates executed in $\mathcal{O}(n)$. The matrix product between matrices 
$(AA^{\top})^q$ and $AR$ can be computed in $\mathcal{O}(q\ell)$. The second phase can be then 
executed in a digital environment in $\mathcal{O}(n\ell^2)$. This is a significant cost gain over efficient digital algorithms, which can become even greater when matrix $A$ is 
updated periodically. In particular, if the matrix $A$ receives row and column updates, say $C\in\mathbb{R}^{m\times p}$ and $D\in\mathbb{R}^{n\times p}$, with $p\ll \min(m,n)$, then the update $A\rightarrow A+CD^{\top}$ can be computed in-place \Lior{is there a technical term for that (such as in-situ) ?}
in the analog device using the OP update strategy, in just $\mathcal{O}(p)$ time. 
Therefore, the use of analog hardware allows us to perform PCA analysis of evolving datasets 
in a rapid fashion. 
\fi

Algorithm~\ref{app:alg:algo1} presents our general 
problem-solving framework of randomized subspace iteration on MFMWA devices, and also supports problem instances where the rows 
and columns of the matrix $A$ are updated dynamically. 
The first phase computes the matrix sketch on MFMWA hardware where 
the matrix product $AR=\sum\limits_{j=1}^{j=n}A_{:,j}R_{j,:}$ is computed 
as a series of $n$ OP updates.
For randomized subspace iteration with $q>0$, the additional products 
of matrices $A$ and $A^{\top}$ with $AR$ can be also computed using 
the MFMWA OP updates. 
The second phase is then executed on a digital computer. 
Algorithm~\ref{app:alg:algo1}
provides significant runtime and energy complexity improvements over the most efficient digital computer algorithms. Such improvements can become even greater when the matrix $A$ is updated periodically: If the matrix $A$ receives row and column updates, e.g., $C\in\mathbb{R}^{m\times p}$ and $D\in\mathbb{R}^{n\times p}$ with $p \ll \min(m,n)$,
then the update $A\rightarrow A+CD^{\top}$ can be computed in-place on the MFMWA device with our OP update strategy having a runtime complexity
of $\mathcal{O}(p)$. 
% of order $p$.

%\begin{minipage}{0.9\textwidth}
\begin{algorithm}[H]
\centering
% \caption{\it Hybrid Analog-Digital Randomized PCA Algorithm}
\caption{\it Hybrid Analog-Digital Randomized PCA}
\label{app:alg:algo1}
\footnotesize
\begin{algorithmic}[1]
\STATE {\bfseries Input:} Matrix $A\in\mathbb{R}^{m\times n}$ ($m\geq n$), rank $k\in \mathbb{Z}$, sketch size $\ell\in \mathbb{Z}$, power parameter $q\in \mathbb{Z}$
 \STATE {\bfseries Output:} Approximate Top $k$ principal components matrix $\widehat{U}_k$
% \newline
%\STATE {\bfseries Stage A:} Randomization phase (hybrid analog-digital)
\STATE Allocate on digital: $Y\in\mathbb{R}^{m\times \ell}$
\STATE  Allocate on analog device: $\bar{A}\in\mathbb{R}^{m\times n}, \bar{r}\in \mathbb{R}^{n}, \bar{z}\in \mathbb{R}^{m}$
\STATE $\bar{A} \gets \textrm{MLoad}(A,\bar{A})$
\FOR {$j=1,\ldots, \ell$}
\STATE  Generate random vector $r\in\mathbb{R}^n$ of i.i.d.\ standard Gaussian entries $\mathcal{N}(0, 1)$
\STATE $\textrm{VWrite}(r, \bar{r})$
\FOR {$i=1,\dots,q$}
\STATE $\textrm{MV-analog}(\bar{A},\bar{r}, \bar{z})$
\STATE $\textrm{M}^{\top}\textrm{V-analog}(\bar{A}, \bar{z},\bar{r})$
\ENDFOR
\STATE $\bar{z} \gets \textrm{MV-analog}(\bar{A},\bar{r},\bar{z})$
\STATE $\textrm{VRead}(Y_{:,j}, \bar{z})$
\ENDFOR
% \newline
%\STATE {\bfseries Stage B:} Projection phase (digital only)
\STATE  Compute $Q = \mathtt{orth}(Y)$
\STATE  Compute $B = Q^{\top}A$
\STATE  Set $[\widetilde{U}_k, \widetilde{\Sigma}_k, \widetilde{V}_k^{\top}] = \mathtt{svd}(B,k)$
\STATE  Set $\widehat{U}_k = Q\widetilde{U}_k$
\STATE  If required, update the analog crossbar array $A\rightarrow A+CD^{\top}$.
Repeat steps 6 to 19
\end{algorithmic}
\end{algorithm}

Table~\ref{tab:cost0-pca}
% %Table~\ref{tab:cost}
% Table~{\color{red} 2} of the main paper
% %of the supplement 
presents an asymptotic analysis of the number of 
arithmetic operations of a fully digital approach and our hybrid MFMWA approach, in terms of the matrix size $n$, the number of nonzero entries ${\tt nnz}(A)$, the desired number of 
%principal components
PCs
$k$ ($\sim  100$), the degree exponent $q$ ($\sim  5$) and the update rank $p$, when applied toward computational tasks associated with eigenvalue decomposition.
In the first problem instance the dataset is fixed and resides in primary memory, whereas in the second problem instance the dataset is streamed and evolves over time. 
%The hybrid analog-digital implementation performs the MV products on an analog array.
We observe that our hybrid approach provides significant reductions in terms of computational overhead over the fully digital approach, rendering speedups in both problem instances with respect to a multiplicative ${\tt nnz}(A)$ factor and an additional speedup in the second problem instance with respect to a multiplicative factor of the problem size (see, e.g., the discussion in 
% Section~\ref{app:ssec:model}).
Appendix~\ref{app:ssec:model}).
% % Section~{\color{red} A.5}).
% Section~{\color{red} A.E}).

\begin{table*}[htb]
\caption{\it Asymptotic analysis of the number of arithmetic operations required to complete various computational tasks associated with randomized PCA using purely digital-based and hybrid-based methods.}\label{tab:cost0-pca} 
%{\small
\begin{center} 
%\resizebox{1\textwidth}{!}{
\begin{tabular}{|c|c|}
\hline
Methods & Number of arithmetic operations \\
\hline
%\hline
{\bf Computation of eigenvalue decomposition in primary memory} & \\
Direct full eigenvalue decomposition solution & $\mathcal{O}(n^3)$ \\
Fully digital randomized PCA solution &$\mathcal{O}({\tt nnz}(A)qk+ nk^2 )$\\ 
Hybrid MFMWA randomized PCA solution & $\mathcal{O}(qk + nk^2 )$\\ 
%\hline
{\bf Computation of eigenvalue decomposition for evolving data}&\\
Fully digital update and sketch solution & $\mathcal{O}({\tt nnz}(A)qk+ nk^2+np^2 )$\\ 
Hybrid MFMWA update and sketch solution & $\mathcal{O}(qk+ nk^2+p^2 )$ \\
\hline
\end{tabular}
%}
\end{center} 
\end{table*}

\section{Experimental and Analytical Results} \label{sec:results}
We now turn to an evaluation of our
% % %framework
% % design of a
% new class of 
novel randomized algorithms tailored to MFMWA architectures that demonstrate
the potential for
orders of magnitude improvements
%with respect to
in
both computation and energy
%in comparison
compared
with
%classical algorithms on
digital computers.
%in order to break the existing performance bottlenecks in computational science
%Many of our
Our evaluation is primarily based on a combination of physical experiments
% , mathematical analysis, 
and simulation experiments.
%The physical experiments are based on measurements from a real Multi-Function Multi-Way Analog crossbar array device, 
% confirming
% %and supporting
% these more detailed simulation-based experiments and also
%demonstrating for the {\em first time} an implementation of an OP update mechanism on {\em physical analog hardware}.
%In order to consider much larger problem instances and allow us to demonstrate results for analog devices on our planned roadmap, the
%The
The
% simulation
% %Most
% experiments
latter
are based on our simulation engine 
incorporating models from current and near-term
%Multi-Function Multi-Way Analog 
MFMWA
arrays.
% ,
% %, the technical details of which are described in Appendix~\ref{app:ssec:exp}.
This
% which
includes simulated execution times of computations on digital devices that are obtained by assuming the digital hardware operates constantly at a peak performance of $10$
% teraflop/s (TFLOPS)
tera-floating-point operations per second (TFLOPS)
or higher\footnote{$10$ TFLOPS is representative of the peak performance of the NVIDIA A100 for double-precision computations.}
and ignoring any time spent to transfer data between processing unit(s) and system memory
(thus rendering optimistic digital execution times).
%; refer to Appendix~\ref{app:ssec:exp} for additional assumptions and settings of both the analog and digital computers.
% We also include measurements from a real Multi-Function Multi-Way Analog crossbar array device, confirming
% %and supporting
% these more detailed simulation-based experiments and also demonstrating for the {\em first time} an implementation of an OP update mechanism on {\em physical analog hardware}.
% Additional technical details are provided in
% %\ref{app:ssec:model}, \ref{app:ssec:algs} and \ref{app:ssec:exp}.
% the supplement.
These detailed simulation experiments of large-scale MFMWA devices on
% our planned 
% %device
% roadmap 
planned hardware roadmaps
confirm the results of the  physical experiments. %confirm these more detailed simulation-based experiments
%for large-scale analog devices.
%, where the latter allow us to consider much larger problem instances.
% Additional technical details are provided in
% %\ref{app:ssec:model}, \ref{app:ssec:algs} and \ref{app:ssec:exp}.
% the supplement.

The presentation of a collection of our performance evaluation results is organized as follows.
We first provide the outcomes from physical experiments of our randomized OLLS algorithm on a real chip prototype using analog devices, marking the first use of MFMWA operations in such an algorithm.
Then we present simulation-based experiments comparing our hybrid randomized OLLS algorithm against a fully digital version under different parameter settings, followed by simulation-based experiments comparing our hybrid randomized PCA algorithm against a fully digital version within the context of background subtraction in video surveillance.
Additional details on our experiments are provided in the supplement, together with additional analytical and experimental results.

\begin{figure*}
\centering
% \begin{subfigure}{0.48\textwidth}
\begin{tikzpicture}[thick,scale=0.77, every node/.style={scale=1.3}]
\definecolor{clr2}{RGB}{31,182,83}
	\begin{axis}[grid=both,
    grid style={line width=.1pt, draw=gray!10},
    major grid style={line width=.2pt,draw=gray!50},
		xlabel= \# of pulses,ylabel near ticks,
		ylabel= Classification accuracy,
		legend pos=south east,
		legend style={font=\fontsize{6}{4}\selectfont},
		title={Physical experiment},
		title style={at={(0.5,1.09)},anchor=north,yshift=-0.1},]
	\addplot[color=olive,mark=square,line width=2pt,mark size=4pt] coordinates {
		(15,0.9833)
		(31,0.9833)
		(63,0.9833)
	};	\addlegendentry{Digital baseline}
	\addplot[color=orange,mark=o,line width=2pt,mark size=4pt] coordinates {
		(15,0.9631)
		(31,0.9631)
		(63,0.9631)
	};\addlegendentry{Digital streaming}
	\addplot[color=violet,mark=triangle,line width=2pt,mark size=4pt] coordinates {
		(15,0.647)
		(31,0.70)
		(63,0.936)
	};\addlegendentry{Analog streaming}
	\end{axis}
\end{tikzpicture}
% \end{subfigure}\hfill
\qquad\qquad
% \begin{subfigure}{0.48\textwidth}
% \includegraphics[width=8.0cm, height=5.68cm]{physical3.eps}
\includegraphics[width=8.4cm, height=6.0cm]{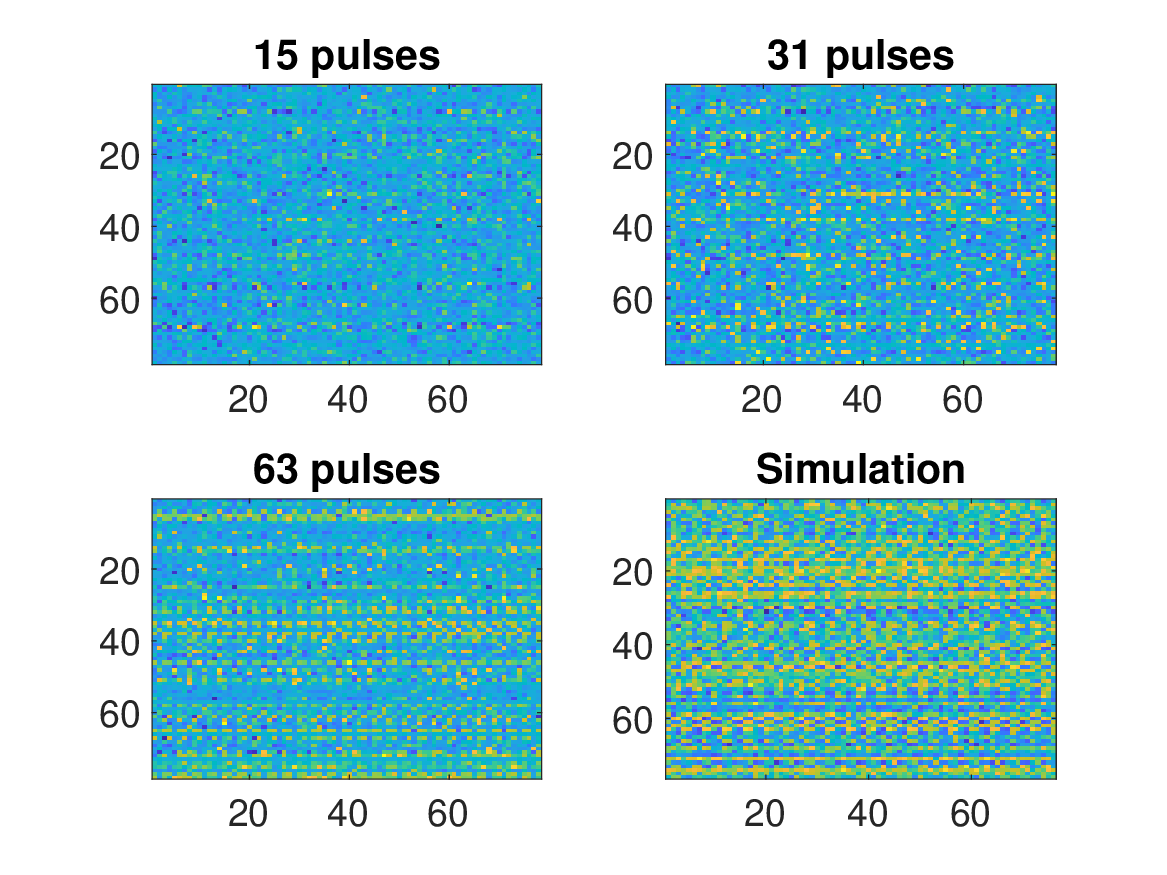}
% \end{subfigure}
\caption{\textit{Left: Classification error rate achieved by the baseline and randomized 
streaming approaches on a purely digital computational environment versus that of randomized 
streaming implemented through an 
%Multi-Function Multi-Way Analog
MFMWA
%crossbar 
array. Right: 
% Time ratio between digital and hybrid randomized PCA as the matrix size and number $k$ of sought PCs vary.
Cosine similarity between the ideal solution and the regressors (solutions) returned by simulation and actual hardware. In the latter case, the number of pulses was varied as 
15, 31, and 63 (corresponding to 4-bit, 5-bit, and 6-bit resolution, respectively).}}
\label{fig:physical3}
\end{figure*}
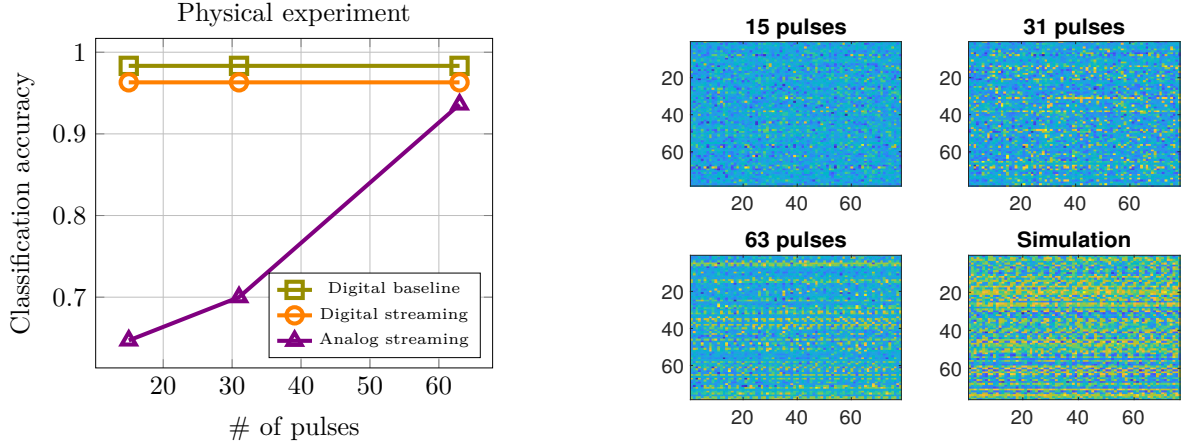

\subsection{Physical Experiments}
%\paragraph{Physical Experiments.}
We first consider a set of classification experiments with our randomized streaming
Algorithm~\ref{app:alg:algo0},
% Algorithm~$1$,
%that leverages our core analog compute primitives, 
marking the first time an OP update mechanism, with a randomized streaming algorithm, is implemented on actual analog hardware. 
These experiments are conducted on a physical
%Multi-Function Multi-Way Analog 
MFMWA
array
%device 
originally designed to perform deep
%neural network 
NN
training using parallel MV
%multiplies
products
and OP updates~\cite{kohda2020IEDM}, where the matrix entries are stored on Metal-Insulator-Metal Capacitors (MIMCAP).
%; refer to Appendix~\ref{app:ssec:exp} for all technical details on the experiments including the analog devices and the algorithmic approaches.
% Our set of physical experiments leverages the same analog compute primitives, marking the first time an OP update mechanism together with a randomized streaming algorithm are implemented on actual analog hardware. 
%
The fully analog ASIC
%DNN
training chip is designed and fabricated in CMOS 90nm technology, where the MIMCAP unit cell area is $114$/um$^2$ and the ASIC clock frequency is $100$MHz.
%A detailed structure of the crossbar together with peripheral circuit (ASIC block diagram) is discussed  in our supplement.
The MFMWA chip has two sets of SRAM and pulse-width modulation to provide pulses to the analog multiplier/adder elements.
% In addition, the 
The chip also
includes two integrators and ADC elements to convert total 
capacitor voltages to digital values.

To solve the classification problem,
we apply 
%linear least squares 
LLS
regression.  
The experimental dataset consists of $8192$ three-dimensional points distributed in the unit cube centered around the origin of a three-dimensional Cartesian coordinate system,
% \footnote{A visualization of these three-dimensional points can be found in
% %Figure~\ref{app:fig:physical1}.
% Figure~{\color{red} 9} of the supplement.
% }
choosing half of the points to have positive first coordinate
% (assigned class ``$+1$'') 
(class ``$+1$'') 
and the remaining to have negative first coordinate
%(assigned class ``$-1$''). 
(class ``$-1$''). 
After randomly dividing these points into training and testing datasets each of size $4096$, we construct a $4096\times 76$ matrix $A=\left[A^{(1)},\ldots,A^{(19)}\right]$ where $A^{(1)}=\cdots=A^{(19)}$, $A^{(i)}=[x,y,z,t]$, and $x,y,z \in \mathbb{R}^{4096}$ respectively denote the first, second and third coordinate of the training samples, while $t\in\{-1,+1\}^{4096}$ denotes the corresponding
%class
assigned
label. A $76\times 4096$ sketching matrix $S$ with $\pm 1$ entries was employed as part of $4096$ OP updates involving the rows and columns of the $S$ and $A$ matrices, respectively. We then extract the resulting $76\times 76$ matrix $Z=SA=[Z^{(1)},\ldots,Z^{(19)}]$ from the 
%Multi-Function Multi-Way Analog 
MFMWA
array and solve the corresponding 
%linear least squares 
LLS
problem on a digital computer to obtain $w_i \in \mathbb{R}^3$ for each $Z^{(i)}\in \mathbb{R}^{76\times 4}$. Upon averaging $w_1,\ldots,w_{19}$ to obtain the global regressor $w$, the label of each sample is set equal to the sign of the dot product
%(i.e., $-1$ if the dot product is negative, and $+1$ otherwise)
between $w$ and each three-dimensional sample of the testing dataset. 

We repeat each experiment for three different input resolutions, namely $6$-bit, $5$-bit, and 
$4$-bit (i.e., $63$, $31$, and $15$ pulses). Note that a $1$-bit decrease in resolution
%1-bit lower resolution
halves the duration of the applied inputs,
while the integration time of the ADCs does not decrease.
Figure~\ref{fig:physical3} (left) plots the classification accuracy achieved by the fully digital and hybrid randomized streaming algorithm implementations, contrasted with the digital baseline where the data samples are classified through
%linear least-squares 
LLS
regression. 
%Randomized streaming algorithms achieve lower accuracy since they project the dataset onto a %lower-dimensional subspace, thus leading to some loss of information. 
Our MFMWA approach achieves essentially the same accuracy as the corresponding digital approach, when the  
number of pulses in the encoding is sufficiently large, while at the same time exhibiting significant
%reductions 
computation and energy improvements
over the digital approach of up to three orders of 
% magnitude;
% % magnitude improvement with respect to computation and energy;
% this is consistent with
% %the results in
% Table~\ref{tab:analog-digital}
% %Table~{\color{red} 4}
% and
% Table~\ref{tab:cost0}.
% % Table~{\color{red} 5}
% % of the supplement. 
magnitude
(which is consistent with
the results in
Table~\ref{tab:analog-digital},
Table~\ref{tab:cost0-olls},
and
Table~\ref{tab:cost0-pca}).
% improvements consistent with the results in
% %Tables~\ref{tab:analog-digital}
% Tables~{\color{red} 4}
% %and~\ref{tab:cost0}
% and~{\color{red} 5}
% %of the supplement.
% (supplement).
In addition, Figure~\ref{fig:physical3} (right) plots the cosine similarity between the  least-squares solution and the individual regressors returned by randomized sketching; 
both by simulation and physical hardware. In agreement with the left subfigure, where 
higher input resolution improved the obtained classification accuracy, increasing the 
input resolution leads to regressors that are better aligned with the simulated solution, 
while also reducing the variance of the individual regressors.

\subsection{Sketching of Streamed Matrices} \label{sec:sketch}
We next consider simulation-based experiments with our randomized streaming
Algorithm~\ref{app:alg:algo0}.
% Algorithm~$1$.
%(supplement).
% Following the computational model and analysis presented in Section A.5 
% of our supplement, Table~{\color{red} 5} (supplement) presents an 
% asymptotic comparison of the number of arithmetic operations required by 
% the fully digital approach and our hybrid  analog-digital approach, in 
% terms of the matrix dimensions $m$, $n$ and the accuracy precision 
% $\epsilon$, when applied toward two computational tasks: $(a)$ column subspace embedding of matrix $A$; and $(b)$ ordinary linear least-squares. % problem $\widehat{x}=\argmin\limits_{x\in \mathbb{R}^n}\|Ax-b\|_2$.
% %where $A\in \mathbb{R}^{m\times n}$ and $b\in \mathbb{R}^m$.
An asymptotic comparison of the number of arithmetic operations required by the fully digital approach and our hybrid MFMWA approach, in terms of the matrix dimensions $m$, $n$ and the accuracy precision $\epsilon$, shows a reduction from
$\mathcal{O}(m\ell\, \mathtt{log}\, \ell+n\mathtt{poly}(\ell/\epsilon))$ 
% on the order of $m\ell\, \mathtt{log}\, \ell+n\mathtt{poly}(\ell/\epsilon)$ 
to 
$\mathcal{O}(m +n \ell^2)$ 
% on the order of $m +n \ell^2$ 
for the computational task of column subspace embedding of matrix $A$;
and shows a reduction from 
$\mathcal{O}(m\ell\,  \mathtt{log}\, \ell+\mathtt{poly}(\ell/\epsilon) + \ell n^2+n^3)$ 
% on the order of $m\ell\,  \mathtt{log}\, \ell+\mathtt{poly}(\ell/\epsilon) + \ell n^2+n^3$
to 
$\mathcal{O}(m+\ell n^2+n^3)$ 
% on the order of $m+\ell n^2+n^3$ 
for the computational task of
% %ordinary linear least-squares.
% OLLS.
OLLS; refer to Table~\ref{tab:cost0-olls}.
In both computational tasks, the data matrix $A$ is streamed by rows.
%The hybrid analog-digital implementation performs the sketch of $A$ (i.e., series of $m$ OP updates) on an Multi-Function Multi-Way Analog crossbar array.
We observe that our hybrid approach provides significant reductions over the 
fully digital approach where the speedups achieved by the hybrid approach grow super-linearly with respect to the sketching dimension $\ell$. Moreover, 
% following the time and energy complexities described in Table~{\color{red} 4} of our supplement,
these reductions translate to up to three orders of magnitude improvements in both computation and energy; refer to Table~\ref{tab:analog-digital}.
% ;
% % rendering speedups in both tasks with respect to a multiplicative $\mathcal{O}(\ell\log\ell)$ factor and an additive $\mathtt{poly}(\ell/\epsilon)$ factor. 
% % It is important to further emphasize that the running time of the hybrid approach is independent of the precision $\epsilon \in \mathbb{R}$ featured in the running time of the fully digital approach, \Haim{Unclear here what $\epsilon$ is? Also: how can it be independent of it?} with the additive $\mathtt{poly}(\ell/\epsilon)$ factor blowing up as the desired precision grows large (i.e., as $\epsilon\rightarrow 0$).
% % \Haim{Shouldn't the results section focus on the experimental results, and any complexity results should actually be be in the previous section?}
% refer to Table~\ref{tab:analog-digital}.
% % %Table~{\color{red} 4}
% % % and Table~{\color{red} 5}
% % of the supplement.

\begin{figure*}[tb!]
    \centering
    \includegraphics[width=0.37\textwidth]{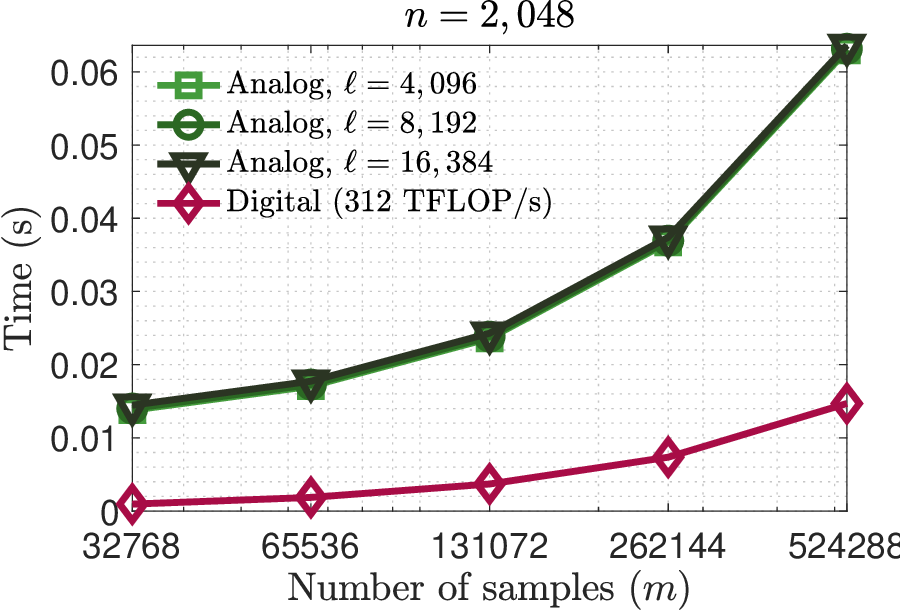}\hspace*{0.5in}
    \includegraphics[width=0.37\textwidth]{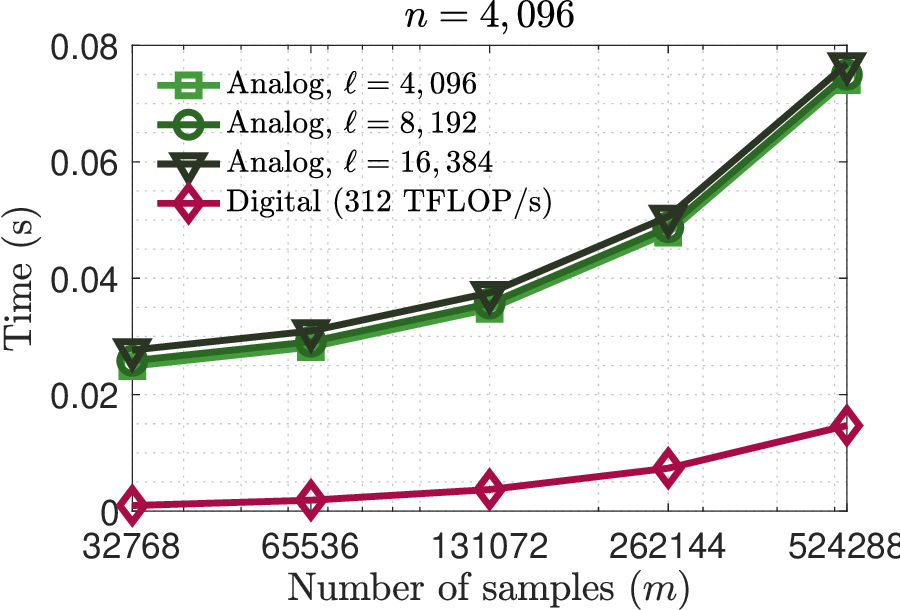}\vspace*{0.1in}
    \includegraphics[width=0.37\textwidth]{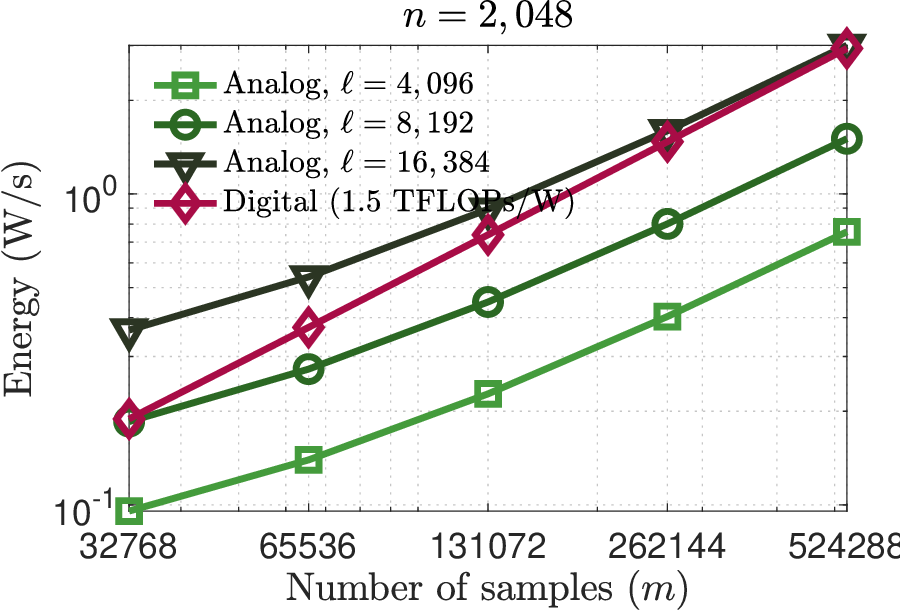}\hspace*{0.5in}
    \includegraphics[width=0.37\textwidth]{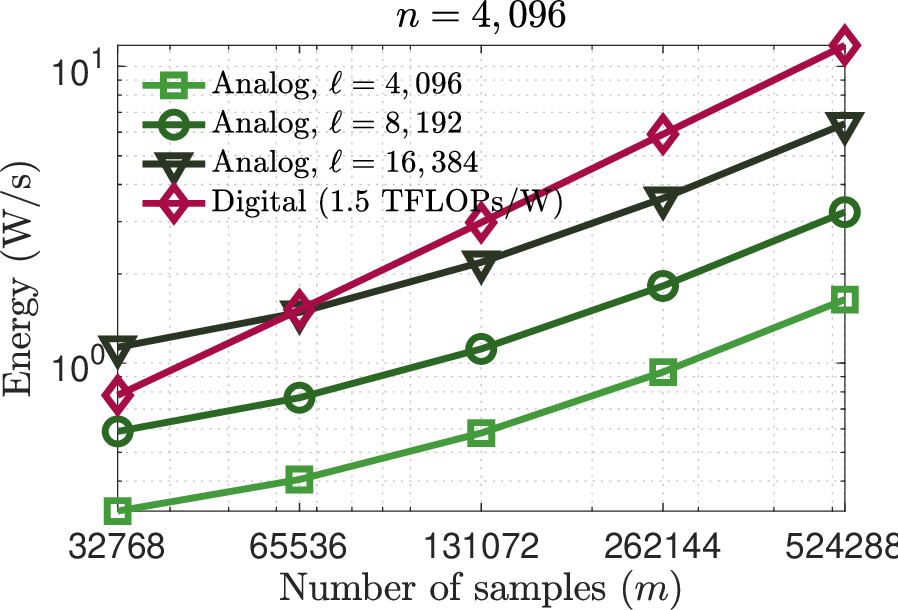}\vspace*{0.1in}
    \includegraphics[width=0.37\textwidth]{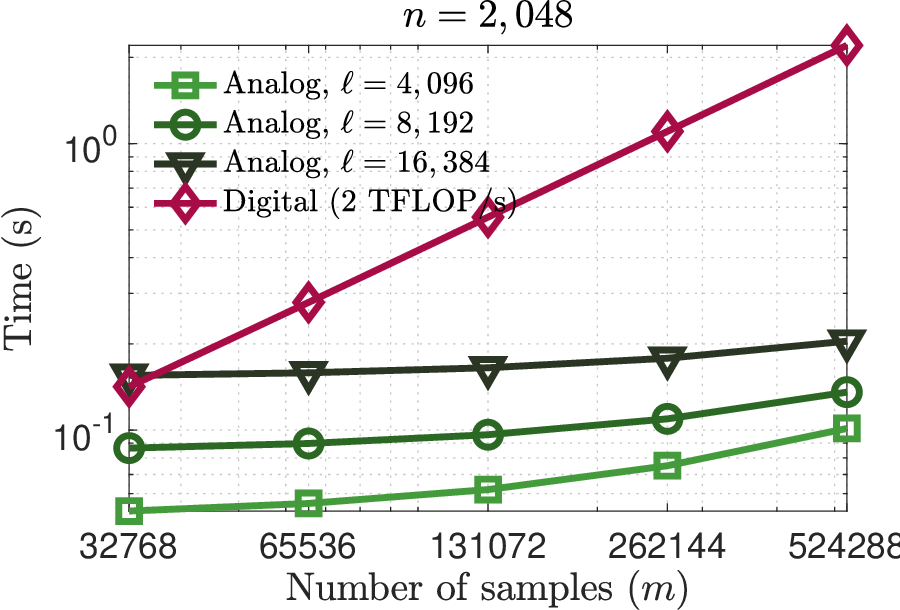}\hspace*{0.5in}
    \includegraphics[width=0.37\textwidth]{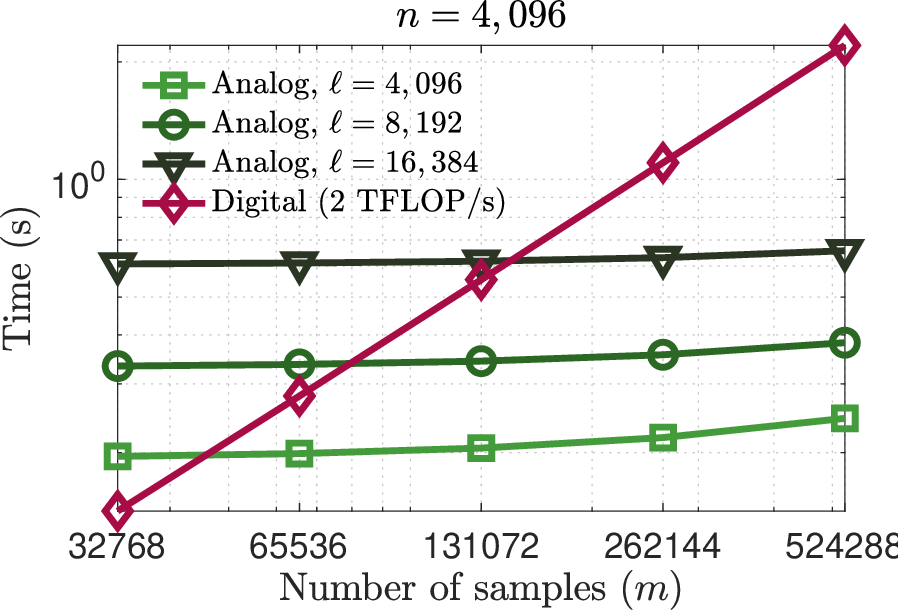}
    \caption{\it Simulated wall-clock times and estimated energy consumption of 
    analog-digital vs.\ digital (without sketching). The number of samples 
    $m$ and sketching dimension $\ell$ vary. The number of features of each sample 
    was set to $n=2048$ (left) and $n=4096$ (right). } \label{sketch_flops1}
\end{figure*}

Figure~\ref{sketch_flops1} plots the simulated execution (wall-clock) times of the fully digital and hybrid sketching algorithms for different numbers of samples $m$ and various sketching dimensions $\ell$. The number of features for each sample is set to $n=2048$ (left) and $n=4096$ (right).
%The simulated execution times of the digital implementation are obtained assuming that the digital hardware operates constantly at a peak performance equal to $10$ teraflops (TFLOPS), which is representative of the peak performance of the NVIDIA A100 for double-precision computations; thus, we ignore any time spent on transferring data between the processing unit(s) and the system memory. 
%On the other hand, the write overhead of the analog device was set equal to $10^{-2}$ seconds, while the time to perform a single OP update was set equal to $10^{-7}$ seconds.
% Consistent with the results in
% %Tables~\ref{tab:analog-digital}
% Tables~{\color{red} 4}
% %and~\ref{tab:cost0}
% and~{\color{red} 5}
% %of the supplement,
% (supplement),
% we 
We
observe from Figure~\ref{sketch_flops1} that our hybrid MFMWA approach provides significant  improvements over the fully digital approach, where the
%performance
computational
and energy performance gaps between the two implementations widen as $m$ grows much larger than $n$. 
In particular, for sketching dimensions requiring access to the main system memory, our hybrid approach can be about $20$ times faster and $10$ times more
%economical in energy consumption
energy efficient
than a state-of-the-art digital computer, with both factors increasing as $m$ grows. 
%\CWW{What is the sketching dimension $l$ for the fully digital algorithm in Figure 3?}

 \begin{figure*}[tb!]
 \centering
 \includegraphics[height=0.2\textwidth,trim={6cm 1cm 6cm 1cm}]{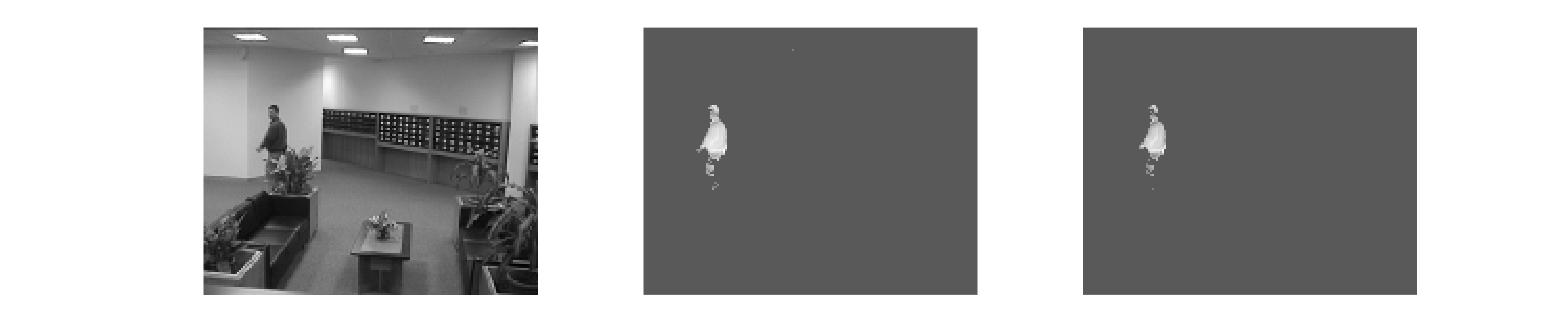}
\includegraphics[height=0.2\textwidth,trim={6cm 1cm 6cm 1cm}]{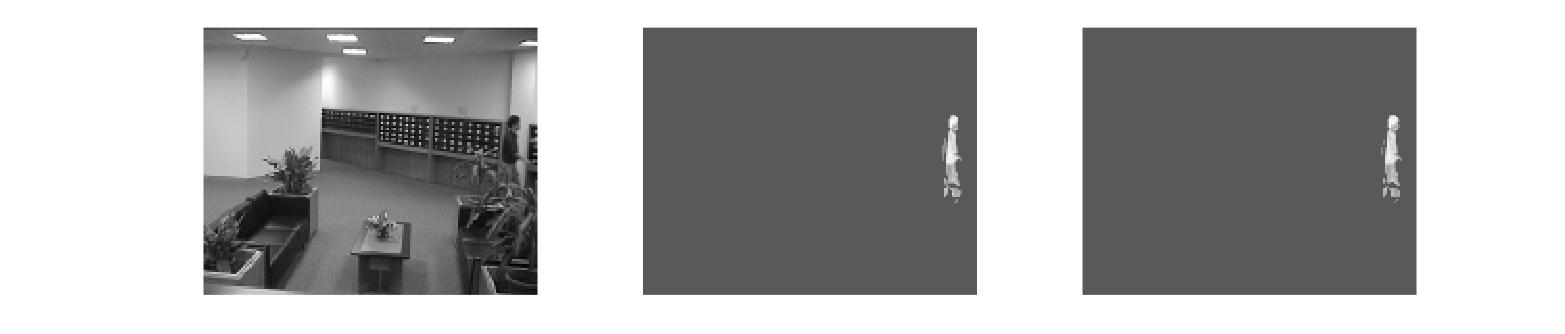}
 \includegraphics[height=0.2\textwidth,trim={6cm 1cm 6cm 1cm}]{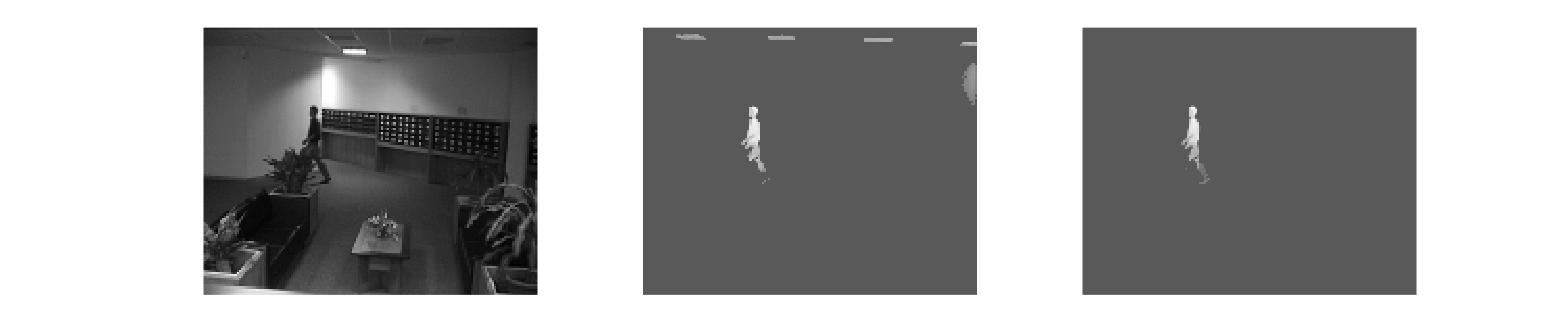}
  \includegraphics[height=0.2\textwidth,trim={6cm 1cm 6cm 1cm}]{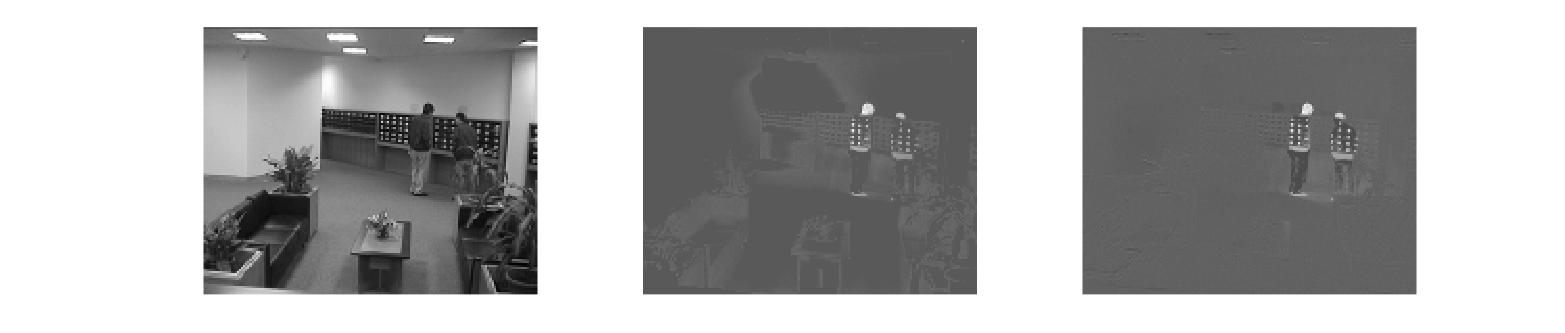}
 \caption{\it Video Background Subtraction (Lobby dataset).  Four sample frames from the video (leftmost column);  Background subtraction obtained using hybrid PCA (center column) and digital PCA (rightmost column).} 
 \label{fig:video}
 \end{figure*}%\\[-.3cm]

\subsection{Principal Component Analysis}
Lastly, we consider simulation-based experiments with our randomized PCA Algorithm~\ref{app:alg:algo1}.
An asymptotic comparison of the number of arithmetic operations required by the fully digital approach and our hybrid MFMWA approach,
in terms of the matrix size $n$, the number of nonzero entries ${\tt nnz}(A)$, the desired number of PCs $k$ ($\sim  100$), the degree exponent $q$ ($\sim  5$) and the update rank $p$, 
shows a reduction from
$\mathcal{O}({\tt nnz}(A)qk+ nk^2 )$
to 
$\mathcal{O}(qk + nk^2 )$
for the computational task of eigenvalue decomposition in primary memory;
and shows a reduction from 
$\mathcal{O}({\tt nnz}(A)qk+ nk^2+np^2 )$
to
$\mathcal{O}(qk+ nk^2+p^2 )$
for the computational task of eigenvalue decomposition for evolving data; refer to Table~\ref{tab:cost0-pca}.
We observe that our hybrid approach provides significant reductions over the fully digital approach where the speedups achieved by the hybrid approach grow linearly with respect to the number of nonzero entries ${\tt nnz}(A)$ and the matrix size $n$. Moreover, 
these reductions translate to up to three orders of magnitude improvements with respect to both computation and energy; refer to Table~\ref{tab:analog-digital}.

Next, we consider a real-world data application of our
Algorithm~\ref{app:alg:algo1}
% Algorithm~$2$
in the context of the computer vision task of \emph{background subtraction} in surveillance videos. %The objective is to exploit PCA to separate the foreground information from the redundant background in videos.
%
%For the  computer vision application, 
In particular,
we study the dataset
% titled: the
% %so-called
``Lobby in an office  building with switching on/off lights'', which
was introduced in~\cite{LiHuGuTi04}, also used in~\cite{candes2011robust}, and is
publicly available at \url{http://perception.i2r.a-star.edu.sg/bk_model/bk_index.html}.  This Lobby video contains $1546$ frames, each of size $160\times128$, and thus the size of the data matrix is equal to $1546\times20480$.
%We consider the video dataset: “Lobby in an office  building  with switching  on/off  lights” available  from \url{http://perception.i2r.a-star.edu.sg/bk_model/bk_index.html}. 
%The Lobby video contains $1,546$ frames each of size $160\times128$, and the data  matrix  size  is  $1,546\times20,480$. 
%Figure~\ref{fig:video} presents the results obtained by our hybrid randomized PCA  and fully digital PCA with rank $k=5$. 
For background subtraction using PCA, we use the following procedure. The video frames are first vectorized into a matrix and mean-centered. Next, the top $k$ singular vectors of this mean-centered matrix are computed using a PCA method (e.g., the hybrid PCA approach presented in this paper). Next, the frames are projected on to the PCs, which capture the dominant background information. Finally, the projected matrix is used to remove (subtract)  the background from the original frames.
% In
% Figure~\ref{fig:video},
% % Figure~4,
% we present the results we obtained after background subtraction.

Figure~\ref{fig:video} presents the corresponding qualitative results obtained using the fully digital and hybrid randomized PCA algorithms with rank $k=5$. The leftmost column of images are (four) true sample frames from the video. The center column of images are the background subtraction results obtained for the four frames using our hybrid PCA approach, and the rightmost column of images are the results from the fully digital PCA approach.
%Details on the procedure used for background subtraction using PCA are provided in the supplement.
We first observe that there are little to no noticeable differences between the results obtained by the digital approach and our hybrid approach in most cases. Moreover, any differences appear to be quite small and do not adversely impact the goal of accurately separating the foreground information from the background. These results show that our hybrid randomized PCA achieves comparable accuracy results, while providing  significant computational and energy improvements over the fully digital approach, as demonstrated above.

\section{Conclusions}
\label{sec:concl}
Numerical computation has been and remains critical to many areas of MI and AI, but the exponential growth in their computation and energy requirements creates a perfect storm for computational sustainability.
This problem is exacerbated by the limitations 
% moving forward 
of digital technology to support such exponential growth
%, which greatly exceeds the rate of improvements in digital hardware.
due to the slowdown in Moore's law.
One potential solution is based on MFMWA technology
% , a promising computing architecture comprising arrays of memristors that 
that can 
support the computation of MV products and OP updates in time that is independent of the matrix size (subject to physical constraints and architecture design) with comparable improvements in energy consumption.
These benefits, however, come at the expense of inherent unpredictability and noise due to hardware imperfections and temporal conductance variations.

In this paper, we design, analyze and implement a general class of randomized numerical linear algebra algorithms for core numerical methods involving MV products and OP updates, tailored to mitigate the detrimental impact of imperfect numerical computations on MFMWA technology while realizing its potential benefits.
We derive 

\newpage
\noindent
an analysis of the core numerical methods of our algorithms on both a hybrid MFMWA accelerator architecture and a purely digital accelerator architecture.
We then present measurements from experiments conducted on a real MFMWA device with our randomized streaming algorithm,
% which marks the first time an OP update within such an algorithm is implemented on actual analog hardware.
followed by detailed simulation-based experiments that compare our randomized streaming and randomized PCA algorithms on larger hybrid MFMWA architectures against a fully digital version, including pattern analysis in computer vision.
This collection of analytical and experimental results demonstrate major advances in numerical speedups, energy reductions and computational sustainability, such as orders of magnitude improvements over state-of-the-art algorithms on digital computers with similar accuracy. 
Our algorithms and results can help to unlock the true potential of sustainable numerical computation and its impact on various application areas of MI and AI.

\onecolumn

\appendix

\section{Supplemental Material}
This supplement supports the main body of the paper by providing additional background material and additional technical details and results.
Standard mathematical notation is employed throughout portions of the paper, which we first briefly define.
Additional background on numerical computation, MFMWA technology and hybrid MFMWA architecture technology is  provided in turn.
We then present a computational model and analysis of the core primitives of our algorithms executed on the hybrid MFMWA accelerator architecture in comparison with execution on a purely digital accelerator.
Additional technical details on our 
% randomized sketching and randomized PCA algorithms tailored to hybrid MFMWA architectures
physical experiments
are provided next.
% , followed by additional technical details on our physical experiments.
We then present our simulation engine and additional details on our simulation experiments for randomized sketching and randomized principal component analysis (PCA) obtained with the use of this engine,
including results from a real-world data application of our randomized PCA algorithm in the analysis of population structure in human genetics.
A complementary mathematical analysis of the relationship between analog noise and algorithm convergence when an approximate inverse preconditioner is used on an MFMWA array to solve sparse linear systems is provided in~\cite{Kalantzis2021-nb}.

\subsection{Mathematical Notation}
\label{app:ssec:notation}
Intuitively, big-$\mathcal{O}$ notation states that the absolute value of a function is bounded above by another function asymptotically.
More formally, we say that $f(x) = \mathcal{O}(g(x))$ if there exists a positive real number $M$ and a real number $x_0$ such that $|f(x)| \leq M g(x)$ for all $x \geq x_0$.
Similarly, big-$\Omega$ notation states that the absolute value of a function is bounded below by another function asymptotically.
More formally, we say that $f(x) = \Omega(g(x))$ if there exists a positive real number $M$ and a real number $x_0$ such that $|f(x)| \geq M g(x)$ for all $x \geq x_0$.
Lastly, big-$\Theta$ notation states that the absolute value of a function is bounded both above and below by another function asymptotically.
More formally, we say that $f(x) = \Theta(g(x))$ if there exists positive real numbers $M$ and $N$ and a real number $x_0$ such that $M g(x) \leq |f(x)| \leq N g(x)$ for all $x \geq x_0$.

The notation $\mathtt{nnz(A)}$ denotes the number of nonzero entries in the matrix $A$.
The notation $\mathtt{poly}(x)$ denotes a polynomial value of $x$.
Lastly, the notation $\mathcal{N}(\mu, \sigma)$ denotes the normal distribution with mean $\mu$ and standard deviation $\sigma$.

\subsection{Numerical Computation}
\label{app:ssec:numerical}
Numerical computation
is central to understanding
complex phenomena throughout our world, lying at the forefront of machine intelligence (MI), artificial intelligence (AI), scientific,  and other disciplines.
Such
computation
is often concerned with the deep analysis of various aspects of our universe at large scales.
%\Lior{can we somehow also claim AI?}.
In 
AI~\cite{lecun2015deep},
numerical
computation similarly involves the analysis and learning of various high-dimensional aspects of our daily life.
These diverse forms of analysis are all performed through the execution of computer programs that run on high-end computing infrastructures with the goals of advancing 
discovery and understanding of
complex systems.
Our ability to efficiently perform these computations and predictions of fundamental behaviors at large scales and high dimensions is crucially important for major advances across many fields.

The ever-increasing, inconceivable levels of detail and realism
in the numerical computations
across
disciplines
have been achieved
through amazing advances in the speed of digital computers,
from
the dawn of early computers up to the computing facilities of today.
An important well-known example of numerical computation playing an integral part in successfully solving many grand challenges at the forefront of various disciplines is 
the analysis of statistical patterns within privacy-preserving collaborative 
genome-wide association studies 
%(GWAS)~\cite{bush2012genome,kuleshov2019machine}.
(GWAS)~\cite{Uffelmann2021,kuleshov2019machine,Mbatchou2021}.
Additional important well-known examples include
the analysis of features and patterns in computer imaging and vision~\cite{turk1991eigenfaces,LiHuGuTi04,
candes2011robust,gittens2016matrix,bouwmans2016handbook},
and 
the analysis of population structure in human genetics~\cite{paschou2007pca,bose2019terapca}.
Other successful high-profile applications of numerical computation include those in areas related to
geography~\cite{armstrong2000geography},
astronomy~\cite{yang2015implementing} and cosmology~\cite{almgren2013nyx,bedorf201424},
seismology~\cite{werner2005monitoring},
climate science~\cite{giorgi1991approaches},
signal processing~\cite{comon1990tracking,yang1995projection},
% image and video processing~\cite{turk1991eigenfaces,candes2011robust,bouwmans2016handbook},
natural language processing~\cite{berry1995using,vaswani2017attention},
and drug discovery~\cite{kendall2000high,sliwoski2014computational}.
The computational requirements of these  large-scale applications necessitate the processing and analysis of very large volumes of data. Such large data volumes can be generated at a rate that exceeds current advances in memory and/or processor hardware, which in turn create challenging performance and energy bottlenecks. It is therefore critically important to be able to analyze data and infer patterns by processing large data collections in order to tackle the major performance and energy bottlenecks arising in
computational
applications
across 
MI, AI, 
scientific,
and other disciplines.
As described in the main body of the paper, however, this purely digital approach can no longer sustain the growing requirements of modern applications for computational speed and energy consumption.

We consider an important class of algorithms at the core of these numerical computations.
More specifically, we consider randomized algorithms and randomized numerical linear algebra, the latter being a multi-disciplinary field that exploits randomization as a mean to design faster and more practical algorithms for the solution of large-scale linear algebra problems.
The main body of the paper primarily focuses on the randomized algorithms of matrix sketching and PCA, which represent two of the most commonly employed approaches among the various numerical algorithms supporting computational applications.
Successful applications of these randomized algorithms are quite numerous.
This includes matrix sketching applied to cryptography of privacy-preserving collaborative GWAS that improved the runtime of the most advanced cryptographic protocols by two orders of magnitude~\cite{kockan2020sketching}; another example concerns matrix sketching becoming a critical component for the practical solution of large-scale ordinary least squares problems~\cite{avron2010blendenpik}.
This also includes PCA applied to 
mass spectrometry imaging and vision~\cite{gittens2016matrix,LiHuGuTi04,candes2011robust},
population genetics~\cite{paschou2007pca,bose2019terapca}, cosmology, and climate science~\cite{gittens2016matrix}; 
PCA is routinely employed to identify structure and ancestry relations in the distribution of genetic variation across geographical location and ethnic background~\cite{paschou2007pca,reich2008principal}; and PCA is popularly used in image and video processing tasks such as face recognition~\cite{turk1991eigenfaces} and video background subtraction~\cite{candes2011robust,ubaru2019find}.
It is for these reasons that we focus primarily in this paper on accelerating matrix sketching and PCA randomized algorithms.
By speeding up these core algorithms, we can have dramatic widespread impact across many fields, including MI and AI.

\subsection{MFMWA Technology}
\label{app:ssec:analog}
The use of arrays of Non-Volatile-Memory (NVM) to perform multiply-accumulate operations was recently proposed for the acceleration of deep learning
workloads \cite{gokmen2016acceleration,burr2015experimental,burr2017neuromorphic}, although the basic concept of using resistive arrays for matrix-vector multiplications dates back to the early 1960s \cite{cite-key}. NVM materials are intrinsically noisy and come with some amount of material variations that can be tolerated by certain algorithms \cite{gokmen2016acceleration,gokmen2020algorithm}. The basic mathematical operations that can be performed on these arrays are (1) matrix-vector (MV)
%multiplication
products
and (2) rank-$1$ or 
outer-product (OP)
updates. For a dense $N\times N$ matrix, both these operations require roughly $2N^2$ multiply and add operations on conventional digital computers.
%For
%%matrix-vector (MV) multiplication 
%MV products
%with a dense $N\times N$ matrix, roughly $2N^2$ multiply and add operations are required;
%%a rank-$1$
%an OP
%update on such a matrix will amount to the same. 
If these mathematical operations are performed on an NVM array, %the MV
%multiplication
%products
%and
%rank-$1$
%OP
%updates 
they can be performed in constant time irrespective of the matrix size subject to physical constraints and architectural design.

This is illustrated in
Figure~\ref{revised-figna1}(a)
% Figure~1(a)
which shows a two-dimensional array with horizontal and vertical metal lines that are connected at the cross points by an NVM element. Each NVM element has an associated conductance $g_{kl}$.  
In order to perform an MV operation, the digital input vector $x$ is converted into time-encoded signals of maximum duration $t_{meas}$ and applied along the rows of the NVM array. Ohm's and Kirchhoff's laws result in integrated charges in charge capacitors along the columns yielding a voltage vector $V$. 
%An MV operation is accomplished by converting the digital input vector $x$ into a time encoded signal of maximum length $t_{meas}$, and the physics according to Ohm’s and Kirchhoff’s laws will result in an integrated charge at each column represented by a voltage vector $V$ in a charged capacitor. 
Simple network analysis renders
\begin{equation}\label{app:na1a}
%. \sum\limits_{k=1}^N x_k g_{kl} \; = \; V_l ,
  x_1 g_{1l} + x_2 g_{2l} + \ldots + x_N g_{Nl} \; = \; V_l ,
\end{equation}
where this equation is the result of an MV %multiplication;
product;
the operation, which is equivalent to $2N^2$ operations, is performed  in time $t_{meas}$. With increasing arrays size, this analog operation will become superior to that of the digital execution $2N^2/O_D=t_{dgig} > t_{meas}$ where $O_D$ is the operations per second of the digital engine. %\Anshul{Can remove next sentence if we need more space} -- space in the supplement is not an issue
A similar comparison holds true for the MV
%multiplication
product
with the transposed matrix by simply interchanging the roles of the rows and columns.
% The
% %rank-$1$
% OP
% update is a bit more complicated, where the digital row $\delta$ and column vector $x$ are converted into a stream of bits with length $L$. These bit streams are simultaneously applied to the rows and columns such that, wherever two pulses coincide, an incremental change of the conductivity is obtained which is modulated by the material and the operation condition dependent learning rate $\eta$ according to $g_{kl} =g_{kl} -\eta x_k  \delta_l$.

% \begin{figure}[htb]
% \centering    
% \includegraphics[width=0.95\textwidth]{cartoon3b.pdf}
% \caption{\it Hybrid architecture consisting of a digital chip including a (multi-core) CPU and 
% RAM system memory, connected to one or more analog 
%  arrays through a dedicated system bus. The CPU is focused on executing low-complexity 
% operations as well as transferring data to the analog  array(s) through the dedicated 
% system bus. } \label{fig:schematic}
% \end{figure}

The OP update is a bit more complicated \cite{gokmen2016acceleration}.
Here the digital row $\delta$ and column vector $x$ are converted into a stream of bits with length $L$. These bit streams are simultaneously applied to the rows and columns such that, wherever two pulses coincide, an incremental change of the conductivity is obtained which is modulated by the material resulting in 
\begin{equation}\label{app:na1b}
g_{kl} =g_{kl} + x_k  \delta_l.
\end{equation}
This
%rank-$1$
OP
update is performed in time $t_{update}=LN_{bit}t_{bit}$,
where $t_{bit}$ is the length of an individual bit in the bit stream of 
length $L$. The same argument as above for MV 
%multiplication
products
holds for the OP update, namely with increasing array size we see an increased benefit. A viable implementation of both the MV and OP operations need to take into account the effects of many sources of noise, such as 
digital-to-analog conversion,
analog-to-digital conversion, %programming noise to write the conductances and operating noise and non-idealities of the peripheral circuits and cross point elements.
and both the operational and non-ideal behaviors of the devices (i.e., peripheral circuits and cross-point elements).
%All these are currently under investigation in the literature and rapid progress has been made over the last years.
Rapid progress has been made on these issues over the past few years.

Each cross-point element in MFMWA arrays can be viewed as a processing unit as each one is capable of performing a single multiplication and addition operation needed for the %vector-matrix
MV
%multiply
products
and
%rank-$1$ outer-product
OP
updates. Therefore, larger arrays imply more compute units and hence increases the total throughput of the system. However, it is not possible to design arbitrarily large arrays due to physical constraints, such as increasing voltage drops and signal distortion issues and limited resistance range on devices. There are a few studies  
\cite{IBMJRD2019,gokmen2016acceleration, marinella2018} suggesting that it will be possible to design arrays with realistic technology constraints at a throughput of a single array around $400$ TOp/s while consuming less than 
$3$mm$^2$ chip area and about $2$W of energy. Including the peripheral circuit costs, the energy efficiency is therefore estimated as $200$ TOp/J. These energy and computation metrics are orders of magnitude better than the alternative digital solutions performing computations using high precision floating point numbers with $32$- or $64$-bit resolution. In order to improve the energy and computational metrics for deep neural network (DNN) workloads, highly optimized digital application-specific integrated circuit (ASIC) implementations have been proposed and developed~\cite{bfloat16,NvidiaA100,graphcore,Sun2019,Sun2019a}.
%[Ref's by Tayfun: 1)Cloud Tpu, . (2007). Available: https://cloud.google.com/tpu/docs/bfloat16. 2) Nvidia, . (2021). Available: https://www.nvidia.com/en-us/data-center/a100/. 3) Graphcore, . (2021). Available: https://www.graphcore.ai/. 4) Sun, X., Choi, J., Chen, C.-Y., Wang, N., Venkataramani, S., Srinivasan, V., et al. (2019). Hybrid 8-bit Floating point (HFP8) Training and Inference for Deep Neural Networks. Adv. Neural Inf. Process. Syst. 32, 4901–4910. 5) Sun, X., Wang, N., Chen, C.-Y., Ni, J.-M., Agrawal, A., Cui, X., et al. (2020). Ultra-Low Precision 4-bit Training of Deep Neural Networks. Adv. Neural Inf. Process. Syst. 33, 1796–1807.]. 
Using advanced methods it is shown that even $2$-bit integer arithmetic can be used for DNN inference workloads~\cite{Choi2019}.
%[Ref by Tayfun: https://mlsys.org/Conferences/2019/doc/2019/168.pdf]. 

On the other hand for training workloads, the current state of the art is to use hybrid $8$-bit floating point arithmetic, which is significantly more costly compared to inference solutions. These examples illustrate that depending on the task in hand, the benefit of using analog may diminish if there exists a digital solution running at aggressively reduced precision as in the case of DNN inference. However, algorithms such as those detailed in this manuscript use the same matrix primitives but require relatively high numeric precision for convergence. For these set of algorithms, MFMWA hardware still provides significant energy and computational benefits and hence can be used as a building block for general purpose hardware that addresses algorithms and tasks well beyond deep learning.

% \MSS{Malte to draft high-level description of how MV and OP are executed on analog, based in part on the description in the simulation section, with appropriate references.}

%Matrix-vector
MV
products can be achieved on MFMWA hardware by storing the elements of the matrix into resistive device elements of an MFMWA array~\cite{sebastian2020memory}. The input vector is typically pulse-width modulated, so that longer pulses code for larger input values.
As illustrated in
Figure~\ref{revised-figna1}(a),
% Figure~1(a),
the voltage pulses result in a current at the output and coductances on the devices get multiplied and added up by means of Ohm's and Kirchoff's laws, respectively. The output current is proportional to the row-matrix and input-vector scalar product and multiple columns of the array can integrate the $n$ currents necessary for the  MV product in parallel. The current is then read out and converted to digital integer numbers in parallel using a dedicated
analog-digital converter (ADC) 
for each output line. Negative and positive input values and weight values can be handled in different phases.
%Outer product
OP
updates can be performed on an MFMWA array as follows. Since for independent stochastic processes (i.e., Poisson process) the coincidence probability in a finite time window is given by the product of the individual pulse probabilities, one can compute the OP update by encoding the input values into the rate of stochastic voltage pulse trains. The voltage pulse trains are then delivered at either ends of the array, while requiring that only coinciding voltage pulses of the same polarity change the device conductance by a small finite amount; see
Figure~\ref{revised-figna1}(a).
% Figure~1(a).
In this manner, the consecutive coincident pulses will add to the devices in parallel by an overall amount that is proportional to the product of the pulse probabilities of the input values, and therefore implementing the OP update operation in MFMWA memory. Negative values again can be handled in different updated phases. For more details refer to \cite{gokmen2016acceleration}.

\subsection{Hybrid MFMWA Architecture Technology}
\label{app:ssec:stochastic}
To reduce power and increase performance, a common approach is to use computing substrates that operate at low precision. Examples include binary neural networks (BNN) with single bit weights \cite{Qin2020}, stochastic computing \cite{Alaghi2013} and dither computing \cite{dithercomputing2021} where data are represented as single-bit pulses. 
As described in the main body of the paper, 
stochastic computing represents and processes information as a series of binary pulses which allows simple logic gates to perform computations such as addition and multiplication. 
While this coding format is used in analog computing to provide the input data to the MFMWA array during the OP operation, the decoding of the stochastic computing output is done via the physics of each array element \cite{gokmen2016acceleration}.

A closely related concept is stochastic rounding \cite{connolly2021} where a value is quantized to one of the two nearest quantization levels with probability inversely proportional to the distance to these levels. This technique is used to develop Deep Learning hardware with low precision arithmetic \cite{wang:2018}. It has been shown \cite{Gupta2015a} that randomization in the form of stochastic rounding in Deep Learning and in stochastic computing can compensate, to some degree, the lack of precision. In digital technology this necessitates the inclusion of costly pseudo-random number generators.
For analog computing, however, in contrast to digital computing with limited precision, additional randomization units are not needed as the randomization is inherent in the technology in the form of read noise and write noise. We therefore expect the benefits of stochastic computing and stochastic rounding for digital technologies to occur naturally in MFMWA technologies without the need for additional hardware.  

Another common method to exploit randomness to trade off performance versus precision and accuracy is to only process a random subset of data/computations under the assumption that suitably selected ``average'' behavior is close to the ``optimal'' behavior. Such randomized algorithms \cite{motwani_raghavan_1995} are the workhorse of modern high performance computing where the amount of data and/or computations are increasing faster than the improvement in digital computer performance.
As outlined in
% %Sections~\ref{sec:algs}
% Sections~{\color{red} 3.2} 
% and \ref{app:ssec:algs}, 
the main body of the paper,
analog computing can exploit such algorithms and reap the benefits of a reduced number of operations while retaining the advantages of high speed and parallelism.

% In stochastic computing, information is represented and processed as a series of binary pulses. This allows simple arithmetic units for computation such as addition and multiplication, and additionally provides error tolerance. Stochastic computing is also used in analog computing to provide the input data to the analog array.
% The decoding of the stochastic computing output, done via a summation of the pulses in digital technology, is performed in the analog array via the inherent physics of output capacitors whose voltage is proportional to the integral of the current (i.e., the accumulated charge). 

% \alert{Chai Wah:  please add one more paragraph that ties analog and randomized algorithms as a transition to the next section}

\subsection{Computational Model for MFMWA Accelerator Devices}
\label{app:ssec:model}

In this section, we consider a computational model that supports an analysis of the time complexity and energy complexity of the core primitives of our algorithms executed on both a hybrid MFMWA accelerator architecture and a purely digital accelerator architecture.
More specifically, the algorithms proposed in this paper are based on a computational machine equipped with a hybrid MFMWA accelerator comprising a large digital main-memory and a collection of MFMWA
%accelerator
devices each consisting of multiple tiles of memresistive arrays,
additional digital storage, logic to quickly convert between digital and analog representations (ADC and DAC), and additional digital logic to support the matrix operations supplied by the MFMWA accelerator (e.g., comparators to find the maximum matrix entry in order to normalize inputs, and logic to combine results across mutltiple tiles);
refer to
Figure~\ref{revised-figna1}(b)
% Figure~1(b)
for a schematic abstraction of this architecture (see also \cite{le2018mixed}).
Such a hybrid MFMWA accelerator
%device
is designed towards performing certain matrix algebra computations quickly.
In particular, the corresponding MFMWA devices support the in-memory computation of MV products and OP updates;
%, where these computations are performed across the tiles in parallel.
matrices are stored across tiles of MFMWA arrays,
%each acting as a weight-stationary systolic array,
whereas vectors are primarily stored in the digital cache-memory of the MFMWA accelerator.
%which includes logic for very fast conversion to/from analog representations.

For comparison with a computational model of a corresponding purely digital accelerator architecture (defined below), 
we consider a computational model of the aforementioned hybrid MFMWA accelerator architecture
%with analog
%%accelerator
%devices
that
%comprises
includes
$m$ tiles, each consisting of an $n \times n$ MFMWA array, and digital cache-memory arrays which are in total on the order of $n \times m$.
This computational model assumes an architecture in which access to and computation on the $m$ tiles is performed in parallel, and thus the time complexity for accessing the $m$ tiles and performing computational operations on the $m$ tiles does not grow with $m$.
On the other hand, the energy complexity does grow with $m$ because the level of energy expended is basically the same irrespective of whether the $m$ tiles operate sequentially or in parallel.
Moreover, the time complexity for combining the results from across each of the $m$ tiles grows logarithmically with $m$ (see below).
%We focus on the steps of the core primitives that have different time and energy complexities within the computational models for the hybrid-analog and purely-digital accelerators.
The values of $m$ and $n$ will depend on various design constraints and tradeoffs involving the technology, physics and engineering used to implement the MFMWA devices.

We first focus on the computation of an MV product between an $(n \sqrt{m}) \times (n \sqrt{m})$ matrix and an $(n \sqrt{m}) \times 1$ vector, which we separate into an initial matrix writing phase and a computation phase.
The initial writing phase comprises the mapping of the $(n \sqrt{m}) \times (n \sqrt{m})$ input matrix from the digital main-memory of the MFMWA accelerator to the arrays of the analog devices.
This consists of $n$ sequential mappings of order-$n$ column vectors across the $m$ tiles having a time complexity of $\mathcal{T}_W \cdot n$ and an energy complexity of $\mathcal{E}_W \cdot n m$, where $\mathcal{T}_W$ and $\mathcal{E}_W$ respectively denote the time and energy required to map a single order-$n$ column vector onto an array.
Hence, the total time complexity and total energy complexity of the initial matrix writing phase is given by
\begin{equation}
    T_{\mathcal{MW}} = \mathcal{T}_W \cdot n
    \qquad \mbox{  and  } \qquad
    E_{\mathcal{MW}} = \mathcal{E}_W \cdot n m ,
    \label{eq:WR-T+E}
\end{equation}
after which the $(n \sqrt{m}) \times (n \sqrt{m})$ input matrix will be mapped onto the arrays of the MFMWA devices.

The next step is to transport the $(n \sqrt{m}) \times 1$ input vector from the digital main-memory of the accelerator to the MFMWA devices, consisting of a single transport of order-$n$ column vectors across the $m$ tiles.
This step has a time complexity of $\mathcal{T}_I$ and an energy complexity of $\mathcal{E}_I \cdot m$, where $\mathcal{T}_I$ and $\mathcal{E}_I$ respectively denote the time and energy required to transport a single order-$n$ column vector onto a digital cache-memory array.
Then the actual MV product is performed in parallel across the $m$ MFMWA tiles having a time complexity of $\mathcal{T}_M$ and an energy complexity of $\mathcal{E}_M \cdot m$,
where $\mathcal{T}_M$ and $\mathcal{E}_M$ respectively denote the time required and the energy expended to perform the parallel operation on a single tile.
The last step is to combine and transport the $n \times 1$ result vectors from each of the $m$ MFMWA tiles to obtain the final order-$(n \sqrt{m}) \times 1$ result vector in a digital cache-memory array,
consisting of logarithmic pairwise combinations and transports of the result vectors across the $m$ tiles in parallel.
This step has a time complexity of $\mathcal{T}_R \cdot \log_2 m$ and an energy complexity of $\mathcal{E}_R \cdot m$, where $\mathcal{T}_R$ and $\mathcal{E}_R$ respectively denote
the time and energy required to combine and transport each pair of result vectors from an MFMWA tile to a corresponding digital cache-memory array.

Hence, once the input matrix has been mapped onto the arrays of the MFMWA devices, the time complexity of an MV product between an $(n \sqrt{m}) \times (n \sqrt{m})$ matrix and an $(n \sqrt{m}) \times 1$ vector is given by
\begin{equation}
T_{\mathcal{MV}} =
\mathcal{T}_I +
\mathcal{T}_M +
\mathcal{T}_R \cdot \log_2 m ,
\label{eq:MV-T}
\end{equation}
and the corresponding energy complexity of the MV product is given by
\begin{equation}
E_{\mathcal{MV}} =
\mathcal{E}_I \cdot m +
\mathcal{E}_M \cdot m +
\mathcal{E}_R \cdot m .
\label{eq:MV-E}
\end{equation}
%An important assumption above is that the time and energy required to read the elements of the result vector from the co-processor components to the main digital host system is either
%negligble or simply the same cost is incurred on the corresponding digital system.
This complexity analysis assumes that the $(n \sqrt{m}) \times (n \sqrt{m})$ input matrix has already been mapped from the digital main-memory of the accelerator onto the arrays of the MFMWA devices, either after the initial writing phase or as the result of a previous operation (see OP update below).
However, if the input matrix does not already reside on the arrays, then the total time complexity and total energy complexity for the MV product would be given by
%$T_{\mathcal{MW}} + T_{\mathcal{MV}}$ and $E_{\mathcal{MW}} + E_{\mathcal{MV}}$, 
$$T_{\mathcal{MW}} + T_{\mathcal{MV}} \qquad \mbox{  and  } \qquad E_{\mathcal{MW}} + E_{\mathcal{MV}},$$ 
respectively.
We note that the values of the variables $\mathcal{T}_W$, $\mathcal{E}_W$, $\mathcal{T}_I$, $\mathcal{E}_I$, $\mathcal{T}_M$, $\mathcal{E}_M$, $\mathcal{T}_R$ and $\mathcal{E}_R$
depend on the values of $m$ and $n$ and on the related technology, physics and engineering used to implement the MFMWA devices.
Representative values for these variables will be discussed below (see Table~\ref{tab:variable-values}).

We now focus on the computation of an OP update between an order-$(n \sqrt{m})$ column vector and an order-$(n \sqrt{m})$ row vector.
Some algorithms might require an initialization of the array to certain desired values, such as setting the elements of the array to zero, in which case the time and energy complexity of this initialization step is as given in~\eqref{eq:WR-T+E}.
Then, the $(n \sqrt{m}) \times 1$ input vector is transported from the digital main-memory of the accelerator to the MFMWA devices.
This consists of a single mapping of order-$n$ column vectors across the $m$ tiles having a time complexity of $\mathcal{T}_I$ and an energy complexity of $\mathcal{E}_I \cdot m$,
where $\mathcal{T}_I$ and $\mathcal{E}_I$
%again respectively denote the time and energy required to transport a single order-$n$ column vector onto a digital cache-memory array.
are as defined above.
The next step is to transport the $1 \times (n \sqrt{m})$ input vector from the digital main-memory of the accelerator to the MFMWA devices, consisting of a single transport of order-$n$ row vectors across the $m$ tiles.
This step 
again
has a time complexity of $\mathcal{T}_I$ and an energy complexity of $\mathcal{E}_I \cdot m$.
%, where $\mathcal{T}_I$ and $\mathcal{E}_I$ again respectively denote the time and energy required to transport a single order-$n$ row vector onto a digital cache-memory array.
Then the actual OP update is performed in parallel across the $m$ MFMWA tiles having a time complexity of $\mathcal{T}_O$ and an energy complexity of $\mathcal{E}_O \cdot m$,
where $\mathcal{T}_O$ and $\mathcal{E}_O$ respectively denote the time required and the energy expended to perform the parallel operation on a single tile.
At this point, the resulting matrix is left in the arrays of the MFMWA tiles for use in subsequent computations.
In particular, if the OP update is followed by an MV product, then
%there would be a savings of $\mathcal{T}_W \cdot n$ and $\mathcal{E}_W \cdot n m$, and
the total time complexity and total the energy complexity of the subsequent MV product operation in this case would be as given in \eqref{eq:MV-T} and \eqref{eq:MV-E}, respectively.

Hence, once the arrays of the MFMWA devices have been appropriately initialized as needed, the time complexity of an OP update between an order-$(n \sqrt{m})$ column vector and an order-$(n \sqrt{m})$ row vector is given by
\begin{equation}
T_{\mathcal{OP}} =
\mathcal{T}_I \cdot 2 +
\mathcal{T}_O ,
\label{eq:OP-T}
\end{equation}
and the corresponding energy complexity of the OP update is given by
\begin{equation}
E_{\mathcal{OP}} =
\mathcal{E}_I \cdot 2 m +
\mathcal{E}_O \cdot m .
\label{eq:OP-E}
\end{equation}
This complexity analysis assumes that the order-$(n \sqrt{m})$ input vectors need to be transported from the digital main-memory of the accelerator to the digital cache-memory arrays of the MFMWA devices.
However, if the input vectors already reside on the MFMWA devices (either as part of or resulting from a previous operation), then the above complexity analysis would simplify as follows
%$T_{\mathcal{OP}} = \mathcal{T}_O$ and $E_{\mathcal{OP}} = \mathcal{E}_O \cdot m$.
$$T_{\mathcal{OP}} = \mathcal{T}_O \qquad \mbox{  and  } \qquad E_{\mathcal{OP}} = \mathcal{E}_O \cdot m.$$
The above complexity analysis further assumes that the arrays either do not need to be initialized or have already been appropriately initialized for the algorithm of interest.
However, if the OP update requires an initialization of the array to certain desired values, then the total time complexity and total energy complexity for the OP update would be given by
%$T_{\mathcal{MW}} + T_{\mathcal{OP}}$ and $E_{\mathcal{MW}} + E_{\mathcal{OP}}$, 
$$T_{\mathcal{MW}} + T_{\mathcal{OP}} \qquad \mbox{  and  } \qquad E_{\mathcal{MW}} + E_{\mathcal{OP}},$$ 
respectively.
Once again, we note that the values of the variables $\mathcal{T}_I$, $\mathcal{E}_I$, $\mathcal{T}_O$ and $\mathcal{E}_O$
depend on the values of $m$ and $n$ and on the related technology, physics and engineering used to implement the MFMWA devices.
Representative values for these variables will be discussed below (see Table~\ref{tab:variable-values}).

Our above complexity analysis given in 
% \eqref{eq:MV-T}, \eqref{eq:MV-E}, \eqref{eq:OP-T} and \eqref{eq:OP-E}
\eqref{eq:MV-T}~--~\eqref{eq:OP-E}
does not consider the time complexity and energy complexity for the return transfer of the final results from the MFMWA devices back to the digital main-memory of the accelerator.
In the case of an MV product of an $(n \sqrt{m}) \times (n \sqrt{m})$ matrix and an $(n \sqrt{m}) \times 1$ vector, the return transfer step involves transporting the $(n \sqrt{m}) \times 1$ final result vector from the
digital cache-memory array
%analog devices
to the digital main-memory of the accelerator.
This vector read consists of a single transport of order-$n$ column vectors across the $m$ tiles having a time complexity of $T_{\mathcal{RV}} = \mathcal{T}_I$ and an energy complexity of $E_{\mathcal{RV}} = \mathcal{E}_I \cdot m$, where $\mathcal{T}_I$ and $\mathcal{E}_I$ 
are as defined above.
%respectively denote the time and energy required to transport a single order-$n$ column vector from a digital cache-memory array back to the digital main-memory of the analog accelerator.
In the case of an OP update of an order-$(n \sqrt{m})$ column vector and an order-$(n \sqrt{m})$ row vector, the return transfer step consists of transporting the $(n \sqrt{m}) \times (n \sqrt{m})$ result matrix from the MFMWA devices to the digital main-memory of the accelerator.
This matrix read is performed through a sequence of $(n \sqrt{m})$ MV products, each producing one $(n \sqrt{m}) \times 1$ column of the result matrix which is then transported from the
%analog devices
digital cache-memory array
to the digital main-memory of the accelerator.
The total time complexity and total energy complexity is then given by
$$T_{\mathcal{RM}} = (n \sqrt{m}) (T_{\mathcal{MV}} + \mathcal{T}_I) \qquad \mbox{ and } \qquad E_{\mathcal{RM}} = (n \sqrt{m}) (E_{\mathcal{MV}} + \mathcal{E}_I 
\cdot m),$$ respectively.
A complete summary of the time and energy complexities for the core primitives of our algorithms is presented in Table~\ref{tab:time-energy}.

\begin{table}[htb]
\caption{\it Summary of time and energy complexities of the hybrid MFMWA accelerator architecture for core primitives of our algorithms with respect to an MV product of an $(n \sqrt{m}) \times (n \sqrt{m})$ matrix and an $(n \sqrt{m}) \times 1$ vector, and an OP update of an order-$(n \sqrt{m})$ column vector and an order-$(n \sqrt{m})$ row vector. The MFMWA accelerator
%architecture 
%%consists of
includes
$m$ tiles of $n\times n$ of memresistive arrays, additional digital storage on the order of $n\times m$, ADC and DAC logic, and additional digital logic to support the MV product and OP update operations. The matrix write primitive is from the digital main-memory to the arrays, the vector read primitive is from the digital cache-memory to the digital main-memory, and the matrix read primitive is from the arrays to the digital main-memory.}\label{tab:time-energy}
%{\small
\begin{center} 
\begin{tabular}{||c|c|c||}
\hline\hline
\textbf{Core Primitive} & \textbf{Time Complexity} & \textbf{Energy Complexity} \\
\hline
Matrix Write &
%\begin{equation} T_{\mathcal{MW}} = \mathcal{T}_W \cdot n \quad \tag{\ref{eq:WR-T+E}} \end{equation}
$T_{\mathcal{MW}} = \mathcal{T}_W \cdot n$
&
%\begin{equation} E_{\mathcal{MW}} = \mathcal{E}_W \cdot n m \quad \tag{\ref{eq:WR-T+E}} \end{equation} 
$E_{\mathcal{MW}} = \mathcal{E}_W \cdot n m$
\\
\hline
MV Product &
%\begin{equation} T_{\mathcal{MV}} = \mathcal{T}_I + \mathcal{T}_M + \mathcal{T}_R \cdot \log_2 m \quad \tag{\ref{eq:MV-T}} \end{equation}
$T_{\mathcal{MV}} = \mathcal{T}_I + \mathcal{T}_M + \mathcal{T}_R \cdot \log_2 m$
&
%\begin{equation} E_{\mathcal{MV}} = \mathcal{E}_I \cdot m + \mathcal{E}_M \cdot m + \mathcal{E}_R \cdot m \quad \tag{\ref{eq:MV-E}} \end{equation}
$E_{\mathcal{MV}} = \mathcal{E}_I \cdot m + \mathcal{E}_M \cdot m + \mathcal{E}_R \cdot m$
\\
\hline
% Zero-Out Matrix &
% $T_{\mathcal{MZ}} = \mathcal{T}_Z$ &
% $E_{\mathcal{MZ}} = \mathcal{E}_Z \cdot m$ \\
% \hline
OP Update &
%\begin{equation} T_{\mathcal{OP}} = \mathcal{T}_I \cdot 2 + \mathcal{T}_O \quad \tag{\ref{eq:OP-T}} \end{equation}
$T_{\mathcal{OP}} = \mathcal{T}_I \cdot 2 + \mathcal{T}_O$
&
%\begin{equation} E_{\mathcal{OP}} = \mathcal{E}_I \cdot 2 m + \mathcal{E}_O \cdot m \quad \tag{\ref{eq:OP-E}} \end{equation} 
$E_{\mathcal{OP}} = \mathcal{E}_I \cdot 2 m + \mathcal{E}_O \cdot m$
\\
\hline
Vector Read &
$T_{\mathcal{RV}} = \mathcal{T}_I$ &
$E_{\mathcal{RV}} = \mathcal{E}_I \cdot m$ \\
\hline
Matrix Read &
$T_{\mathcal{RM}} = (n \sqrt{m}) (T_{\mathcal{MV}} + \mathcal{T}_I)$ &
$E_{\mathcal{RM}} = (n \sqrt{m}) (E_{\mathcal{MV}} + \mathcal{E}_I 
\cdot m)$ \\
\hline\hline
\end{tabular}
\end{center} 
\end{table}

Next, we seek to compare the above time complexity and energy complexity analysis against a computational model for a corresponding purely digital
%solution based on a computational machine 
architecture
equipped with a GPU accelerator comprising a large digital main-memory, a collection of hundreds of processing cores, additional digital storage, and additional digital logic to support the various operations supplied by the digital accelerator.
For the same area efficiency as such a single digital GPU accelerator, we consider the foregoing computational model of a hybrid MFMWA accelerator architecture with $m = 64$ and $n = 2K$, where $K := 1024$.
This means that our comparisons will be based on an MV product between a $16K \times 16K$ matrix and a $16K \times 1$ vector, and based on an OP update between an order-$16K$ column vector and an order-$16K$ row vector.
For these core operations on a baseline design and technology of MFMWA
%accelerator
devices, the time and energy complexity  will be as given in Table~\ref{tab:time-energy} with respect to the variable ranges and values given in Table~\ref{tab:variable-values}.
This then yields a matrix write time and energy of
\begin{equation}
    T_{\mathcal{MW}} \in [2048,20480] \mbox{ microseconds} \quad \mbox{ and } \quad E_{\mathcal{MW}} \in [262144,13107200] \mbox{ micro-Joules}, \label{eq:MW-values}
\end{equation}
respectively;
a MV product time and energy of 
\begin{equation}
    T_{\mathcal{MV}} \in [0.135,0.240] \mbox{ microseconds} \qquad \mbox{ and } \qquad E_{\mathcal{MV}} \in [12.928,33.28] \mbox{ micro-Joules}, \label{eq:MV-values}
\end{equation}
respectively;
an OP update time and energy of 
\begin{equation}
    T_{\mathcal{OP}} \in [0.11,0.14] \mbox{ microseconds} \qquad \mbox{ and } \qquad E_{\mathcal{OP}} \in [12.928,33.28] \mbox{ micro-Joules}, \label{eq:OP-values}
\end{equation}
respectively;
a vector read time and energy of
\begin{equation}
    T_{\mathcal{RV}} \in [0.005,0.020] \mbox{ microseconds} \qquad \mbox{ and } \qquad E_{\mathcal{RV}} \in [64000,640000] \mbox{ micro-Joules}, \label{eq:RV-values}
\end{equation}
respectively;
and a matrix read time and energy of 
\begin{equation}
    T_{\mathcal{RM}} \in [2293.76,4259.84] \mbox{ microseconds} \quad \mbox{ and } \quad E_{\mathcal{RM}} \in [212860.928,555745.28] \mbox{ micro-Joules}, \label{eq:RM-values}
\end{equation}
respectively.

\begin{table}[htb]
\caption{\it Representative ranges and values for the variables of the time and energy complexity given in Table~\ref{tab:time-energy} for a baseline design and technology of a hybrid MFMWA accelerator architecture with $m=64$, $n=2K$, $K=1024$, and $M=1,048,576$. These ranges and values are based on the design and technology results presented in~\cite{gokmen2016acceleration}.}\label{tab:variable-values}
%{\small
\begin{center} 
\begin{tabular}{||c|c||}
\hline\hline
\textbf{Time Variables} & \textbf{Energy Variable} \\
\hline
$\mathcal{T}_W \in [1,10] \mbox{ microseconds}$
&
$\mathcal{E}_W \in [2,100] \mbox{ microJoules}$
\\
\hline
$\mathcal{T}_I \in [5,20] \mbox{ nanoseconds}$
&
$\mathcal{E}_I \in [1,10] \mbox{ nanoJoules}$
\\
\hline
$\mathcal{T}_M \approx 100 \mbox{ nanoseconds}$
&
$\mathcal{E}_M \in [200,500] \mbox{ nanoJoules}$
\\
\hline
$\mathcal{T}_R \in [5,20] \mbox{ nanoseconds}$
&
$\mathcal{E}_R \in [1,10] \mbox{ nanoJoules}$
\\
\hline
$\mathcal{T}_O \approx 100 \mbox{ nanoseconds}$
&
$\mathcal{E}_O \in [200,500] \mbox{ nanoJooules}$
\\
\hline\hline
\end{tabular}
\end{center} 
\end{table}

Now, for comparison with the above time complexity and energy complexity of the hybrid MFMWA accelerator architecture, we consider the foregoing purely digital accelerator architecture based on an A100 GPU with $40$-$80$MB digital cache-memory~\cite{NvidiaA100}.
Assuming $16$-bit precision\footnote{We note that the standard computational codes for many of the applications in
% %Sections~\ref{sec:exp}
% Sections~{\color{red} 4}
% and \ref{app:ssec:exp}
the Experimental Results section of the main body of the paper and
% Section~\ref{app:ssec:exp}
Appendix~\ref{app:ssec:exp}
either use $32$-bit or $64$-bit precision, which would further increase the time and energy complexity benefits of the hybrid MFMWA accelerator architecture.}, the computation of the core primitives will be memory-bound (or bandwidth-bound) for both an MV product between a $16K \times 16K$ matrix ($2 \times 16K \times 16K$ bytes) and a $16K \times 1$ vector ($2 \times 16K \times 1$ bytes), and an OP update between an order-$16K$ column vector ($2 \times 16K$ bytes) and an order-$16K$ row vector ($2 \times 16K$ bytes).
Then, in comparison with the time and energy complexity analysis of the hybrid MFMWA accelerator architecture, the purely digital accelerator architecture has a write time and energy of
\begin{equation}
    T_{\mathcal{MW}}^d \approx 250 \mbox{ microseconds} \qquad \mbox{ and } \qquad E_{\mathcal{MW}}^d \approx 12000 \mbox{ micro-Joules}, \label{eq:MWd-values}
\end{equation}
respectively;
a MV product time and energy of 
\begin{equation}
    T_{\mathcal{MV}}^d \in [250.005,250.02] \mbox{ microseconds} \quad \mbox{ and } \quad E_{\mathcal{MV}}^d \in [12000.064,12000.64] \mbox{ micro-Joules}, \label{eq:MVd-values}
\end{equation}
respectively,
which also includes the comparable time (with respect to $\mathcal{T}_I$) and energy (with respect to $\mathcal{E}_I$) to transport the input vector from the digital main-memory of the accelerator to the digital cache-memory array;
an OP update time and energy of 
\begin{equation}
    T_{\mathcal{OP}}^d \in [250.01,250.04] \mbox{ microseconds} \quad \mbox{ and } \quad E_{\mathcal{OP}}^d \in [12000.128,12001.28] \mbox{ micro-Joules}, \label{eq:OPd-values}
\end{equation}
respectively,
which also includes the comparable time (with respect to $\mathcal{T}_I$) and energy (with respect to $\mathcal{E}_I$) to transport both input vectors from the digital main-memory of the accelerator to the digital cache-memory array;
% a vector read time and energy of
% \begin{equation}
%     T_{\mathcal{RV}}^d \approx 0.0153 \mbox{ microseconds} \qquad \mbox{ and } \qquad E_{\mathcal{RV}}^d \approx 0.7324 \mbox{ micro-Joules}, \label{eq:RVd-values}
% \end{equation}
% respectively;
and a matrix read time and energy of 
\begin{equation}
    T_{\mathcal{RM}}^d \approx 250 \mbox{ microseconds} \qquad \mbox{ and } \qquad E_{\mathcal{RM}}^d \approx 12000 \mbox{ micro-Joules}, \label{eq:RMd-values}
\end{equation}
respectively.
The vector read time and energy on the purely digital accelerator should be comparable to the vector read time and energy on the hybrid MFMWA accelerator, i.e., $T_{\mathcal{RV}}^d \approx T_{\mathcal{RV}}$ and  $E_{\mathcal{RV}}^d \approx E_{\mathcal{RV}}$,
and therefore we do not consider further these corresponding time and energy costs.

From equations~\eqref{eq:MW-values}~--~\eqref{eq:RM-values} and %equations~
\eqref{eq:MWd-values}~--~\eqref{eq:RMd-values}, we summarize in
Table~\ref{tab:analog-digital}
%Table~{\color{red} 1} 
of the main paper
our time complexity and energy complexity comparison between the hybrid MFMWA accelerator architecture and the purely digital accelerator architecture.
While \eqref{eq:MW-values},\eqref{eq:MWd-values} and \eqref{eq:RM-values},\eqref{eq:RMd-values} respectively show up to two orders of magnitude larger time complexity and up to three orders of magnitude larger energy complexity for the matrix write (primarily due to mapping the matrix onto the arrays) and matrix read (primarily due to the $(n\sqrt{m})$ MV products) on the MFMWA accelerator, equations \eqref{eq:MV-values},\eqref{eq:MVd-values} and \eqref{eq:OP-values},\eqref{eq:OPd-values} respectively show three orders of magnitude larger time complexity and energy complexity for an MV product and OP update on the digital accelerator.
In particular, once the $(n\sqrt{m}) \times (n\sqrt{m})$ input matrix has been mapped onto the arrays of the MFMWA devices, the time and energy complexity of an MV product is approximately $0.188\mu$s and $23.104\mu$J, respectively, on the hybrid MFMWA accelerator in comparison with approximately $250\mu$s and $12000\mu$J, respectively, on the purely digital accelerator.
Similarly, assuming there is no requirement to initialize the arrays of the MFMWA devices, the time and energy complexity of an OP update is approximately $0.125\mu$s and $23.104\mu$J, respectively, on the hybrid MFMWA accelerator in comparison with approximately $250\mu$s and $12000\mu$J, respectively, on the purely digital accelerator.

\begin{table}[htb]
\caption{\it Summary of time complexity (in microseconds, or $\mu$s) and energy complexity (in micro-Joules, or $\mu$J) comparisons between the hybrid MFMWA accelerator architecture and purely digital accelerator architecture for the core primitives of our algorithms with respect to an MV product of an $(n \sqrt{m}) \times (n \sqrt{m})$ matrix and an $(n \sqrt{m}) \times 1$ vector, and an OP update of an order-$(n \sqrt{m})$ column vector and an order-$(n \sqrt{m})$ row vector. The MFMWA accelerator %architecture 
%%consists of
includes
$64$ tiles of $2K \times 2K$ of memresistive arrays, additional digital storage,
%on the order of $n\times m$,
ADC and DAC logic, and additional digital logic to support the MV product and OP update operations;
the matrix write primitive is from the digital main-memory to the arrays, and the matrix read primitive is from the arrays to the digital main-memory. The digital accelerator
%architecture
is based on an A100 GPU with $40$-$80$MB digital cache-memory;
the matrix write primitive is from the digital main-memory to the digital cache-memory, and the matrix read primitive is from the digital cache-memory to the digital main-memory.}\label{tab:analog-digital:app}
%{\small
\begin{center} 
\begin{tabular}{||l|l|l||}
\hline\hline
\textbf{Core Primitive} & \textbf{Time Complexity} & \textbf{Energy Complexity} \\
\hline\hline
MFMWA Matrix Write &
$T_{\mathcal{MW}} \in [2048\mu\mbox{s},20480\mu\mbox{s}]$ 
&
$E_{\mathcal{MW}} \in [262144\mu\mbox{J},13107200\mu\mbox{J}]$
\\
\hline
Digital Matrix Write &
$T_{\mathcal{MW}}^d \approx 250\mu\mbox{s}$ 
&
$E_{\mathcal{MW}}^d \approx 12000\mu\mbox{J}$
\\
\hline\hline
MFMWA MV Product &
$T_{\mathcal{MV}} \in [0.135\mu\mbox{s},0.240\mu\mbox{s}]$
&
$E_{\mathcal{MV}} \in [12.928\mu\mbox{J},33.28\mu\mbox{J}]$
\\
\hline
Digital MV Product &
$T_{\mathcal{MV}}^d \in [250.005\mu\mbox{s},250.02\mu\mbox{s}]$
&
$E_{\mathcal{MV}}^d \in [12000.064\mu\mbox{J},12000.64\mu\mbox{J}]$
\\
\hline\hline
MFMWA OP Update &
$T_{\mathcal{OP}} \in [0.11\mu\mbox{s},0.14\mu\mbox{s}]$
&
$E_{\mathcal{OP}} \in [12.928\mu\mbox{J},33.28\mu\mbox{J}]$
\\
\hline
Digital OP Update &
$T_{\mathcal{OP}}^d \in [250.01\mu\mbox{s},250.04\mu\mbox{s}]$
&
$E_{\mathcal{OP}}^d \in [12000.128\mu\mbox{J},12001.28\mu\mbox{J}]$
% \\
% \hline\hline
% MFMWA Vector Read &
% $T_{\mathcal{RV}} \in [0.005\mu\mbox{s},0.020\mu\mbox{s}]$
% &
% $E_{\mathcal{RV}} \in [64000\mu\mbox{J},640000\mu\mbox{J}]$
% \\
% \hline
% Digital Vector Read &
% $T_{\mathcal{RV}}^d \approx 0.0153\mu\mbox{s}$
% &
% $E_{\mathcal{RV}}^d \approx 0.7324\mu\mbox{J}$
\\
\hline\hline
MFMWA Matrix Read &
$T_{\mathcal{RM}} \in [2293.76\mu\mbox{s},4259.84\mu\mbox{s}]$
&
$E_{\mathcal{RM}} \in [212860.928\mu\mbox{J},555745.28\mu\mbox{J}]$
\\
\hline
Digital Matrix Read &
$T_{\mathcal{RM}}^d \approx 250\mu\mbox{s}$
&
$E_{\mathcal{RM}}^d \approx 12000\mu\mbox{J}$
\\
\hline\hline
\end{tabular}
\end{center} 
\end{table}

For many algorithms involving OP updates, there is no need to initialize the arrays of the MFMWA devices;
and for many algorithms involving MV products, the input matrix may not need to be written to the arrays due to the result(s) of a previous operation(s), such as an OP update.
In these cases, we clearly observe significant time and energy improvements from the hybrid MFMWA accelerator architecture over the purely digital accelerator architecture for these two core primitives.
On the other hand, when an MV product requires the input matrix to be written to the arrays of the MFMWA devices or when an OP update requires the arrays to be initialized, then the time and energy complexity of these two operations is approximately $11264\mu$s and $6684695\mu$J, respectively, on the hybrid MFMWA accelerator in comparison with approximately $250\mu$s and $12000\mu$J, respectively, on the purely digital accelerator.
In such cases, the purely digital accelerator architecture provides significant time and energy improvements over the hybrid MFMWA accelerator architecture for a \emph{single execution} of either of the two core primitives.
However, for many algorithms of interest including those considered in this paper, the overhead of a matrix write or matrix initialization with respect to $T_{\mathcal{MW}}$ can be amortized by multiple subsequent MV products with the same matrix or multiple subsequent OP updates or a combination of both.
In particular, the significant time and energy benefits observed in the experimental results of
% %Sections~\ref{sec:exp}
% Sections~{\color{red} 4}
% and~\ref{app:ssec:exp} 
the main body of the paper
and
% Section~\ref{app:ssec:exp} 
Appendix~\ref{app:ssec:exp} 
are primarily due to multiple MV products and OP updates that overwhelmingly amortized the overhead related to $T_{\mathcal{MW}}$.
For example, the experiments based on streaming only employed OP updates with no matrix initialization, while the experiments based on PCA included an initial matrix write followed by combinations of MV products and OP updates (where an OP update was performed after every few MV products).

\subsection{Experimental Results}
\label{app:ssec:exp}
% %
% \textbf{Purpose}: present main experiments of key algorithms that highlight the significant performance (computation and/or energy) benefits of the analog computing technology over the existing digital computing technology
%
We first present additional details on our physical experiments.
Then we describe our simulation engine, followed by some additional details on our simulation experiments using this engine for randomized sketching and randomized PCA,
together with an additional set of experimental results for our randomized PCA.

% \textit{Have presentation of experimental results follow the organization of the key algorithms in Section~\ref{sec:algs}.}

\subsubsection{Physical Experiments}
We tested the randomized streaming technique described in
Algorithm~\ref{app:alg:algo0}
% Algorithm~$1$
on a real MFMWA array device. For these physical experiments, we use the same hardware designed to perform DNN training exploiting CMOS-based components where charge stored in a capacitor holds the state of each array element. The hardware, whose details can be found in \cite{kohda2020IEDM}, enables two fundamental operations for DNN training, i.e., MV
%multiplication
products
and OP updates, that are also needed to perform randomized streaming. This is the first time a randomized 
streaming algorithm is implemented on actual MFMWA hardware. The streaming algorithm 
computes the OP updates using the stochastic pulsing scheme proposed in \cite{gokmen2016acceleration}. Our experiment features for the first time 
this OP update mechanism implemented on real hardware for applications beyond deep learning training.

\begin{figure}
    \centering
    \includegraphics[width=0.450\textwidth]{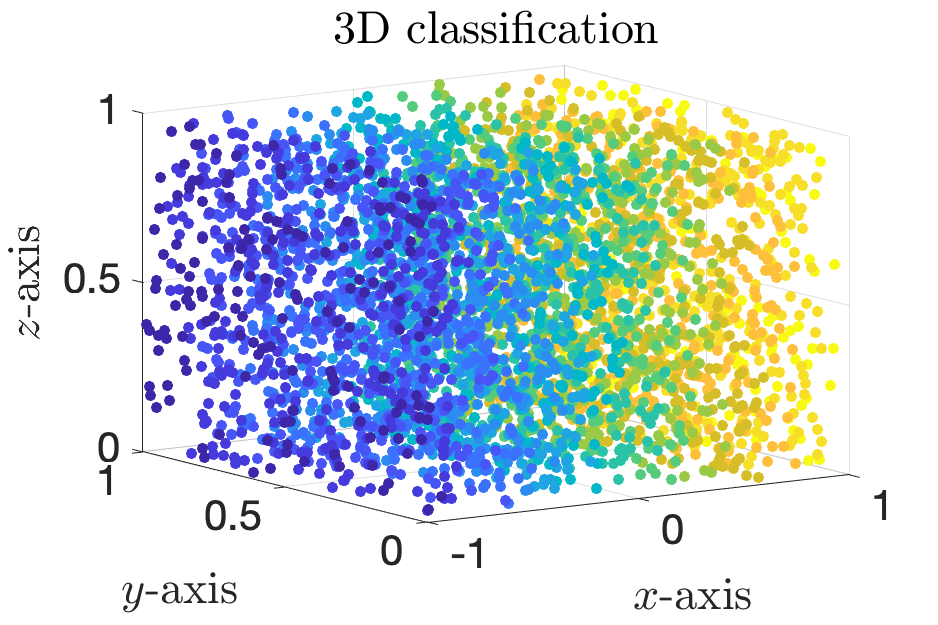}
    \includegraphics[width=0.450\textwidth]{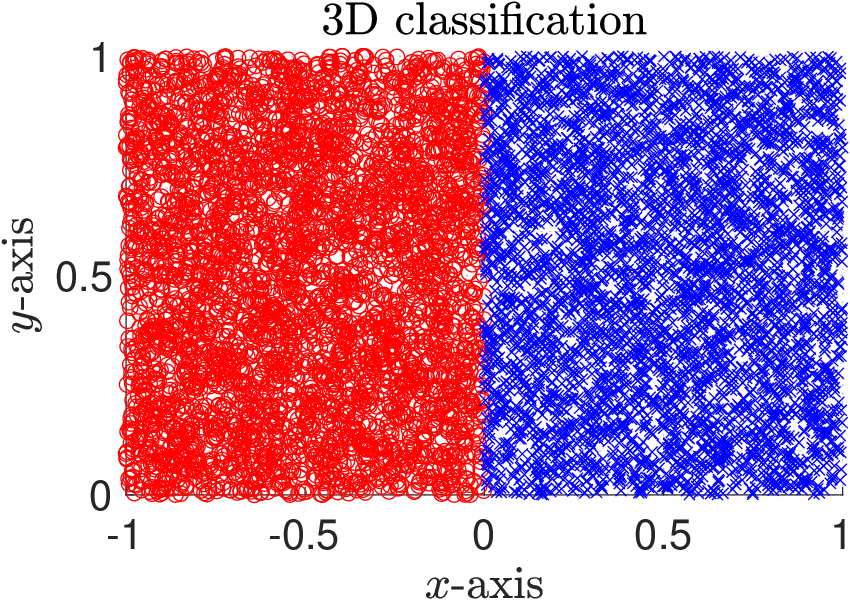}
    \caption{\it Dataset used for the physical experiment. Left: the dataset of three-dimensional coordinates. 
    Right: two-dimensional orthogonal projection along the 
two leading Cartesian coordinates. } \label{app:fig:physical1}
\end{figure}

{\bf Dataset:} The dataset for our experiments consists of $8192$ three-dimensional points distributed in the unit cube that is centered around the origin of a three-dimensional Cartesian coordinate system.
Half of these points were chosen such that their first coordinate is positive, and thus assigned the label ``$+1$''.
The remaining points were chosen so that their first coordinate is negative, and thus assigned the label ``$-1$''.
These $8192$ three-dimensional points were randomly divided equally in two disjoint sets, namely training and testing each of size $4096$.
A visualization of these three-dimensional points can be found in the left subfigure 
of Figure~\ref{app:fig:physical1}, whereas a two-dimensional orthogonal projection along the two leading Cartesian coordinates can be found in the right subfigure of Figure~\ref{app:fig:physical1}. 

{\bf Dimensionality reduction:} Using the above datasets, we construct a $4096\times 76$ matrix $A=[A^{(1)},\ldots,A^{(19)}]$, where $A^{(1)}=A^{(2)}=\ldots=A^{(19)}$ and $A^{(i)}=[x,y,z,t]$, where $x,y,z,t\in \mathbb{R}^{4096}$.
The vectors $x$, $y$, and $z$ respectively denote the first, second, and third coordinate of the $4096$ training samples, while $t$ denotes the corresponding label $\pm 1$. The sketching matrix $S$ was set equal to a $76\times 4096$ matrix 
with entries $\pm 1$. After initializing the conductance of each array at zero, 
we performed $4096$ OP updates where, during the $i$th update, the inputs 
to the word-lines and bit-lines of the MFMWA array were set equal to the 
$i$th row and column of the matrices $S$ and $A$, respectively. Three different experiments were performed, each with a different number of pulses used during the translation of inputs to the stochastic bit streams. As described in \cite{gokmen2016acceleration}, the number of pulses used during stochastic translation is a parameter that can be controlled by the user. 

{\bf Classification:} Once the update procedure is completed, we extract the resulting $76\times 76$ weight matrix $Z=[Z^{(1)},\ldots,Z^{(19)}]$ represented by the conductance of the arrays.
Then, for each $Z^{(i)}\in \mathbb{R}^{76\times 4}$, we solve the corresponding
ordinary linear least squares (OLLS)
problem $w_i=\argmin\limits_{w\in \mathbb{R}^3}\|Z^{(i)}_{:,1:3}w-Z^{(i)}_{:,4}\|_2$ on a digital computer. 
We next use a series of MV
%multiplications
products
performed using one-hot encoded inputs to extract the weights from the array.
The three-dimensional vectors $w_1,\ldots,w_{19}$ are the regressors used for classification. 
A global regressor $w$ is finally computed by averaging the regressors $w_1,\ldots,w_{19}$.
We then compute the dot product between $w$ and each three-dimensional 
sample of the test set. The label assigned to each sample is equal to the sign of this 
dot product (i.e., minus one if the dot product is negative, and plus one otherwise).
Three different experiments were conducted, each with a different number of pulse lengths to encode the inputs.

{\bf Results:}
Figure~\ref{fig:physical3} (left) 
% Figure~2 (left) 
plots the classification error rate achieved by the 
digital-based and MFMWA-based implementations of the randomized streaming 
algorithm versus that of the standard linear least squares approach based on 
normal equations (i.e., the baseline). Naturally, the latter method achieves 
the highest accuracy among the algorithms tested since it does not compress any 
information. On the other hand, streaming algorithms project the dataset onto a 
lower-dimensional subspace, thus leading to some loss of information. Notice that the
encoding through higher number of pulses leads to higher accuracy for the MFMWA-based 
implementation.

\begin{figure}[htb]
     \centering
     \includegraphics[width=0.950\textwidth]{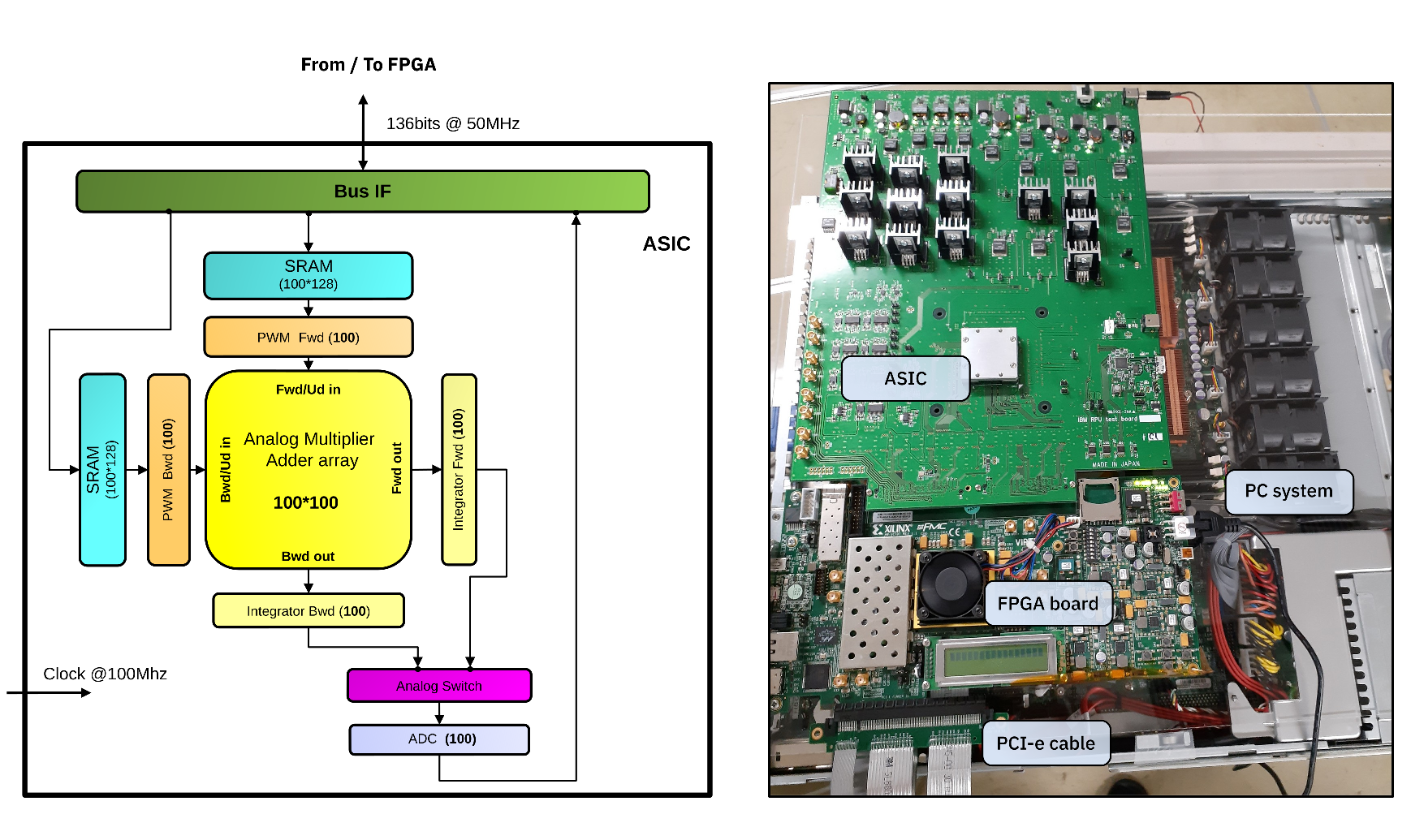}
     \caption{\it Experimental setup. Left panel shows the structure of the MFMWA chip used for the experiment. The MFMWA array using capacitors as weight elements is complimented with input drivers (pulse-width modulators) and output integration, as well as analog-to-digital converters. The chip is attached via a bus to FPGA and memory. The right panels shows the system level setup of the test board. } \label{app:fig:chip}
 \end{figure}

{\bf Experimental setup:}
For our physical experiment, we use a fully analog ASIC DNN training chip 
previously described in ~\cite{kohda2020IEDM}. 
% The 
% chip
% is designed and fabricated in CMOS 90nm technology using MIMCAP as capacitor weight. The unit cell area of the MIMCAP is 114/um$^2$ and the ASIC clock is running at 100 MHz. 
%From an architectural viewpoint, t
The chip consists of 3 independent
arrays, all of which use capacitor base unit cells to store weight
elements. The size of each array is $529\times 100,\ 100\times 100$, 
and $100\times 10$, respectively. Our experiment utilizes only the middle 
array (i.e., size $100\times 100$). A detailed structure of the array 
together with peripheral circuit (ASIC block diagram) is presented in 
Figure~\ref{app:fig:chip}. Our reduced analog chip has two sets of SRAM and pulse width 
modulation (PWM) to provide pulses to the analog multiplier/adder elements. In 
addition, our chip includes two integrators and ADC elements to convert total 
capacitor voltages to digital values.

Our randomized sketching experiment relies on the computation of OP 
updates. Input vectors for the OP computation were first stored on 
the FPGA portion of the chip, and pulse patterns were then generated for both input 
vectors for each OP update operation. The pulse pattern of each input value
was chosen to be random and binary, in such a way that the floating
point value of the input vector was coded into the probability of
having a pulse in one of 63 slots (of the total 127 slots with signs
and zero).  The two binary pattern were applied to column and row
drivers on the west and to the north, respectively. Voltage pulse
coincidence increased or decreased the charge of the capacitor at each
cross point of the array (in correspondence with the sign of
the voltage input). A single pulse was 10$ns$ wide and the overall
input cycle time during this stochastic OP computation was
$6.35\mu s$ (for the total of 6-bit pulses plus sign). This OP update was sequentially computed for all 4096 vectors in the
dataset.

After these OP updates, the value of the weights stored in
the MFMWA devices was read with two reads, one in the forward and one
in the backward (transposed) direction. We used 10000 random vectors to read and 
reconstruct the stored weight matrix in digital (by linear regression). 
In contrast to the binary
pulse pattern during the OP update, the 10000 (random)
input vectors for the MV read operations were encoded into
pulses of variable widths, where the length of the pulse was
proportional to the values of the input vectors (PWM).  We here used the same input resolution as during
the OP update (6-bit plus sign). The input vectors were
applied to the PWM of the column or row drivers, and then integrated
at the respective outputs, where the output current of the
multiply-and-add operations was integrated by charging the output
capacitors. After applying the input, the capacitor
charge was converted from analog to digital (with 9-bit ADC resolution
and in parallel for each of the 100 output lines and additionally
averaging two capacitor reads) and then further processed according to
the need of the implemented algorithms. The cycle time of this
MV operation was 8.97 $\mu s$ (in total for forward and
transposed read). It was repeated for the 10000 input vectors in the
dataset.

\subsubsection{Simulation Engine}
\label{app:ssec:sim}

%
%\alert{Malte/Tayfun:  Please write a draft of this section for the %supplement.}\\

%\textbf{Purpose}: describe our simulation engine for the analog hardware device and system of interest, including recording of operations, wall-clock time, and energy.
%\textbf{\textit{Simulator: Malte/Tayfun}}
%\begin{itemize}
%\item Key points to make:
%\begin{itemize}
%    \item Describe how we can simulate linear algebra computations %executed on analog hardware.
%    \end{itemize}
%\item Prepare a high-level write-up to be incorporated in the main %body
%\item Prepare a detailed write-up for the supplementary material
%\item Include pictures/schematics to support both write-ups, to %the extent possible
%\end{itemize}

Our simulator of the MFMWA MV and OP operations is based on the open source IBM Analog Hardware Acceleration Kit described in detail elsewhere~\cite{rasch2021flexible}.  However, for the convenience of the applications reported here, we ported the simulation of the MFMWA tiles to a custom Matlab implementation. 
The simulation is centered around an ``MFMWA tile'' that
corresponds to a 2D matrix that is stored on a non-volatile resistive
array and includes most pre- and post-processing (such as
dynamic input scaling and
%analog-digital converters)
ADCs)
that are based on the resistive processing unit (RPU) definitions
as described in detail in~\cite{gokmen2016acceleration}. Briefly, the basic operation of an MFMWA tile is
%a matrix-vector
an MV
product $y_i = \sum_j w_{ij}x_j$, 
where the matrix elements $w_{ij}$ are assumed to be stored in the
array (optionally together with the biases). The MFMWA
%matrix-vector
MV
product is, however, corrupted by noise and
non-idealities. In particular, following \cite{gokmen2016acceleration}, we restrict the input and output ranges, and assume that all the noise sources are referred to a Gaussian additive noise process at the output. Thus, the simulated MFMWA MV product can be described by   
\begin{equation}
  \label{eq:mat-vec-analog}
    y_i = f_\text{adc}\big(\sum_j w_{ij} f_\text{dac}(x_j)  + \sigma_\text{out}\xi_i\big),
\end{equation}
where $f_\text{adc}$ and $f_\text{dac}$ are discretization functions
(together with, e.g., dynamic input scaling and bound clipping) of the
digital-analog conversion process (as given by the hardware design),
and $\xi_i$ are standard normal noise processes (see \cite{gokmen2016acceleration} for details). More specifically, during an MV product, the injected additive Gaussian noise strength is assumed at 10\% of the nominal weight maximum and the MV result is bounded to 20 times the nominal weight maximum to account for signal saturation at output lines. The input signals are bounded in the range $[-1,\ldots, 1]$ with a 7-bit input resolution, whereas the outputs are quantized assuming a 9-bit ADC. To mitigate the shortcomings of the signal bounds, we use noise and update management techniques described in \cite{Gokmen2017CNN}, but no bound management as it would potentially cause runtime overhead. We also use the appropriate weight scaling (described in \cite{RaschIEEEDT2019}) when needed. Because we assume that analog quantities are converted to digital after each MV product with parallel ADCs, all other operations not performed on the MFMWA array are assumed to be in floating point precision. 

% \Malte{Vasileios: Double-check whether you used the default behavior for the update/noise/bound management as assumed in the text.}

The second operation simulated on an MFMWA tile is the OP update which is done in-memory by generating independent random pulse trains for each input and output line, where the probability of pulse occurrence codes for the two respective vector elements to add to the MFMWA matrix is as suggested and described in \cite{gokmen2016acceleration}. The incremental changes of the device conductances are assumed to be triggered by a pulse coincidence, and therefore is proportional to the product of the two vector elements, as required by the OP update.  We assume that pulse trains are finite, and each incremental update is of constant average size, but might vary in the up versus down direction and from crosspoint-to-crosspoint in a systematic fashion, in addition to a pulse-to-pulse variation as detailed in \cite{gokmen2016acceleration}. Note that this update pulse model for the update pass is the default model implemented in the IBM Hardware Acceleration Kit~\cite{rasch2021flexible}, which is called the \texttt{ConstantStep} device model.

\paragraph{Wall clock time estimation}

Since the MFMWA tile operation (OP or MV) is effectively constant in terms of the size of the matrix, we simply count the number of MFMWA operations for runtime estimation and scaled it by the expected time of an MFMWA operation (which will, e.g., scale with the time it takes to convert the input quantities to analog and back to digital and with the length of the presented stochastic pulse trains).   

%\Malte{Maybe still needs expansion, e.g. how are digital operations estimated? }

In addition, the wall-clock times take into consideration the write time required to populate the MFMWA hardware device. This takes place only once in the beginning, and its expected completion time is currently set at about ten milliseconds.

\subsubsection{Sketching of Streamed Matrices}
\label{app:ssec:sketch}

In this section we consider the computational advantages stemming  
from exploiting MFMWA arrays to reduce the runtime 
complexity and wall-clock time of matrix sketching computations 
in numerical linear algebra. Such tasks include, for example, the 
computation of a column subspace embedding of a $m\times n$ matrix 
$A$, or the solution of an
% ordinary linear least-squares 
OLLS
problem 
$\hat{x}=\argmin\limits_{x\in \mathbb{R}^n}\|Ax-b\|_2$ where 
$A\in \mathbb{R}^{m\times n}$ and $b\in \mathbb{R}^m$. Without loss 
of generality, we assume that $A$ is a dense matrix.

Our analytical results are summarized in 
Table \ref{tab:cost0-olls}
% Table~{\color{red} 2} 
of the main paper,
listing an asymptotic analysis of the number of arithmetic operations required by the above two tasks for: $(i)$ the standard approach typically used to perform the corresponding task in the published literature; $(ii)$ randomized sketching with JLT transformations, implemented on digital hardware; and $(iii)$ randomized sketching implemented on digital-analog hardware.

\subsubsection{Principal Component Analysis}

In this section we consider the computational gains achieved by applying 
PCA using the hybrid variant of randomized SVD described in
Algorithm~\ref{app:alg:algo1}.
% Algorithm~$2$.
We consider both simulated and real-world datasets.

Our analytical results are summarized in 
Table \ref{tab:cost0-pca}
% Table~{\color{red} 3} 
of the main paper, presenting an asymptotic analysis of the number of arithmetic operations of a fully digital approach and our hybrid MFMWA approach, 
% in terms of the matrix size $n$, the number of nonzero entries ${\tt nnz}(A)$, the desired number of PCs $k$ ($\sim  100$), the degree exponent $q$ ($\sim  5$) and the update rank $p$, 
when applied toward computational tasks associated with eigenvalue decomposition.

We next consider a real-world data application of our 
Algorithm~\ref{app:alg:algo1}
% Algorithm~$2$
in the context of \emph{population genetics}, namely
%We start by considering
the tasks of single nucleotide polymorphisms (SNPs) tagging and eigengene analysis in deoxyribonucleic acid (DNA) microarray gene analysis~\cite{paschou2007intra}. 
% The objective is to select a subset of gene expressions or SNPs from a table of SNPs for different populations that captures the spectral information (variations) of the population of
%interest~\cite{paschou2007intra}. 
In particular, for the Genetics application, we studied two different genetics datasets:
the ``Yale'' dataset containing a total of $248$ SNPs of around $2000$ unrelated individuals 
from $38$ populations around the world, where we consider four genomic regions (SORCS3, 
PAH, HOXB, 17q25); and the ``Hapmap'' dataset representing a public database from the HapMap project2 (phase I) that comprises $1,000,000$ SNPs from different populations, where we consider the same four regions following the approach
%in~\cite{paschou2007intra,ubaru2019sampling} 
in~\cite{paschou2007intra}
to encode the matrix $A$ for each dataset. 
Table~\ref{tab:genetics}
%Figure~\ref{fig:test} (right)
% %of the supplement
% (supplement)
presents the corresponding errors $\|A- \widehat{U}_k\widehat{U}_k^{\top}A\|/\|A\|$ from the fully digital randomized approach and our hybrid analog-digital PCA approach
% %(as discussed in Section~\ref{ssc:randPCA})
% (Section~\ref{ssc:randPCA})
(as discussed in the main body of the paper)
averaged over ten trials for each instance of the genetic datasets (and corresponding matrix $A$).
%, see details in the supplement.
The number of PCs was fixed to $k=5$, and we chose the sketch size $\ell=k$ for the fully digital approach.
Our hybrid approach is considered under the parameter settings $\ell\in\{k,2k,3k\}$. 
% Note that, with the use of analog devices in our hybrid approach, computing a larger sketch (i.e., increasing the size $\ell$) does not significantly increase the runtime; in particular, the sketching cost remains the same irrespective of the sketch size, but there will be a small additional cost for orthogonalization.
The magnitude of the error results shows how well the eigengenes (i.e., PCs) capture the variance of the underlying dataset.
We observe that the significant runtime improvements from our hybrid randomized 
PCA approach demonstrated above can be realized with an error that is comparable to the error of the fully digital approach (particularly for $\ell=3k$), even performing better than the digital approach in several cases.

\begin{table}[tb]
 \caption{\it Genetics Data Results for the Fully Digital Randomized Approach and Our Hybrid Analog-Digital Approach. Average error over ten trials for each instance of  datasets.}\label{tab:genetics}
%\resizebox{1\textwidth}{!}{
 \begin{center} 
 \begin{tabular}{|l|c|c|c|c|c|}
\hline
Genetics Dataset &
Matrix Size &
Digital PCA &
\multicolumn{3}{c|}{Hybrid PCA}\\\cline{3-6}
&&$\ell = k=5$ &
 $\ell= k =5$&
$\ell= 2k $&
 $\ell=3k$\\
\hline
\hline
Yaledataset/SORCS3&
$1966 \times 53$&
0.1247 &
0.2156 &
 0.1540&
0.1258 \\
Yaledataset/PAH&
$1979 \times 32$&
0.1175&
0.1558&
0.1231&
 {\bf 0.1095}\\
Yaledataset/HOXB&
$1953 \times 96$&
0.2424&
0.3117&
 0.2775&
0.2453\\
Yaledataset/17q25&
$1962 \times  63$&
0.3186&
0.3764&
0.3362&
{\bf 0.2627}\\
\hline
HapMap/SORCS3&
$268 \times  307$&
0.1268&
0.2982&
0.2309&
0.1929\\
HapMap/PAH&
$266 \times 88$&
0.1251&
0.2856&
0.1685&
{\bf 0.1202}\\
HapMap/HOXB&
$269 \times 571$&
0.2750&
0.3815&
0.3310&
0.3082\\
HapMap/17q25&
$265 \times 370$&
0.3143&
0.4043&
0.3542&
0.3363\\
\hline
\end{tabular}
 \end{center} 
 %}
\end{table}

\bibliographystyle{IEEEtran}
\bibliography{local}

\begin{thebibliography}{10}
\providecommand{\url}[1]{#1}
\csname url@samestyle\endcsname
\providecommand{\newblock}{\relax}
\providecommand{\bibinfo}[2]{#2}
\providecommand{\BIBentrySTDinterwordspacing}{\spaceskip=0pt\relax}
\providecommand{\BIBentryALTinterwordstretchfactor}{4}
\providecommand{\BIBentryALTinterwordspacing}{\spaceskip=\fontdimen2\font plus
\BIBentryALTinterwordstretchfactor\fontdimen3\font minus
  \fontdimen4\font\relax}
\providecommand{\BIBforeignlanguage}[2]{{%
\expandafter\ifx\csname l@#1\endcsname\relax
\typeout{** WARNING: IEEEtran.bst: No hyphenation pattern has been}%
\typeout{** loaded for the language `#1'. Using the pattern for}%
\typeout{** the default language instead.}%
\else
\language=\csname l@#1\endcsname
\fi
#2}}
\providecommand{\BIBdecl}{\relax}
\BIBdecl

\bibitem{Sustainability2040}
P.~Dockrill, ``Computers will require more energy than the world generates by
  2040,'' \emph{ScienceAlert}, July 2016.

\bibitem{online}
``{AI} can do great things - if it doesn't burn the planet,''
  \url{https://www.wired.com/story/ai-great-things-burn-planet/}, accessed:
  2021-06-09.

\bibitem{https://doi.org/10.1002/aisy.202000085}
\BIBentryALTinterwordspacing
A.~Mehonic, A.~Sebastian, B.~Rajendran, O.~Simeone, E.~Vasilaki, and A.~J.
  Kenyon, ``Memristors - from in-memory computing, deep learning acceleration,
  and spiking neural networks to the future of neuromorphic and bio-inspired
  computing,'' \emph{Advanced Intelligent Systems}, vol.~2, no.~11, p. 2000085,
  2020. [Online]. Available:
  \url{https://onlinelibrary.wiley.com/doi/abs/10.1002/aisy.202000085}
\BIBentrySTDinterwordspacing

\bibitem{mooreslaw}
J.~Shalf, ``The future of computing beyond {M}oore's law,'' \emph{Phil. Trans.
  R. Soc. A.}, vol. 378, no. 20190061, 2020.

\bibitem{dongarra2017design}
J.~Dongarra, S.~Hammarling, N.~J. Higham, S.~D. Relton, P.~Valero-Lara, and
  M.~Zounon, ``The design and performance of batched blas on modern
  high-performance computing systems,'' \emph{Procedia Computer Science}, vol.
  108, pp. 495--504, 2017.

\bibitem{yang2013memristive}
J.~J. Yang, D.~B. Strukov, and D.~R. Stewart, ``Memristive devices for
  computing,'' \emph{Nature nanotechnology}, vol.~8, no.~1, pp. 13--24, 2013.

\bibitem{agarwal2016energy}
S.~Agarwal, T.-T. Quach, O.~Parekh, A.~H. Hsia, E.~P. DeBenedictis, C.~D.
  James, M.~J. Marinella, and J.~B. Aimone, ``Energy scaling advantages of
  resistive memory crossbar based computation and its application to sparse
  coding,'' \emph{Frontiers in neuroscience}, vol.~9, p. 484, 2016.

\bibitem{chua:1971}
L.~Chua, ``Memristor-the missing circuit element,'' \emph{IEEE Transactions on
  Circuit Theory}, vol.~18, no.~5, pp. 507--519, 1971.

\bibitem{ambrogio:analog:2023}
S.~Ambrogio, P.~Narayanan, A.~Okazaki, A.~Fasoli, C.~Mackin, K.~Hosokawa,
  A.~Nomura, T.~Yasuda, A.~Chen, A.~Friz, M.~Ishii, J.~Luquin, Y.~Kohda,
  N.~Saulnier, K.~Brew, S.~Choi, I.~Ok, T.~Philip, V.~Chan, C.~Silvestre,
  I.~Ahsan, V.~Narayanan, H.~Tsai, and G.~W. Burr, ``An analog-{AI} chip for
  energy-efficient speech recognition and transcription,'' \emph{Nature}, vol.
  620, pp. 768--775, 2023.

\bibitem{sebastian2020memory}
A.~Sebastian, M.~Le~Gallo, R.~Khaddam-Aljameh, and E.~Eleftheriou, ``Memory
  devices and applications for in-memory computing,'' \emph{Nature
  nanotechnology}, vol.~15, no.~7, pp. 529--544, 2020.

\bibitem{analogjackson60}
A.~S. Jackson, \emph{Analog Computation}.\hskip 1em plus 0.5em minus
  0.4em\relax McGraw-Hill, 1960.

\bibitem{ambrogio2018equivalent}
S.~Ambrogio, P.~Narayanan, H.~Tsai, R.~M. Shelby, I.~Boybat, C.~Di~Nolfo,
  S.~Sidler, M.~Giordano, M.~Bodini, N.~C. Farinha \emph{et~al.},
  ``Equivalent-accuracy accelerated neural-network training using analogue
  memory,'' \emph{Nature}, vol. 558, no. 7708, pp. 60--67, 2018.

\bibitem{feinberg2018enabling}
B.~Feinberg, U.~K.~R. Vengalam, N.~Whitehair, S.~Wang, and E.~Ipek, ``Enabling
  scientific computing on memristive accelerators,'' in \emph{2018 ACM/IEEE
  45th Annual International Symposium on Computer Architecture (ISCA)}.\hskip
  1em plus 0.5em minus 0.4em\relax IEEE, 2018, pp. 367--382.

\bibitem{golub2013matrix}
G.~H. Golub and C.~F. Van~Loan, \emph{Matrix computations}.\hskip 1em plus
  0.5em minus 0.4em\relax JHU press, 2013.

\bibitem{gittens2016matrix}
A.~Gittens, A.~Devarakonda, E.~Racah, M.~Ringenburg, L.~Gerhardt, J.~Kottalam,
  J.~Liu, K.~Maschhoff, S.~Canon, J.~Chhugani \emph{et~al.}, ``Matrix
  factorizations at scale: A comparison of scientific data analytics in spark
  and {C+MPI} using three case studies,'' in \emph{2016 IEEE International
  Conference on Big Data (Big Data)}.\hskip 1em plus 0.5em minus 0.4em\relax
  IEEE, 2016, pp. 204--213.

\bibitem{LiHuGuTi04}
L.~Li, W.~Huang, I.~Y.-H. Gu, , and Q.~Tian, ``Statistical modeling of complex
  backgrounds for foreground object detection,'' \emph{IEEE Transactions on
  Image Processing}, vol.~13, pp. 1459--1472, 2004.

\bibitem{candes2011robust}
E.~J. Cand{\`e}s, X.~Li, Y.~Ma, and J.~Wright, ``Robust principal component
  analysis?'' \emph{Journal of the ACM (JACM)}, vol.~58, no.~3, pp. 1--37,
  2011.

\bibitem{bush2012genome}
W.~S. Bush and J.~H. Moore, ``Genome-wide association studies,'' \emph{PLoS
  Comput Biol}, vol.~8, no.~12, p. e1002822, 2012.

\bibitem{Uffelmann2021}
E.~Uffelmann, Q.~Q. Huang, N.~S. Munung, J.~de~Vries, Y.~Okada, A.~R. Martin,
  H.~C. Martin, T.~Lappalainen, and D.~Posthuma, ``Genome-wide association
  studies,'' \emph{Nature Reviews}, pp. 1--21, 2021.

\bibitem{Mbatchou2021}
J.~Mbatchou, L.~Barnard, J.~Backman, A.~Marcketta, J.~A. Kosmicki,
  A.~Ziyatdinov, C.~Benner, C.~O'Dushlaine, M.~Barber, B.~Boutkov, L.~Habegger,
  M.~Ferreira, A.~Baras, J.~Reid, G.~Abecasis, E.~Maxwell, and J.~Marchini,
  ``Computationally efficient whole-genome regression for quantitative and
  binary traits,'' \emph{Nature Genetics}, vol.~53, pp. 1097--1103, 2021.

\bibitem{paschou2007pca}
P.~Paschou, E.~Ziv, E.~G. Burchard, S.~Choudhry, W.~Rodriguez-Cintron, M.~W.
  Mahoney, and P.~Drineas, ``{PCA}-correlated {SNP}s for structure
  identification in worldwide human populations,'' \emph{PLoS Genet}, vol.~3,
  no.~9, p. e160, 2007.

\bibitem{bose2019terapca}
A.~Bose, V.~Kalantzis, E.-M. Kontopoulou, M.~Elkady, P.~Paschou, and
  P.~Drineas, ``Tera{PCA}: a fast and scalable software package to study
  genetic variation in tera-scale genotypes,'' \emph{Bioinformatics}, vol.~35,
  no.~19, pp. 3679--3683, 2019.

\bibitem{drineas2016randnla}
P.~Drineas and M.~W. Mahoney, ``Randnla: randomized numerical linear algebra,''
  \emph{Communications of the ACM}, vol.~59, no.~6, pp. 80--90, 2016.

\bibitem{kockan2020sketching}
C.~Kockan, K.~Zhu, N.~Dokmai, N.~Karpov, M.~O. Kulekci, D.~P. Woodruff, and
  S.~C. Sahinalp, ``Sketching algorithms for genomic data analysis and querying
  in a secure enclave,'' \emph{Nature Methods}, vol.~17, no.~3, pp. 295--301,
  2020.

\bibitem{Kalantzis2021-nb}
V.~Kalantzis, A.~Gupta, L.~Horesh, T.~Nowicki, M.~S. Squillante, and C.~W. Wu,
  ``Solving sparse linear systems with approximate inverse preconditioners on
  analog devices,'' \url{https://arxiv.org/abs/2107.06973}, Jul. 2021.

\bibitem{woodruff2014sketching}
D.~P. Woodruff, ``Sketching as a tool for numerical linear algebra,''
  \emph{Theoretical Computer Science}, vol.~10, no. 1-2, pp. 1--157, 2014.

\bibitem{sarlos2006improved}
T.~Sarlos, ``Improved approximation algorithms for large matrices via random
  projections,'' in \emph{2006 47th Annual IEEE Symposium on Foundations of
  Computer Science (FOCS'06)}.\hskip 1em plus 0.5em minus 0.4em\relax IEEE,
  2006, pp. 143--152.

\bibitem{ubaru2017low}
S.~Ubaru, A.~Mazumdar, and Y.~Saad, ``Low rank approximation and decomposition
  of large matrices using error correcting codes,'' \emph{IEEE Transactions on
  Information Theory}, vol.~63, no.~9, pp. 5544--5558, 2017.

\bibitem{avron2010blendenpik}
H.~Avron, P.~Maymounkov, and S.~Toledo, ``Blendenpik: Supercharging lapack's
  least-squares solver,'' \emph{SIAM Journal on Scientific Computing}, vol.~32,
  no.~3, pp. 1217--1236, 2010.

\bibitem{pearson1901liii}
K.~Pearson, ``{LIII}. on lines and planes of closest fit to systems of points
  in space,'' \emph{The London, Edinburgh, and Dublin Philosophical Magazine
  and Journal of Science}, vol.~2, no.~11, pp. 559--572, 1901.

\bibitem{hotelling1933analysis}
H.~Hotelling, ``Analysis of a complex of statistical variables into principal
  components.'' \emph{Journal of educational psychology}, vol.~24, no.~6, p.
  417, 1933.

\bibitem{gu2015subspace}
M.~Gu, ``Subspace iteration randomization and singular value problems,''
  \emph{SIAM Journal on Scientific Computing}, vol.~37, no.~3, pp.
  A1139--A1173, 2015.

\bibitem{halko2011finding}
N.~Halko, P.-G. Martinsson, and J.~A. Tropp, ``Finding structure with
  randomness: Probabilistic algorithms for constructing approximate matrix
  decompositions,'' \emph{SIAM review}, vol.~53, no.~2, pp. 217--288, 2011.

\bibitem{kohda2020IEDM}
Y.~Kohda, Y.~Li, K.~Hosokawa, S.~Kim, R.~Khaddam-Aljameh, Z.~Ren, P.~Solomon,
  T.~Gokmen, S.~Rajalingam, C.~Baks, W.~Haensch, and E.~Leobandung,
  ``Unassisted true analog neural network training chip,'' in \emph{2020 IEEE
  International Electron Devices Meeting (IEDM)}, 2020, pp. 36.2.1--36.2.4.

\bibitem{lecun2015deep}
Y.~LeCun, Y.~Bengio, and G.~Hinton, ``Deep learning,'' \emph{nature}, vol. 521,
  no. 7553, pp. 436--444, 2015.

\bibitem{kuleshov2019machine}
V.~Kuleshov, J.~Ding, C.~Vo, B.~Hancock, A.~Ratner, Y.~Li, C.~R{\'e},
  S.~Batzoglou, and M.~Snyder, ``A machine-compiled database of genome-wide
  association studies,'' \emph{Nature communications}, vol.~10, no.~1, pp.
  1--8, 2019.

\bibitem{turk1991eigenfaces}
M.~Turk and A.~Pentland, ``Eigenfaces for recognition,'' \emph{Journal of
  cognitive neuroscience}, vol.~3, no.~1, pp. 71--86, 1991.

\bibitem{bouwmans2016handbook}
T.~Bouwmans, N.~S. Aybat, and E.-h. Zahzah, \emph{Handbook of robust low-rank
  and sparse matrix decomposition: Applications in image and video
  processing}.\hskip 1em plus 0.5em minus 0.4em\relax CRC Press, 2016.

\bibitem{armstrong2000geography}
M.~P. Armstrong, \emph{Geography and computational science}.\hskip 1em plus
  0.5em minus 0.4em\relax Taylor \& Francis, 2000.

\bibitem{yang2015implementing}
J.~Yang, X.~Meng, and M.~W. Mahoney, ``Implementing randomized matrix
  algorithms in parallel and distributed environments,'' \emph{Proceedings of
  the IEEE}, vol. 104, no.~1, pp. 58--92, 2015.

\bibitem{almgren2013nyx}
A.~S. Almgren, J.~B. Bell, M.~J. Lijewski, Z.~Luki{\'c}, and E.~Van~Andel,
  ``Nyx: A massively parallel amr code for computational cosmology,'' \emph{The
  Astrophysical Journal}, vol. 765, no.~1, p.~39, 2013.

\bibitem{bedorf201424}
J.~B{\'e}dorf, E.~Gaburov, M.~S. Fujii, K.~Nitadori, T.~Ishiyama, and S.~P.
  Zwart, ``24.77 pflops on a gravitational tree-code to simulate the milky way
  galaxy with 18600 gpus,'' in \emph{SC'14: Proceedings of the International
  Conference for High Performance Computing, Networking, Storage and
  Analysis}.\hskip 1em plus 0.5em minus 0.4em\relax IEEE, 2014, pp. 54--65.

\bibitem{werner2005monitoring}
G.~Werner-Allen, J.~Johnson, M.~Ruiz, J.~Lees, and M.~Welsh, ``Monitoring
  volcanic eruptions with a wireless sensor network,'' in \emph{Proceeedings of
  the Second European Workshop on Wireless Sensor Networks, 2005.}\hskip 1em
  plus 0.5em minus 0.4em\relax IEEE, 2005, pp. 108--120.

\bibitem{giorgi1991approaches}
F.~Giorgi and L.~O. Mearns, ``Approaches to the simulation of regional climate
  change: a review,'' \emph{Reviews of Geophysics}, vol.~29, no.~2, pp.
  191--216, 1991.

\bibitem{comon1990tracking}
P.~Comon and G.~H. Golub, ``Tracking a few extreme singular values and vectors
  in signal processing,'' \emph{Proceedings of the IEEE}, vol.~78, no.~8, pp.
  1327--1343, 1990.

\bibitem{yang1995projection}
B.~Yang, ``Projection approximation subspace tracking,'' \emph{IEEE
  Transactions on Signal processing}, vol.~43, no.~1, pp. 95--107, 1995.

\bibitem{berry1995using}
M.~W. Berry, S.~T. Dumais, and G.~W. O'Brien, ``Using linear algebra for
  intelligent information retrieval,'' \emph{SIAM Review}, vol.~37, no.~4, pp.
  573--595, 1995.

\bibitem{vaswani2017attention}
A.~Vaswani, N.~Shazeer, N.~Parmar, J.~Uszkoreit, L.~Jones, A.~N. Gomez,
  {\L}.~Kaiser, and I.~Polosukhin, ``Attention is all you need,'' in
  \emph{Proceedings of the 31st International Conference on Neural Information
  Processing Systems}, 2017, pp. 6000--6010.

\bibitem{kendall2000high}
R.~A. Kendall, E.~Apr{\`a}, D.~E. Bernholdt, E.~J. Bylaska, M.~Dupuis, G.~I.
  Fann, R.~J. Harrison, J.~Ju, J.~A. Nichols, J.~Nieplocha \emph{et~al.},
  ``High performance computational chemistry: An overview of nwchem a
  distributed parallel application,'' \emph{Computer Physics Communications},
  vol. 128, no. 1-2, pp. 260--283, 2000.

\bibitem{sliwoski2014computational}
G.~Sliwoski, S.~Kothiwale, J.~Meiler, and E.~W. Lowe, ``Computational methods
  in drug discovery,'' \emph{Pharmacological reviews}, vol.~66, no.~1, pp.
  334--395, 2014.

\bibitem{reich2008principal}
D.~Reich, A.~L. Price, and N.~Patterson, ``Principal component analysis of
  genetic data,'' \emph{Nature genetics}, vol.~40, no.~5, pp. 491--492, 2008.

\bibitem{ubaru2019find}
S.~Ubaru, A.-K. Seghouane, and Y.~Saad, ``Find the dimension that counts: Fast
  dimension estimation and krylov {PCA},'' in \emph{Proceedings of the 2019
  SIAM International Conference on Data Mining}.\hskip 1em plus 0.5em minus
  0.4em\relax SIAM, 2019, pp. 720--728.

\bibitem{gokmen2016acceleration}
T.~Gokmen and Y.~Vlasov, ``Acceleration of deep neural network training with
  resistive cross-point devices: Design considerations,'' \emph{Frontiers in
  neuroscience}, vol.~10, p. 333, 2016.

\bibitem{burr2015experimental}
G.~W. Burr, R.~M. Shelby, S.~Sidler, C.~Di~Nolfo, J.~Jang, I.~Boybat, R.~S.
  Shenoy, P.~Narayanan, K.~Virwani, E.~U. Giacometti \emph{et~al.},
  ``Experimental demonstration and tolerancing of a large-scale neural network
  (165 000 synapses) using phase-change memory as the synaptic weight
  element,'' \emph{IEEE Transactions on Electron Devices}, vol.~62, no.~11, pp.
  3498--3507, 2015.

\bibitem{burr2017neuromorphic}
G.~W. Burr, R.~M. Shelby, A.~Sebastian, S.~Kim, S.~Kim, S.~Sidler, K.~Virwani,
  M.~Ishii, P.~Narayanan, A.~Fumarola \emph{et~al.}, ``Neuromorphic computing
  using non-volatile memory,'' \emph{Advances in Physics: X}, vol.~2, no.~1,
  pp. 89--124, 2017.

\bibitem{cite-key}
\BIBentryALTinterwordspacing
K.~Steinbuch, ``Die lernmatrix,'' \emph{Kybernetik}, vol.~1, no.~1, pp. 36--45,
  1961. [Online]. Available: \url{https://doi.org/10.1007/BF00293853}
\BIBentrySTDinterwordspacing

\bibitem{gokmen2020algorithm}
T.~Gokmen and W.~Haensch, \emph{Frontiers in neuroscience}, vol.~14, 2020.

\bibitem{IBMJRD2019}
H.-Y. Chang, P.~Narayanan, S.~C. Lewis, N.~C.~P. Farinha, K.~Hosokawa,
  C.~Mackin, H.~Tsai, S.~Ambrogio, A.~Chen, and G.~W. Burr, ``{AI} hardware
  acceleration with analog memory: Microarchitectures for low energy at high
  speed,'' \emph{IBM Journal of Research and Development}, vol.~63, no.~6, pp.
  8:1--8:14, 2019.

\bibitem{marinella2018}
M.~J. Marinella, S.~Agarwal, A.~Hsia, I.~Richter, R.~Jacobs-Gedrim, J.~Niroula,
  S.~J. Plimpton, E.~Ipek, and C.~D. James, ``Multiscale co-design analysis of
  energy, latency, area, and accuracy of a reram analog neural training
  accelerator,'' \emph{IEEE Journal on Emerging and Selected Topics in Circuits
  and Systems}, vol.~8, no.~1, pp. 86--101, 2018.

\bibitem{bfloat16}
``The bfloat16 numerical format,''
  \url{https://cloud.google.com/tpu/docs/bfloat16}, accessed: 2021-12-07.

\bibitem{NvidiaA100}
``{NVIDIA A100} tensor core {GPU},''
  \url{https://www.nvidia.com/en-us/data-center/a100/}, accessed: 2021-11-02.

\bibitem{graphcore}
``Graphcore,'' \url{https://www.graphcore.ai/}, accessed: 2021-12-07.

\bibitem{Sun2019}
\BIBentryALTinterwordspacing
X.~Sun, J.~Choi, C.~Chen, N.~Wang, S.~Venkataramani, V.~Srinivasan, X.~Cui,
  W.~Zhang, and K.~Gopalakrishnan, ``Hybrid 8-bit floating point {(HFP8)}
  training and inference for deep neural networks,'' in \emph{Advances in
  Neural Information Processing Systems 32: Annual Conference on Neural
  Information Processing Systems 2019, NeurIPS 2019, December 8-14, 2019,
  Vancouver, BC, Canada}, H.~M. Wallach, H.~Larochelle, A.~Beygelzimer,
  F.~d'Alch{\'{e}}{-}Buc, E.~B. Fox, and R.~Garnett, Eds., 2019, pp.
  4901--4910. [Online]. Available:
  \url{https://proceedings.neurips.cc/paper/2019/hash/65fc9fb4897a89789352e211ca2d398f-Abstract.html}
\BIBentrySTDinterwordspacing

\bibitem{Sun2019a}
\BIBentryALTinterwordspacing
------, ``Hybrid 8-bit floating point {(HFP8)} training and inference for deep
  neural networks,'' in \emph{Advances in Neural Information Processing Systems
  32: Annual Conference on Neural Information Processing Systems 2019, NeurIPS
  2019, December 8-14, 2019, Vancouver, BC, Canada}, H.~M. Wallach,
  H.~Larochelle, A.~Beygelzimer, F.~d'Alch{\'{e}}{-}Buc, E.~B. Fox, and
  R.~Garnett, Eds., 2019, pp. 4901--4910. [Online]. Available:
  \url{https://proceedings.neurips.cc/paper/2019/hash/65fc9fb4897a89789352e211ca2d398f-Abstract.html}
\BIBentrySTDinterwordspacing

\bibitem{Choi2019}
\BIBentryALTinterwordspacing
J.~Choi, S.~Venkataramani, V.~Srinivasan, K.~Gopalakrishnan, Z.~Wang, and
  P.~Chuang, ``Accurate and efficient 2-bit quantized neural networks,'' in
  \emph{Proceedings of Machine Learning and Systems 2019, MLSys 2019, Stanford,
  CA, USA, March 31 - April 2, 2019}, A.~Talwalkar, V.~Smith, and M.~Zaharia,
  Eds.\hskip 1em plus 0.5em minus 0.4em\relax mlsys.org, 2019. [Online].
  Available: \url{https://proceedings.mlsys.org/book/268.pdf}
\BIBentrySTDinterwordspacing

\bibitem{Qin2020}
H.~Qin, R.~Gong, X.~Liu, X.~Bai, J.~Song, and N.~Sebe, ``Binary neural
  networks: {A} survey,'' \emph{Pattern Recognition}, vol. 105, p. 107281,
  2020.

\bibitem{Alaghi2013}
A.~Alaghi and J.~P. Hayes, ``Survey of stochastic computing,'' \emph{ACM
  Transactions on Embedded Computing Systems}, vol.~12, no.~2s, pp. 1--19,
  2013.

\bibitem{dithercomputing2021}
C.~W. Wu, ``Dither computing: a hybrid deterministic-stochastic computing
  framework,'' in \emph{Proceedings of the 28th IEEE Symposium on Computer
  Arithmetic}, 2021.

\bibitem{connolly2021}
M.~P. Connolly, N.~J. Higham, and T.~Mary, ``Stochastic rounding and its
  probabilistic backward error analysis,'' \emph{SIAM Journal on Scientific
  Computing}, vol.~43, no.~1, pp. A566--A585, 2021.

\bibitem{wang:2018}
N.~Wang, J.~Choi, D.~Brand, C.-Y. Chen, and K.~Gopalakrishnan, ``Training deep
  neural networks with 8-bit floating point numbers,'' in \emph{Proceedings of
  the 31st International Conference on Neural Information Processing
  Systems(NIPS 2018)}, 2018, pp. 7675--7684.

\bibitem{Gupta2015a}
S.~Gupta, A.~Agrawal, K.~Gopalakrishnan, and P.~Narayanan, ``Deep learning with
  limited numerical precision,'' in \emph{Proceedings of the 32nd International
  Conference on Machine Learning}, Feb. 2015, pp. 1737--1746.

\bibitem{motwani_raghavan_1995}
R.~Motwani and P.~Raghavan, \emph{Randomized Algorithms}.\hskip 1em plus 0.5em
  minus 0.4em\relax Cambridge University Press, 1995.

\bibitem{le2018mixed}
M.~Le~Gallo, A.~Sebastian, R.~Mathis, M.~Manica, H.~Giefers, T.~Tuma, C.~Bekas,
  A.~Curioni, and E.~Eleftheriou, ``Mixed-precision in-memory computing,''
  \emph{Nature Electronics}, vol.~1, no.~4, pp. 246--253, 2018.

\bibitem{rasch2021flexible}
M.~J. Rasch, D.~Moreda, T.~Gokmen, M.~Le~Gallo, F.~Carta, C.~Goldberg,
  K.~El~Maghraoui, A.~Sebastian, and V.~Narayanan, ``A flexible and fast
  pytorch toolkit for simulating training and inference on analog crossbar
  arrays,'' in \emph{2021 IEEE 3rd International Conference on Artificial
  Intelligence Circuits and Systems (AICAS)}.\hskip 1em plus 0.5em minus
  0.4em\relax IEEE, 2021, pp. 1--4.

\bibitem{Gokmen2017CNN}
\BIBentryALTinterwordspacing
T.~Gokmen, M.~Onen, and W.~Haensch, ``Training deep convolutional neural
  networks with resistive cross-point devices,'' \emph{Frontiers in
  Neuroscience}, vol.~11, p. 538, 2017. [Online]. Available:
  \url{https://www.frontiersin.org/article/10.3389/fnins.2017.00538}
\BIBentrySTDinterwordspacing

\bibitem{RaschIEEEDT2019}
M.~J. Rasch, T.~Gokmen, and W.~Haensch, ``Training large-scale artificial
  neural networks on simulated resistive crossbar arrays,'' \emph{IEEE Design
  Test}, vol.~37, no.~2, pp. 19--29, 2020.

\bibitem{paschou2007intra}
P.~Paschou, M.~W. Mahoney, A.~Javed, J.~R. Kidd, A.~J. Pakstis, S.~Gu, K.~K.
  Kidd, and P.~Drineas, ``Intra-and interpopulation genotype reconstruction
  from tagging snps,'' \emph{Genome Research}, vol.~17, no.~1, pp. 96--107,
  2007.

\end{thebibliography}

\end{document}